\documentclass[11pt]{article}

\usepackage{color}
\usepackage{latexsym}
\usepackage{amssymb}
\usepackage{amsmath, amsfonts,amssymb,theorem,euscript,array,enumerate,amsfonts,mathrsfs}
\usepackage{graphicx}
\usepackage[colorlinks=true,linkcolor=blue,urlcolor=blue]{hyperref}

\newtheorem{Theorem}{Theorem}[section]
\newtheorem{Definition}{Definition}[section]
\newtheorem{Proposition}{Proposition}[section]

\newtheorem{Lemma}{Lemma}[section]

\newtheorem{Remark}{Remark}[section]

\def\esssup_#1{\underset{#1}{\mathrm{ess\,sup\, }}}
\def\essinf_#1{\underset{#1}{\mathrm{ess\,inf\, }}}

\def \R{\mathbb{R}}
\def \E{\mathbb{E}}
\def \F{\mathbb{F}}
\def \G{\mathbb{G}}
\def \P{\mathbb{P}}
\def \Q{\mathbb{Q}}

\def \Ac{{\cal A}}
\def \Bc{{\cal B}}
\def \Dc{{\cal D}}
\def \Ec{{\cal E}}
\def \Fc{{\cal F}}
\def \Gc{{\cal G}}
\def \Pc{{\cal P}}
\def \Mc{{\cal M}}
\def \Vc{{\cal V}}

\def \eps{\varepsilon}

\def \x{\boldsymbol x}

\def \ep{\hbox{ }\hfill$\Box$}
\def \epR{\hbox{ }\hfill$\lozenge$}

\allowdisplaybreaks

\begin{document}

\title{BSDE Representation and Randomized Dynamic Programming Principle for Stochastic Control Problems\\ of Infinite-Dimensional Jump-Diffusions}

\author{
Elena BANDINI\footnote{Department of Mathematics, University of Milano-Bicocca, Via Roberto Cozzi 55, 20125 Milano, Italy, \sf elena.bandini at unimib.it}
\qquad\quad
Fulvia CONFORTOLA\footnote{Department of Mathematics, Politecnico di Milano, Via Bonardi 9, 20133 Milano, Italy, \sf fulvia.confortola at polimi.it}
\qquad\quad
Andrea COSSO\footnote{Department of Mathematics, University of Bologna, Piazza di Porta S. Donato 5, 40126 Bologna, Italy, \sf andrea.cosso at unibo.it}
}

\maketitle

\begin{abstract}
We consider a general class of stochastic optimal control problems, where the state process lives in a real separable Hilbert space and is driven by a cylindrical Brownian motion and a Poisson random measure; no special structure is imposed on the coefficients, which are also allowed to be path-dependent; in addition, the diffusion coefficient can be degenerate. For such a class of stochastic control problems, we prove, by means of purely probabilistic techniques based on the so-called randomization method, that the value of the control problem admits a probabilistic representation formula (known as non-linear Feynman-Kac formula) in terms of a suitable backward stochastic differential equation. This probabilistic representation considerably extends current results in the literature on the infinite-dimensional case, and it is also relevant in finite dimension. Such a representation allows to show, in the non-path-dependent (or Markovian) case, that the value function satisfies the so-called randomized dynamic programming principle. As a consequence, we are able to prove that the value function is a viscosity solution of the corresponding Hamilton-Jacobi-Bellman equation, which turns out to be a second-order fully non-linear integro-differential equation in Hilbert space.
\end{abstract}

\vspace{7mm}

\noindent {\bf Keywords:} Backward stochastic differential equations, infinite-dimensional path-dependent controlled SDEs, randomization method, viscosity solutions.

\vspace{7mm}

\noindent {\bf 2010 Mathematics Subject Classification:} 60H10, 60H15, 93E20, 49L25.

\section{Introduction}

\setcounter{equation}{0} \setcounter{Assumption}{0}
\setcounter{Theorem}{0} \setcounter{Proposition}{0}
\setcounter{Corollary}{0} \setcounter{Lemma}{0}
\setcounter{Definition}{0} \setcounter{Remark}{0}

In the present paper we study a general class of stochastic optimal control problems, where the infinite-dimensional state process, taking values in a real separable Hilbert space $H$, has a dynamics driven by a cylindrical Brownian motion $W$ and a Poisson random measure $\pi$. Moreover, the coefficients are assumed to be path-dependent, in the sense that they depend on the past trajectory of the state process. In addition, the space of control actions $\Lambda$ can be any Borel space (i.e., any topological space homeomorphic to a Borel subset of a Polish space).
More precisely, the controlled state process is a so-called mild solution to the following equation:
\[
\begin{cases}
\vspace{1mm}\displaystyle dX_t = AX_t\,dt + b_t(X,\alpha_t)dt + \sigma_t(X,\alpha_t)dW_t+\!\! \int_{U \setminus \{0\}}\!\!\!\!\gamma_t(X, \alpha_t, z)\big(\pi(dt\,dz)-\lambda_\pi(dz)\,dt\big), \;\;0\leq t\leq T, \\
\displaystyle X_0 = x_0,
\end{cases}
\]
where $A$ is a linear operator  generating a strongly continuous semigroup $\{e^{tA},t\geq0\}$, $\lambda_\pi(dz)dt$ is the compensator of $\pi$, while $\alpha$ is an admissible control process, that is a predictable stochastic process taking values in $\Lambda$. Given an admissible control $\alpha$, the corresponding gain functional is given by
\[
J(\alpha) \ = \ \E\bigg[\int_0^T f_t(X^{x_0,\alpha},\alpha_t)\,dt + g(X^{x_0,\alpha})\bigg],
\]
where the running  and terminal reward  functionals $f$ and $g$ may also depend on the past trajectory of the state process. The value of the stochastic control problem, starting at $t=0$ from $x_0$, is defined as
\begin{equation}\label{Intro_valufunc}
V_0 \ = \ \sup_\alpha J(\alpha).
\end{equation}

Stochastic optimal control problems of infinite-dimensional processes have been extensively studied using the theory of Backward Stochastic Differential Equations (BSDEs); we mention in particular the seminal papers \cite{FuhrmanTessitore02}, \cite{FuhrmanTessitore04} and the last chapter of the recent book \cite{fabbrigozziswiech17}, where a detailed discussion of the literature can be found. Notice however that the current results require a special structure of the controlled state equations, namely that the diffusion coefficient $\sigma=\sigma(t,x)$ is uncontrolled and the drift has the following specific form $b=b_1(t,x)+\sigma(t,x)b_2(t,x,a)$.  Up to our knowledge, only the recent paper \cite{CoGuaTess}, which is devoted to the study of ergodic control problems, applies the BSDEs techniques to a more general class of infinite-dimensional controlled state processes; in \cite{CoGuaTess} the drift has the general form $b=b(x,a)$, however the diffusion coefficient is still uncontrolled and indeed constant, moreover the space of control actions $\Lambda$ is assumed to be a real separable Hilbert space (or, more generally, according to Remark 2.2 in \cite{CoGuaTess}, $\Lambda$ has to be the image of a continuous surjection $\varphi$ defined on some real separable Hilbert space).  Finally, \cite{CoGuaTess}  only addresses the non-path-dependent (or Markovian) case, and  does not treat the Hamilton-Jacobi-Bellman (HJB) equation related to the stochastic control problem.

The stochastic optimal control problem  \eqref{Intro_valufunc} is studied by means of the so-called randomization method. This latter  is  a purely probabilistic methodology which allows to prove directly, starting from the definition of $V_0$, that the value itself admits a representation formula (also known as non-linear Feynman-Kac formula) in terms of a suitable backward stochastic differential equation, avoiding completely analytical tools, as for instance the Hamilton-Jacobi-Bellman equation or viscosity solutions techniques.

This procedure was previously  applied in  \cite{FP15} and \cite{BaCoFuPh}, where  a stochastic control problem in finite dimension for diffusive processes  (without jumps) was addressed.    We also mention \cite{khapha12}, which has inspired \cite{FP15} and \cite{BaCoFuPh}, where a non-linear Feynman-Kac formula for the value function of a  jump-diffusive finite-dimensional stochastic control problem is provided.
Notice, however, that the methodology implemented in \cite{khapha12} (and adapted in various different framework, see e.g. \cite{BandiniFuhrman}, \cite{ChoukrounCosso}, \cite{CoPhXi}) is quite different and requires more restrictive assumptions; as a matter of fact, there the authors find the BSDE representation passing through the Hamilton-Jacobi-Bellman equation, and in particular using viscosity solutions techniques; moreover, in order to apply the techniques in \cite{khapha12},
 one already needs to know that the value function is the unique viscosity solution to the HJB equation.

The randomization method developed in the present paper improves considerably the methodology used in \cite{khapha12} and allows to extend the results in  \cite{FP15} and \cite{BaCoFuPh} to  the infinite dimensional  jump-diffusive framework, addressing, in addition, the path-dependent case. We notice that it would be possible to consider a path-dependence, or delay, in the control variable as well; however, in order to make the presentation more understandable and effective, we assume a path-dependence only in the state variable. We underline that  our results are also relevant  for the finite-dimensional case, as it is the first time the randomization method is implemented when a jump component appears in the state process dynamics.

Roughly speaking, the key idea of the randomization method consists in randomizing the control process $\alpha$, by replacing it with an uncontrolled pure jump process $I$ associated with a Poisson random measure $\theta$, independent of $W$ and $\pi$; for the pair of processes  $(X,I)$, a new randomized intensity-control problem is then introduced in such a way that the corresponding   value coincides with the original one. The  idea of this control randomization procedure
comes from the well-known methodology implemented in \cite{80Krylov} to prove the dynamic programming principle, which is based on the use of piece-wise constant policies. More specifically, in \cite{80Krylov} it is shown (under quite
 general assumptions; the only not usual assumption is the continuity of all coefficients with respect to the control variable) that the supremum over all admissible controls $\alpha$ can be replaced by the supremum over a suitable class of piece-wise constant policies. This allows to prove in a relatively easy but rigorous manner the dynamic programming principle, see Theorem III.1.6 in \cite{80Krylov}. Similarly, in the randomization method we prove (Theorem \ref{Main_thm}), under quite
general assumptions (the only not usual assumption is still the continuity of all coefficients with respect to the control variable), that we can optimize over a suitable class of piece-wise constant policies, whose dynamics is now described by the Poisson random measure $\theta$. This particular class of policies allows to prove the BSDE representation (Theorem \ref{ThmExistence}), as well as the randomized dynamic programming principle. Notice that in the present paper we have made an effort to simplify various arguments in the proof of Theorem \ref{Main_thm} and streamline the exposition.

In the Markovian case (Section \ref{S:HJB}), namely when the coefficients are non-path-dependent, we consider a family of stochastic control problems, one for each $(t,x)\in[0,T]\times H$, and define the corresponding value function. Then, exploiting the BSDE representation derived in Section \ref{S:BSDE}, we are able to prove the so-called randomized dynamic programming principle (Theorem \ref{ThmDPP}), which is as powerful as the classical dynamic programming principle, in the sense that it allows to prove (Proposition \ref{P:ViscPropvm}) that the value function is a viscosity solution to the Hamilton-Jacobi-Bellman equation, which turns out to be a second-order fully non-linear integro-differential equation in the Hilbert space $H$:
\begin{equation}
\!\!\!\!\!\!\!\!\!\!\begin{cases}
v_t + \langle Ax,D_xv\rangle + \sup_{a\in\Lambda} \Big\{\frac{1}{2}\text{Tr}\big(\sigma(t,x,a)\sigma^*(t,x,a) D_x^2v\big) +\langle b(t,x,a),D_xv\rangle +f(t,x,a) \\
+\!\int_{U \setminus \{0\}} (v(t,x+ \gamma(t,x,a,z))-v(t,x)-D_x v(t,x) \gamma(t,x,a,z))\lambda_\pi(dz)\!\Big\} = 0,  \;\text{on }(0,T)\!\times\!H,   \\
v(T,x) = g(x), \hspace{10cm} x\in H.\label{E:HJB_intro}
\end{cases}
\end{equation}
Notice that in the non-diffusive case, namely when $\sigma\equiv0$, the control problem corresponding to equation \eqref{E:HJB_intro} has already been studied in \cite{swiechzabczyk16}. Here the authors prove rigorously the (classical) dynamic programming principle (Theorem 4.2 in \cite{swiechzabczyk16}) and show that the value function solves in the viscosity sense equation \eqref{E:HJB_intro} (with $\sigma\equiv0$), Theorem 5.4 in \cite{swiechzabczyk16}. Then, Theorem \ref{ThmDPP} below, which provides the randomized dynamic programming principle, can be seen as a generalization of Theorem 4.2 in \cite{swiechzabczyk16}; similarly, Proposition \ref{P:ViscPropvm} extends Theorem 5.4 in \cite{swiechzabczyk16} to the case with $\sigma$ not necessarily equal to zero. Finally, we recall \cite{swiechzabczyk13}, which is devoted to the proof of a comparison principle for viscosity solutions to equation \eqref{E:HJB_intro} (with $\sigma$ not necessarily equal to zero), to which we refer in Remark \ref{R:Uniq}.

The paper is organized as follows. In Section \ref{S:Notations} we introduce the notations used in the paper and state the assumptions imposed on the coefficients (notice however that in the last section, namely Section \ref{S:HJB}, concerning the Markovian case, we introduce a different set of assumptions and introduce some additional notations). Section \ref{Sec_OptimalControlProblem} is devoted to the formulation of the stochastic optimal control problem, while in Section \ref{S:Randomized} we introduce the so-called randomized control problem, which allows to prove one of our main results, namely Theorem \ref{Main_thm}. In Section \ref{S:BSDE} we prove the BSDE representation of the value $V_0$ (Theorem \ref{ThmExistence}). Finally, Section \ref{S:HJB} is devoted to the study of the non-path-dependent (or Markovian) case, where we prove that the value function satisfies the randomized dynamic programming principle (Theorem \ref{ThmDPP}) and we show that it is a viscosity solution to the corresponding Hamilton-Jacobi-Bellman equation (Proposition \ref{P:ViscPropvm}).

\section{Notations and assumptions}
\label{S:Notations}

Let $H$, $U$ and  $\Xi$ be two real separable Hilbert spaces equipped with their respective Borel $\sigma$-algebrae. We denote by $|\cdot|$ and $\langle\cdot,\cdot\rangle$ (resp. $|\cdot|_U$, $|\cdot|_\Xi$ and $\langle\cdot,\cdot\rangle_\Xi$, $\langle\cdot,\cdot\rangle_U$) the norm and scalar product in $H$ (resp. in $U$ and $\Xi$). Let $(\Omega,\Fc,\P)$ be a complete probability space on which are defined a random variable $x_0\colon\Omega\rightarrow H$, a cylindrical Brownian motion $W=(W_t)_{t\geq0}$ with values in $\Xi$, and a Poisson random measure  $\pi(dt \,dz)$ on $[0,\,\infty) \times U$ with compensator $\lambda_{\pi}(dz)\,dt$. We assume that $x_0$, $W$, $\pi$ are independent. We denote by $\mu_0$ the law of $x_0$, which is a probability measure  on the Borel subsets of $H$. We also denote by $\F^{x_0,W,\pi}= (\mathcal F_t^{x_0,W,\pi})_{t \geq 0}$ the $\P$-completion of the filtration generated by $x_0$, $W$, $\pi$, which turns out to be also right-continuous, as it follows for instance from Theorem 1 in \cite{he_wang}. So, in particular, $\F^{x_0,W,\pi}$ satisfies the usual conditions. When $x_0$ is deterministic (that is, $\mu_0$ is the Dirac measure $\delta_{x_0}$) we denote $\F^{x_0,W,\pi}$ simply by $\F^{W,\pi}$.

Let $L(\Xi;H)$ be the Banach space of bounded linear operators $P\colon\Xi\rightarrow H$, and let $L_2(\Xi;H)$ be the Hilbert space of Hilbert-Schmidt operators $P\colon\Xi\rightarrow H$.

Let $T>0$ be a finite time horizon. For every $t\in[0,T]$, we consider the Banach space $D([0,t];H)$ of c\`adl\`ag maps $\x\colon[0,t]\rightarrow H$ endowed with the supremum norm $\x_t^*:=\sup_{s\in[0,t]}|\x(s)|$; when $t=T$ we also use the notation $\|\x\|_\infty:=\sup_{s\in[0,T]}|\x(s)|$. On $D([0,T];H)$ we define the canonical filtration $(\mathcal D_t^0)_{t\in[0,T]}$, with $\mathcal D^0_t$ generated by the coordinate maps
\begin{align*}
\Pi_s \ \colon \ D([0,T];H) \ &\rightarrow  H, \\
\x(\cdot) \ &\mapsto \ \x(s),
\end{align*}
for all $s\in[0,t]$. We also define its right-continuous version $(\mathcal D_t)_{t\in[0,T]}$, that is $\mathcal D_{t} = \cap_{s >t} \mathcal D^0_s$ for every $t\in[0,T)$ and $\mathcal D_T=\mathcal D_T^0$. Then, we denote by $Pred(D([0,T];H))$ the predictable $\sigma$-algebra on $[0,T]\times D([0,T];H)$ associated with  the filtration $(\mathcal D_t)_{t\in[0,T]}$.

Let $\Lambda$ be a Borel space, namely a topological space homeomorphic to a Borel subset of a Polish space. We denote by $\Bc(\Lambda)$ the Borel $\sigma$-algebra of $\Lambda$. We also denote by $d_\Lambda$ a bounded distance on $\Lambda$.

Let $A\colon\mathcal{D}(A)\subset H\to H$ be a linear operator and consider the maps $b\colon[0,T]\times D([0,T];H)\times \Lambda\rightarrow H$, $\sigma\colon[0,T]\times D([0,T];H)\times \Lambda\rightarrow L(\Xi;H)$, $\gamma: [0,T] \times D([0,T];H) \times \Lambda \times U \rightarrow H$, $f\colon[0,T]\times D([0,T];H)\times \Lambda\rightarrow\R$, $g\colon D([0,T];H)\rightarrow\R$, on which we impose the following assumptions.

\vspace{3mm}

\noindent {\bf (A)}
\begin{itemize}
\item [(i)] $A$ generates a strongly continuous semigroup $\{e^{tA},\ t\geq 0\}$ in $H$.
\item [(ii)]  $\mu_0$, the law of $x_0$, satisfies $\int_H|x|^{p_0}\mu_0(dx)<\infty$
for some $p_0\ge \max(2,2\bar p)$, with the same $\bar p\geq0$ as in \eqref{PolynomialGrowth_f_g} below.
\item[(iii)] There exists a Borel measurable function $\rho\colon U\rightarrow\R$, bounded on bounded subsets of $U$, such that
\[
\inf_{|z|_{U}> R} \rho(z) \ > \ 0, \quad \text{for every }R \ > \ 0 \qquad\quad \text{ and } \qquad\quad
\int_U |\rho(z)|^2 \lambda_\pi(dz) \ < \ \infty.
\]
\item [(iv)] The maps $b$ and $f$ are $Pred(D([0,T];H))\otimes\Bc(\Lambda)$-measurable. For every $v\in H$, the map $\sigma(\cdot,\cdot,\cdot)v\colon[0,T]\times D([0,T];H)\times \Lambda\rightarrow H$ is $Pred(D([0,T];H))\otimes\Bc(\Lambda)$-measurable. The map  $\gamma$ is $Pred(D([0,T];H))\otimes\Bc(\Lambda)\otimes\Bc(U)$-measurable. The map $g$ is $\mathcal D_T$-measurable.
\item [(v)] The map $g$ is continuous on $D([0,T];H)$ with respect to the supremum norm.
For every $t\in[0,T]$, the maps $b_t(\cdot,\cdot)$ and $f_t(\cdot,\cdot)$ are continuous on $D([0,T];H)\times \Lambda$. For every $(t,z)\in[0,T]\times U$, the map $\gamma_t(\cdot,\cdot,z)$ is continuous on $D([0,T];H)\times \Lambda$. For every $t\in[0,T]$ and any $s\in(0,T]$, we have $e^{sA}\sigma_t(\x,a)\in L_2(\Xi;H)$, for all $(\x,a)\in D([0,T];H)\times \Lambda$, and the map $e^{sA}\sigma_t(\cdot,\cdot)\colon D([0,T];H)\times \Lambda\rightarrow L_2(\Xi;H)$ is continuous.
\item [(vi)] For all $t\in[0,T]$, $s\in(0,T]$, $\x,\x'\in D([0,T];H)$, $a\in \Lambda$,
\begin{align}
|b_t(\x,a) - b_t(\x',a)| + |e^{sA}\sigma_t(\x,a)-e^{sA}\sigma_t(\x',a)|_{L_2(\Xi;H)}
\ & \leq \  L   (\x-\x')^*_t, \notag \\
|\gamma_t(\x,a,z)-\gamma_t(\x',a,z)| \ &\leq \ L\,\rho(z) (\x-\x')^*_t, \notag \\
|b_t(0,a)| + |\sigma_t(0,a)|_{L_2(\Xi;H)}
 \ & \leq \ L, \notag
\\
|\gamma_t(0,a,z)| \ & \leq \ L\,\rho(z), \notag
\\
|f_t(\x,a)| + |g(\x)| \ & \leq \ L \big(1 + \|\x\|_{_\infty}^{\bar p} \big), \label{PolynomialGrowth_f_g}
\end{align}
for some constants $L\geq0$ and $\bar p\geq0$.
\end{itemize}

\section{Stochastic optimal control problem}\label{Sec_OptimalControlProblem}

In the present section we formulate the original stochastic optimal control problem on two different probabilistic settings. More precisely, we begin formulating (see subsection \ref{SubS:Standard} below) such a control problem in a standard way, using the probabilistic setting previously introduced. Afterwards, in subsection \ref{SubS:RandomizedSetting} we formulate it on the so-called randomized probabilistic setting (that will be used for the rest of the paper and, in particular, for the formulation of the randomized control problem in Section \ref{S:Randomized}). Finally, we prove that the two formulations have the same value.

\subsection{Formulation of the control problem}\label{SubS:Standard}

We formulate the stochastic optimal control problem on the probabilistic setting introduced in Section \ref{S:Notations}. An admissible control process will be any $\F^{x_0,W,\pi}$-predictable process $\alpha$ with values in $\Lambda$. The set of all admissible control processes is denoted by $\mathcal A$. The controlled state process satisfies the following equation on $[0,T]$:
\begin{equation}
\begin{cases}
\vspace{1mm}\displaystyle dX_t = AX_t\,dt + b_t(X,\alpha_t)dt + \sigma_t(X,\alpha_t)dW_t+\!\! \int_{U \setminus \{0\}}\!\!\!\!\gamma_t(X, \alpha_t, z)\big(\pi(dt\,dz)-\lambda_\pi(dz)\,dt\big), \\
\displaystyle X_0 = x_0, \label{SDE}
\end{cases}
\end{equation}
We look for a \emph{mild solution} to the above equation \eqref{SDE} in the sense of the following definition.
\begin{Definition}
Let $\alpha \in \mathcal A$. We say that a c\`adl\`ag  $\F^{x_0,W,\pi}$-adapted stochastic process $X=(X_t)_{t\in[0,T]}$ taking values in $H$ is a \textbf{mild solution} to equation \eqref{SDE} if, $\P$-a.s.,
\begin{align*}
X_t \ &= \ e^{tA}\,x_0 + \int_0^t e^{(t-s)A}\,b_s(X,\alpha_s)\,ds + \int_0^t e^{(t-s)A}\,\sigma_s(X,\alpha_s)\,dW_s\\
&\quad \ + \int_0^t\int_{U \setminus \{0\}} e^{(t-s)A}\,\gamma_s(X, \alpha_s, z)\,\big(\pi(ds\,dz)-\lambda_\pi(dz)\,ds\big), \qquad\quad \text{for all }0\leq t\leq T.
\end{align*}
\end{Definition}
\begin{Proposition}\label{P:SDE}
Under assumption {\bf (A)}, for every  $\alpha \in \mathcal A$, there exists a unique mild solution $X^{x_0,\alpha}=(X_t^{x_0,\alpha})_{t\in[0,T]}$ to equation \eqref{SDE}. Moreover, for every $1\leq p\leq p_0$,
\begin{equation}\label{EstXalpha}
\E\Big[\sup_{t\in[0,T]}|X_t^{x_0,\alpha}|^p\Big] \ \leq \ C_p\,\big(1 + \E\left[|x_0|^p\right]\big),
\end{equation}
for some positive constant $C_p$, independent of $x_0$ and $\alpha$.
\end{Proposition}
\textbf{Proof.}
Under assumption {\bf (A)}, the existence of a unique mild solution $X^{x_0,\alpha}=(X_t^{x_0,\alpha})_{t\in[0,T]}$ to equation \eqref{SDE}, for every  $\alpha \in \mathcal A$, can be obtained by a fixed point argument proceeding as in Theorem 3.4 in \cite{swiechzabczyk13}, taking into account the fact that the coefficients of equation \eqref{SDE} are path-dependent.

We now prove estimate \eqref{EstXalpha}. In the sequel, we denote by $C$ a positive constant depending only on $T$ and $p$, independent of $x_0$ and $\alpha$, that may vary from line to line. For brevity we will denote $X^{x_0,\alpha}$ simply by $X$.
 We start by noticing that
\begin{align}\label{Ibis}
&\E\Big[\sup_{t\in[0,T]}|X_t|^p\Big]^{1/p} \leq \E\Big[\sup_{t\in[0,T]}|e^{tA}\,x_0|^p\Big]^{1/p} + \E\Big[\sup_{t\in[0,T]} \Big|\int_0^t e^{(t-s)A}\,b_s(X,\alpha_s)\,ds\Big|^p\Big]^{1/p}\notag\\
& + \E\Big[\sup_{t\in[0,T]} \Big|\int_0^t e^{(t-s)A}\,\sigma_s(X,\alpha_s)\,dW_s\Big|^p \Big]^{1/p}\notag\\
& + \E\Big[\sup_{t\in[0,T]} \Big|\int_0^t\int_{U \setminus \{0\}} e^{(t-s)A}\,\gamma_s(X, \alpha_s, z)\, (\pi(ds\,dz)-\lambda_\pi(dz)\,ds)\Big|^p\Big]^{1/p}.
\end{align}
On the other hand, by the Burk\"older-Davis-Gundy inequalities, we have
\begin{align}\label{sigma1}
&\E\Big[\sup_{t\in[0,T]} \Big|\int_0^t e^{(t-s)A}\,\sigma_s(X,\alpha_s)dW_s\Big|^p \Big]^{1/p} \leq C \E\Big[ \Big(\int_0^T e^{2(t-s)A}|\sigma_s(X,\alpha_s)|^2 ds\Big)^{p/2}\Big]^{1/p} \\
&= C \Big|\Big|\int_0^T e^{2(t-s)A}|\sigma_s(X,\alpha_s)|^2 ds\Big|\Big|^{1/2}_{L^{p/2}(\Omega, \mathcal F, \P)} \leq C \Big(\int_0^T \E\Big[e^{p(t-s)A}|\sigma_s(X,\alpha_s)|^p\Big]^{2/p} ds\Big)^{1/2}, \notag
\end{align}
and
\begin{align}\label{gamma1}
 &\E\Big[\sup_{t\in[0,T]} \Big|\int_0^t\int_{U \setminus \{0\}} e^{(t-s)A}\,\gamma_s(X, \alpha_s, z)\, (\pi(ds\,dz)-\lambda_\pi(dz)\,ds)\Big|^p\Big]^{1/p}\notag\\
 &\leq C \, \E\Big[ \Big(\int_0^T||\phi_s||^2_{L^2(U, \lambda_\pi;H)}\,ds\Big)^{p/2}\Big]^{1/p}= C \, \Big|\Big|\int_0^T||\phi_s||^2_{L^2(U, \lambda_\pi;H)}\,ds\Big|\Big|^{1/2}_{L^{p/2}(\Omega, \mathcal F, \P)}\notag\\
 &\leq C \, \Big(\int_0^T \E\Big[\Big |||\phi_s||^p_{L^2(U, \lambda_\pi;H)}\Big]^{2/p} ds\Big)^{1/2}
\end{align}
where  we have set $||\phi_s||_{L^2(U, \lambda_\pi;H)}=\Big(\int_{U \setminus \{0\}}\,|\phi_s(z)|^2 \lambda_\pi(dz)\Big)^{1/2}$ and  $\phi_s(z)=e^{(t-s)A}\,\gamma_s(X, \alpha_s, z)$.
By \eqref{sigma1}, \eqref{gamma1}, together with assumption {\bf (A)}, we get
\begin{align}\label{Isigma}
&\E\Big[\sup_{t\in[0,T]} \Big|\int_0^t e^{(t-s)A}\,\sigma_s(X,\alpha_s)\,dW_s\Big|^p \Big]^{1/p} \leq  C \, \Big(\int_0^T \E\Big[\Big(1 + \sup_{r \in [0,\,s]} |X_r|\Big)^p\Big]^{2/p} ds\Big)^{1/2} \notag\\
&\leq C \Big( 1  +  \Big(\int_0^T\E\Big[\sup_{r \in [0,\,s]} |X_r|^p\Big]^{2/p} ds\Big)^{1/2}\Big)
\end{align}
and
\begin{align}\label{Igamma}
 &\E\Big[\sup_{t\in[0,T]} \Big|\int_0^t\int_{U \setminus \{0\}} e^{(t-s)A}\,\gamma_s(X, \alpha_s, z)\, (\pi(ds\,dz)-\lambda_\pi(dz)\,ds)\Big|^p\Big]^{1/p}\notag\\
 &\leq C\, \Big(\int_0^T e^{2(t-s)A}\, \E\Big[ \Big(1 + \sup_{r \in [0,\,s]} |X_r|\Big)^p\Big(\int_{U \setminus \{0\}}|\rho(z)|^2  \lambda_\pi(dz)\Big)^{p/2}\Big]^{2/p} ds\Big)^{1/2}\notag\\
  &\leq C\, \Big(\int_0^T  \E\Big[ \Big(1 + \sup_{r \in [0,\,s]} |X_r|\Big)^p \Big]^{2/p} ds\Big)^{1/2} \!\! \leq C \Big( 1  +  \Big(\int_0^T\!\!\E\Big[\sup_{r \in [0,\,s]} |X_r|^p\Big]^{2/p}\!\!ds\Big)^{1/2}\Big).
\end{align}
Moreover, using again assumption {\bf (A)},
\begin{align}\label{Ib}
	&\E\Big[\sup_{t\in[0,T]} \Big|\int_0^t e^{(t-s)A}\,b_s(X,\alpha_s)\,ds\Big|^p\Big]^{1/p} \leq  \int_0^T \E\Big[e^{p(t-s)A}\,|b_s(X,\alpha_s)|^p\Big]^{1/p}\,ds \notag\\
	&\leq C \,\int_0^T  \E\Big[\Big(1 + \sup_{r \in [0,\,s]} |X_r|\Big)^p\Big]^{1/p}ds\leq C \Big(1 +   \int_0^T\E\Big[\sup_{r \in [0,\,s]} |X_r|^p\Big]^{1/p} ds\Big).
\end{align}
Therefore, plugging \eqref{Isigma}, \eqref{Igamma} and \eqref{Ib} in \eqref{Ibis}, we get
\begin{align*}
\E\Big[\sup_{t\in[0,T]}|X_t|^p\Big]^{1/p} \leq \ C \,\E\Big[|x_0|^p\Big]^{1/p} + C \Big(1   +   \int_0^T\E\Big[\sup_{r \in [0,\,s]} |X_r|^p\Big]^{1/p} ds&\Big)\\
 + \, C \Big(\int_0^T\E\Big[\sup_{r \in [0,\,s]} |X_r|^p\Big]^{2/p} ds&\Big)^{1/2}.
\end{align*}
Taking the square of both sides and using the Cauchy-Schwarz inequality, we find (we set $\psi_s = \E[\sup_{r \in [0,\,s]} |X_r|^p]^{2/p}$)
\begin{align*}
\psi_T \leq \,\E\Big[|x_0|^p\Big]^{2/p} + C \Big( 1+    \int_0^T \psi_s\, ds\Big),
\end{align*}
and we conclude by the Gronwall inequality.
\ep

\vspace{3mm}

\noindent The controller aims at maximizing over all $\alpha \in \mathcal A$ the gain functional
\[
J(\alpha) \ = \ \E\bigg[\int_0^T f_t(X^{x_0,\alpha},\alpha_t)\,dt + g(X^{x_0,\alpha})\bigg].
\]
By assumption \eqref{PolynomialGrowth_f_g} and estimate \eqref{EstXalpha}, we notice that $J(\alpha)$ is always finite. Finally, the value of the stochastic control problem is given by
\[
V_0 \ = \ \sup_{\alpha\in\mathcal A} J(\alpha).
\]

\subsection{Formulation of the control problem in the randomized setting}
\label{SubS:RandomizedSetting}

We formulate the stochastic optimal control problem on a new probabilistic setting that we now introduce, to which we refer as randomized probabilistic setting. Such a setting will be used for the rest of the paper and, in particular, in Section \ref{S:Randomized} for the formulation of the randomized stochastic optimal control problem.

We consider a new complete probability space $(\hat \Omega, \hat{\mathcal F}, \hat \P)$ on which are defined a random variable $\hat x_0\colon\hat\Omega\rightarrow H$, a cylindrical Brownian motion $\hat W=(\hat W_t)_{t\geq0}$ with values in $\Xi$, a Poisson random measure $\hat\pi(dt \,dz)$ on $[0,\,\infty) \times U$ with compensator $\lambda_{\pi}(dz)\,dt$ (with $\lambda_\pi$ as in Section \ref{S:Notations}), and also a Poisson random measure $\hat\theta(dt \,da)$ on $[0,\,\infty) \times\Lambda$ with compensator $\lambda_0(da)\,dt$ (on $\lambda_0$ we impose assumption {\bf (A${}_{\mathcal R}$)}(i) below). We assume that $\hat x_0$, $\hat W$, $\hat \pi$, $\hat\theta$ are independent. We denote by $\mu_0$ the law of $\hat x_0$ (with $\mu_0$ as in Section \ref{S:Notations}). We also denote by $\hat\F^{\hat x_0,\hat W,\hat\pi,\hat\theta}= (\hat{\mathcal F}_t^{\hat x_0,\hat W,\hat \pi,\hat\theta})_{t \geq 0}$ (resp. $\hat\F^{\hat\theta}= (\hat{\mathcal F}_t^{\hat\theta})_{t \geq 0}$) the $\hat\P$-completion of the filtration generated by $\hat x_0$, $\hat W$, $\hat\pi$, $\hat\theta$ (resp. $\hat\theta$), which satisfies the usual conditions. Moreover, we define $\Pc(\hat\F^{\hat x_0, \hat W,\hat \pi,\hat\theta})$ as the predictable $\sigma$-algebra on $[0,T]\times\hat\Omega$ associated with $\hat\F^{\hat x_0, \hat W,\hat\pi, \hat\theta}$. Finally, we denote by $\hat{\mathcal A}$ the family of all admissible control processes, that is the set of all $\Pc(\hat\F^{\hat x_0, \hat W,\hat \pi,\hat\theta})$-measurable maps $\hat\alpha\colon[0,T]\times\hat\Omega\rightarrow\Lambda$.

We impose the following additional assumptions.

\vspace{3mm}

\noindent {\bf (A${}_{\mathcal R}$)}
\begin{itemize}
\item [(i)] $\lambda_0$ is a finite positive measure on $\Bc(\Lambda)$, the Borel subsets of $\Lambda$, with full topological support.
\item [(ii)]    $a_0$ is a fixed point in $\Lambda$.
\end{itemize}

Similarly to Proposition \ref{P:SDE}, for every admissible control $\hat\alpha\in\hat\Ac$, we can prove the following result.

\begin{Proposition}
Under assumptions {\bf (A)}-{\bf (A${}_{\mathcal R}$)}, for every  $\hat\alpha \in\hat{\mathcal A}$, there exists a unique mild solution $\hat X^{\hat x_0,\hat \alpha}=(\hat X_t^{\hat x_0,\hat\alpha})_{t\in[0,T]}$ to equation \eqref{SDE} with $x_0$, $W$, $\pi$, $\alpha$ replaced respectively by $\hat x_0$, $\hat W$, $\hat\pi$, $\hat\alpha$. Moreover, for every $1\leq p\leq p_0$,
\begin{equation*}
\hat\E\Big[\sup_{t\in[0,T]}|\hat X_t^{\hat x_0,\hat\alpha}|^p\Big] \ \leq \ C_p\,\big(1 + \hat\E\left[|\hat x_0|^p\right]\big),
\end{equation*}
with the same constant $C_p$ as in Proposition \ref{P:SDE}, where $\hat\E$ denotes the expectation under $\hat\P$.
\end{Proposition}

\vspace{3mm}

In the present randomized probabilistic setting the formulations of the control problem reads as follows: the controller aims at maximizing over all $\hat\alpha \in\hat{\mathcal A}$ the gain functional
\begin{equation}\label{hatJ}
\hat J(\hat\alpha) \ = \ \hat\E\bigg[\int_0^T f_t(\hat X^{\hat x_0,\hat\alpha},\hat\alpha_t)\,dt + g(\hat X^{\hat x_0,\hat\alpha})\bigg].
\end{equation}
The corresponding value is defined as
\begin{equation}\label{value_hat}
\hat V_0 \ = \ \sup_{\hat\alpha\in\hat{\mathcal A}} \hat J(\hat\alpha).
\end{equation}

\begin{Proposition}\label{P:V0_Primal}
Under assumptions {\bf (A)}-{\bf (A${}_{\mathcal R}$)}, the following equality holds:
\[
V_0 \ = \ \hat V_0.
\]
\end{Proposition}
\textbf{Proof.}
The proof is organized as follows:
\begin{enumerate}[1)]
\item firstly we introduce a new probabilistic setting in product form on which we formulate the control problem \eqref{value_hat} and denote the new value function $\bar V_0$; then, we show that $\hat V_0=\bar V_0$;
\item we prove that $V_0=\bar V_0$.
\end{enumerate}
\emph{Step 1.} Let $(\Omega',\Fc',\P')$ be another complete probability space where a Poisson random measure $\theta$ on $[0,\infty)\times\Lambda$, with intensity $\lambda_0(da)dt$, is defined. Denote $\bar\Omega=\Omega\times\Omega'$, $\bar\Fc$ the completion of $\Fc\otimes\Fc'$ with respect to $\P\otimes\P'$, and $\bar\P$ the extension of $\P\otimes\P'$ to $\bar\Fc$.
Notice that $x_0,W,\pi$, which are defined on $\Omega$, as well as $\theta$, which is defined on $\Omega'$, admit obvious extensions to $\bar\Omega$. We denote those extensions by $\bar x_0,\bar W,\bar\pi,\bar\theta$. Let $\bar\F^{\bar x_0,\bar W,\bar\pi}=(\bar\Fc_t^{\bar x_0,\bar W,\bar\pi})_{t\geq0}$ (resp. $\bar\F^{\bar x_0,\bar W,\bar\pi,\bar\theta}=(\bar\Fc_t^{\bar x_0,\bar W,\bar\pi,\bar\theta})_{t\geq0}$) be the $\bar\P$-completion of the filtration generated by $\bar x_0$, $\bar W$, $\bar\pi$ (resp. $\bar x_0$, $\bar W$, $\bar\pi$, $\bar\theta$). Finally, let $\bar\Ac$ (resp. $\bar\Ac^{\bar\theta}$) be the set of $A$-valued $\bar\F^{\bar x_0,\bar W,\bar\pi}$-predictable ($\bar\F^{\bar x_0,\bar W,\bar\pi,\bar\theta}$-predictable) stochastic processes. Notice that $\bar\Ac\subset\bar\Ac^{\bar\theta}$.

For any $\bar\alpha\in\bar\Ac^{\bar\theta}$ define (with $\bar\E$ denoting the expectation under $\bar\P$)
\[
\bar J(\bar\alpha) \ = \ \bar\E\bigg[\int_0^T f_t(\bar X^{\bar x_0,\bar\alpha},\bar\alpha_t)\,dt + g(\bar X^{\bar x_0,\bar\alpha})\bigg],
\]
where $\bar X^{\bar x_0,\bar\alpha}=(\bar X_t^{\bar x_0,\bar\alpha})_{t\geq0}$ denotes the stochastic process on $\bar\Omega$, mild solution to equation \eqref{SDE}, with $\alpha$, $x_0$, $W$, $\pi$ replaced respectively by $\bar\alpha$, $\bar x_0$, $\bar W$, $\bar\pi$. We define the value function
\[
\bar V_0 \ = \ \sup_{\bar\alpha\in\bar{\mathcal A}^{\bar\theta}} \bar J(\bar\alpha).
\]
Finally, we notice that $\hat V_0=\bar V_0$. As a matter of fact, the only difference between the control problems with value functions $\hat V_0$ and $\bar V_0$ is that they are formulated on two different probabilistic settings. Given any $\hat\alpha\in\hat\Ac$, it is easy to see (by a monotone class argument) that there exists $\bar\alpha\in\bar\Ac^{\bar\theta}$ such that $(\hat\alpha,\hat x_0,\hat W,\hat\pi,\hat\theta)$ has the same law as $(\bar\alpha,\bar x_0,\bar W,\bar\pi,\bar\theta)$, so that $\hat J(\hat\alpha)=\bar J(\bar\alpha)$, which implies $\hat V_0\leq\bar V_0$. In an analogous way we get the other inequality $\hat V_0\geq\bar V_0$, from which we deduce that $\hat V_0=\bar V_0$.

\vspace{1mm}

\noindent\emph{Step 2.} Let us prove that $V_0=\bar V_0$. We begin noting that, given any $\alpha\in\Ac$, denoting by $\bar\alpha$ the canonical extension of $\alpha$ to $\bar\Omega$, we have that $\bar\alpha\in\bar\Ac$, moreover $(\alpha,x_0,W,\pi)$ has the same law as $(\bar\alpha,\bar x_0,\bar W,\bar\pi)$, so that $J(\alpha)=\bar J(\bar\alpha)$. Since $\bar\alpha\in\bar\Ac$ and $\bar\Ac\subset\bar\Ac^{\bar\theta}$, $\bar\alpha$ belongs to $\bar\Ac^{\bar\theta}$, hence $J(\alpha)=\bar J(\bar\alpha)\leq\bar V_0$. Taking the supremum over $\alpha\in\Ac$, we conclude that $V_0\leq\bar V_0$.

It remains to prove the other inequality $V_0\geq\bar V_0$. In order to prove it, we begin denoting $\bar\F^{\bar\theta}=(\bar\Fc_t^{\bar\theta})_{t\geq0}$ the $\bar\P$-completion of the filtration generated by $\bar\theta$. Notice that $\bar\Fc_t^{\bar x_0,\bar W,\bar\pi,\bar\theta}=\bar\Fc_t^{\bar x_0,\bar W,\bar\pi}\vee\bar\Fc_t^{\bar\theta}$, for every $t\geq0$. Now, fix $\bar\alpha\in\bar\Ac^{\bar\theta}$ and observe that, for every $\omega'\in\Omega'$, the stochastic process $\alpha^{\omega'}\colon\Omega\times[0,T]\rightarrow A$, defined by
\[
\alpha_t^{\omega'}(\omega) \ = \ \bar\alpha_t(\omega,\omega'), \qquad \text{for all }(\omega,\omega')\in\bar\Omega=\Omega\times\Omega',\;t\geq0,
\]
is $\F^{x_0,W,\pi}$-progressively measurable, as $\bar\alpha$ is $\bar\F^{\bar x_0,\bar W,\bar\pi,\bar\theta}$-predictable and so, in particular, $\bar\F^{\bar x_0,\bar W,\bar\pi,\bar\theta}$-progressively measurable. It is well-known (see for instance Theorem 3.7 in \cite{ChungWilliams}) that, for every $\omega'\in\Omega'$, there exists an $\F^{x_0,W,\pi}$-predictable process $\hat\alpha^{\omega'}\colon\Omega\times[0,T]\rightarrow A$ such that $\alpha^{\omega'}=\hat\alpha^{\omega'}$, $d\P\otimes dt$-a.e..

Now, recall that $\bar X^{\bar x_0,\bar\alpha}=(\bar X_t^{\bar x_0,\bar\alpha})_{t\geq0}$ denotes the mild solution to equation \eqref{SDE} on $\bar\Omega$, with $\alpha,x_0,W,\pi$ replaced respectively by $\bar\alpha,\bar x_0,\bar W,\bar\pi$. Similarly, for every fixed $\omega'\in\Omega'$, let $X^{x_0,\hat\alpha^{\omega'}}=(X_t^{x_0,\hat\alpha^{\omega'}})_{t\geq0}$ denotes the mild solution to equation \eqref{SDE} on $\Omega$, with $\alpha$ replaced by $\hat\alpha^{\omega'}$. It is easy to see that there exists a $\P'$-null set $N'\subset\Omega'$ such that, for every $\omega'\notin N'$, the stochastic processes $\bar X^{\bar x_0,\bar\alpha}(\cdot,\omega')$ and $X^{x_0,\hat\alpha^{\omega'}}(\cdot)$ solve the same equation on $\Omega$. Therefore, by pathwise uniqueness, for every $\omega'\notin N'$ we have that $\bar X^{\bar x_0,\bar\alpha}(\cdot,\omega')$ and $X^{x_0,\hat\alpha^{\omega'}}(\cdot)$ are $\P$-indistinguishable. Then, by Fubini's theorem we obtain
\[
\bar J(\bar\alpha) \ = \ \int_{\Omega'} \E\bigg[\int_0^T f_t\big(X^{x_0,\hat\alpha^{\omega'}},\hat\alpha_t^{\omega'}\big)\,dt + g\big(X^{x_0,\hat\alpha^{\omega'}}\big)\bigg] \, \P'(d\omega') \ = \ \E'\big[J\big(\hat\alpha^{\omega'}\big)\big] \ \leq \ V_0.
\]
The claim follows taking the supremum over all $\bar\alpha\in\bar\Ac^{\bar\theta}$.
\ep

\vspace{3mm}

We end this section stating a result slightly stronger than Proposition \ref{P:V0_Primal}. More precisely, we fix a $\sigma$-algebra $\hat\Gc$ independent of $(\hat x_0,\hat W,\hat\pi)$ and such that $\hat\Fc_\infty^{\hat\theta}\subset\hat\Gc$. We denote by $\hat\F^{\hat x_0,\hat W,\hat\pi,\hat\Gc}=(\hat\Fc_t^{\hat x_0,\hat W,\hat\pi,\hat\Gc})_{t\geq0}$ the $\hat\P$-completion of the filtration generated by $\hat x_0$, $\hat W$, $\hat\pi$, $\hat\Gc$ and satisfying $\hat\Gc\subset\hat\Fc_0^{\hat x_0,\hat W,\hat\pi,\hat\Gc}$. Then, we define $\hat{\mathcal A}^{\hat\Gc}$ as the family of all $\hat\F^{\hat x_0, \hat W,\hat \pi,\hat\Gc}$-predictable processes $\hat\alpha\colon[0,T]\times\hat\Omega\rightarrow\Lambda$. Notice that $\hat\Ac\subset\hat\Ac^{\hat\Gc}$.

\begin{Proposition}\label{P:V_0_G}
Under assumptions {\bf (A)}-{\bf (A${}_{\mathcal R}$)}, the following equality holds:
\[
V_0 \ = \ \sup_{\hat\alpha\in\hat\Ac^{\hat\Gc}}\hat J(\hat\alpha).
\]
\end{Proposition}
\textbf{Proof.}
We begin observing that there exists measurable space $(M,\Mc)$ and a random variable $\hat\Gamma\colon(\hat\Omega,\hat\Fc)\rightarrow(M,\Mc)$ such that $\hat\Gc=\sigma(\hat\Gamma)$ (for instance, take $(M,\Mc)=(\hat\Omega,\hat\Gc)$ and $\hat\Gamma$ the identity map). Then, the proof can be done proceeding along the same lines as in the proof of Proposition \ref{P:V0_Primal}, simply noting that the role played by $\hat\theta$ in the proof of Proposition \ref{P:V0_Primal} is now played by $\hat\Gamma$.
\ep

\section{Formulation of the randomized  control problem}
\label{S:Randomized}

We now formulate the randomized stochastic optimal control problem on the probabilistic setting introduced in subsection \ref{SubS:RandomizedSetting}. Our aim is then to prove that the value of such a control problem coincides with $V_0$ or, equivalently (by Proposition \ref{P:V0_Primal}), with $\hat V_0$. Here we simply observe that the randomized problem may depend on $\lambda_0$ and $a_0$, but its value will be independent of these two objects, as it will coincide with the value $V_0$ of the original stochastic control problem (which is independent of $\lambda_0$ and $a_0$).

We begin introducing some additional notation. We firstly notice that there exists a double sequence $(\hat T_n,\hat\eta_n)_{n\geq1}$ of $\Lambda\times(0,\infty)$-valued pairs of random variables, with $(\hat{T}_{n})_{n \geq 1}$ strictly increasing, such that the random measure $\hat \theta$ can be represented as $\hat \theta(dt\,da)=\sum_{n \geq 1} \delta_{(\hat T_n,\hat\eta_n)}(dt\,da)$. Moreover, for every Borel set $\mathscr B \in \mathcal B(\Lambda)$, the stochastic process $(\hat \theta((0,\,t] \times \mathscr B) - t \, \lambda_0(\mathscr B))_{t \geq 0}$ is a martingale under $\hat\P$. Now, we introduce the pure jump stochastic process taking values in $\Lambda$ defined as
\begin{equation}\label{I}
\hat I_t \ = \ \sum_{n \geq 0} \hat \eta_n\,1_{[\hat T_n , \hat T_{n+1})}(t), \qquad\qquad \text{for all }t \geq 0,
\end{equation}
where we set $\hat T_0:=0$ and $\hat\eta_0:=a_0$ (notice that, when $\Lambda$ is a subset of a vector space, we can write \eqref{I} simply as $\hat I_t = a_0 + \int_{0}^{t}\int_A (a- \hat I_{s-}) \, \hat \theta(ds \,da)$).

We use $\hat I$ to randomize the control in equation \eqref{SDE}, which then becomes:
\begin{equation}\label{SDE_random}
\begin{cases}
\vspace{1mm}\displaystyle d\hat X_t = A\hat X_t\,dt + b_t(\hat X,\hat I_t)dt + \sigma_t(\hat X,\hat I_t)d\hat W_t  + \!\int_{U \setminus \{0\}}\!\!\!\gamma_t(\hat X,\hat I_{t-}, z)\big(\hat\pi(dt\,dz)-\lambda_\pi(dz)dt\big),\\
\displaystyle \hat X_0 = \hat x_0.
\end{cases}
\end{equation}
As for equation \eqref{SDE}, we look for a mild solution to \eqref{SDE_random}, namely an $H$-valued c\`adl\`ag $\hat\F^{\hat x_0,\hat W,\hat\pi, \hat \theta}$-adapted  stochastic process $\hat X=(\hat X_t)_{t\in[0,T]}$ such that, $\hat\P$-a.s.,
\begin{align}\label{SDE_random2}
\hat X_t &=  e^{tA}\,\hat x_0 + \int_0^t e^{(t-s)A}\,b(\hat X,\hat I_s)\,ds + \int_0^t e^{(t-s)A}\,\sigma(\hat X,\hat I_s)\,d\hat W_s\\
&\quad \  + \, \int_0^t\int_{U \setminus \{0\}} e^{(t-s)A}\,\gamma(\hat X, \hat I_{s-}, z)\, (\hat\pi(ds\,dz)-\lambda_\pi(dz)\,ds), \qquad \text{for all }0\leq t\leq T. \notag
\end{align}
Under assumptions {\bf (A)}-{\bf (A${}_{\mathcal R}$)}, proceeding as in Proposition \ref{P:SDE}, we can prove the following result.
\begin{Proposition}
Under assumptions {\bf (A)}-{\bf (A${}_{\mathcal R}$)}, there exists a unique mild solution $\hat X=(\hat X_t)_{t\in[0,T]}$ to equation \eqref{SDE_random}, such that, for every $1\leq p\leq p_0$,
\begin{equation}\label{Est_X_rand}
\hat\E\Big[\sup_{t\in[0,T]}|\hat X_t|^p\Big] \ \leq \ C_p\,\big(1 + \hat\E\left[|\hat x_0|^p\right]\big),
\end{equation}
with the same constant $C_p$ as in Proposition \ref{P:SDE}. In addition, for every $t\in[0,T]$ and any $1\leq p\leq p_0$, we have
\begin{equation}\label{Est_X_rand_conditional}
\hat\E\Big[\sup_{s\in[t,T]}|\hat X_s|^p\Big|\hat\Fc_t^{\hat x_0,\hat W,\hat\pi,\hat\theta}\Big] \ \leq \ C_p\,\Big(1 + \sup_{s\in[0,t]}|\hat X_s|^p\Big), \qquad \hat\P\text{-a.s.}
\end{equation}
with the same constant $C_p$ as in Proposition \ref{P:SDE}.
\end{Proposition}
\textbf{Proof.}
Concerning estimate \eqref{Est_X_rand}, the proof can be done proceeding along the same lines as in the proof of Proposition \ref{P:SDE}. On the other hand, regarding estimate \eqref{Est_X_rand_conditional} we begin noting that given any two integrable $\hat\Fc_t^{\hat x_0,\hat W,\hat\pi,\hat\theta}$-measurable random variables $\eta$ and $\xi$, then the following property holds: $\eta\leq\xi$, $\hat\P$-a.s., if and only if $\hat\E[\eta\,1_E]\leq\hat\E[\xi\,1_E]$, for every $E\in\hat\Fc_t^{\hat x_0,\hat W,\hat\pi,\hat\theta}$. So, in particular, estimate \eqref{Est_X_rand_conditional} is true if and only if the following estimate holds:
\begin{equation}\label{Est_X_rand_conditional2}
\hat\E\Big[\sup_{s\in[t,T]}|\hat X_s|^p\,1_E\Big] \ \leq \ C_p\,\Big(\hat\E[1_E] + \hat\E\Big[\sup_{s\in[0,t]}|\hat X_s|^p\,1_E\Big]\Big), \quad \text{for every }E\in\hat\Fc_t^{\hat x_0,\hat W,\hat\pi,\hat\theta}.
\end{equation}
The proof of estimate \eqref{Est_X_rand_conditional2} can be done proceeding along the same lines as in the proof of Proposition \ref{P:SDE}, firstly multiplying equation \eqref{SDE_random2} by $1_E$.
\ep

\vspace{3mm}

We can now formulate the randomized control problem. The family of all admissible control maps, denoted by $\hat \Vc$, is the set of all $\Pc(\hat\F^{\hat x_0, \hat W,\hat\pi, \hat \theta})\otimes\Bc(\Lambda)$-measurable functions $\hat \nu\colon[0,T]\times\hat\Omega\times \Lambda\rightarrow(0,\infty)$ which are bounded from above and bounded away from zero, namely $0<\inf_{[0,T]\times\hat\Omega\times\Lambda}\hat\nu\leq\sup_{[0,T]\times\hat\Omega\times\Lambda}\hat\nu<+\infty$. Given $\hat \nu\in\hat \Vc$, we consider the probability measure $\hat\P^{\hat \nu}$ on $(\hat\Omega,\hat\Fc_T^{\hat x_0, \hat W,\hat\pi, \hat \theta})$ given by $d\hat\P^{\hat \nu}=\hat\kappa_T^{\hat \nu}\,d\hat\P$, where $(\hat\kappa_t^{\hat\nu})_{t\in[0,T]}$ denotes the Dol\'eans-Dade exponential
\begin{equation}\label{Doleans}
\hat\kappa_t^{\hat\nu} \ = \
\Ec_t\bigg(\int_0^\cdot\int_\Lambda \big(\hat \nu_s(a) - 1\big)\,\big(\hat\theta(ds\, da)- \lambda_0(da)\,ds\big)\bigg).
\end{equation}
By Girsanov's theorem (see e.g. Theorem 15.2.6 in \cite{Coh-Ell}), under $\hat\P^{\hat \nu}$ the $\hat\F^{\hat x_0, \hat W,\hat\pi, \hat \theta}$-compensator of $\hat \theta$ on $[0,\,T] \times \Lambda$ is $\hat \nu_s(a)\lambda_0(da)ds$.

Notice that, under $\hat\P^{\hat \nu}$, $\hat W$ remains a Brownian motion and the $\hat\F^{\hat x_0, \hat W,\hat\pi, \hat \theta}$-compensator of $\hat \pi$ on $[0,\,T] \times \Lambda$ is $\lambda_{\pi}(dz)ds$ (see e.g. Theorem 15.3.10 in \cite{Coh-Ell} or Theorem 12.31 in \cite{HeWangYan}).

As a consequence, the following generalization of estimate \eqref{Est_X_rand} holds: for every $1\leq p\leq p_0$,
\begin{equation}\label{Est_X_rand_nu}
\sup_{\hat \nu\in\hat \Vc}\,\hat\E^{\hat \nu}\Big[\sup_{t\in [0,T]}|\hat X_t|^p\Big] \ \leq \ C_p\,\big(1 + \hat\E^{\hat \nu}\big[|x_0|^p\big]\big),
\end{equation}
with the same constant $C_p$ as in \eqref{Est_X_rand}, where $\hat\E^{\hat\nu}$ denotes the expectation with respect to $\hat\P^{\hat\nu}$.

The controller aims at maximizing over all $\hat \nu\in\hat \Vc$ the gain functional
\[
\hat J^{\mathcal R}(\hat \nu) \ = \ \hat\E^{\hat \nu}
\bigg[\int_0^T f_t(\hat X,\hat I_t)\,dt + g(\hat X)\bigg].
\]
By assumption \eqref{PolynomialGrowth_f_g} and estimate \eqref{Est_X_rand_nu}, it follows that $\hat J^{\mathcal R}(\hat \nu)$ is always finite. Finally, the value function of the randomized control problem is given by
\[
\hat V_0^{\mathcal R} \ = \ \sup_{\hat \nu\in\hat \Vc} \hat J^{\mathcal R}(\hat \nu).
\]
In the sequel, we denote the probabilistic setting we have adopted for the randomized control problem shortly by the tuple $(\hat\Omega,\hat\Fc,\hat\P;\hat x_0,\hat W,\hat\pi,\hat\theta;\hat I,\hat X;\hat\Vc)$.

Our aim is now to prove that $\hat V_0^{\mathcal R}$ coincides with the value $V_0$ of the original control problem. Firstly, we state three auxiliary results:
\begin{enumerate}[1)]
\item the first result (Lemma \ref{L:NewSetting}) shows that the value $\hat V_0^{\mathcal R}$ of the randomized control problem is independent of the probabilistic setting on which the problem is formulated;
\item in Lemma \ref{L:hatJR=hatJ} we prove that there exists a probabilistic setting for the randomized control problem where $\hat J^{\mathcal R}$ can be expressed in terms of the gain functional $\hat J$ in \eqref{hatJ}; as noticed in Remark \ref{R:hatJR=hatJ}, this result allows to formulate the randomized control problem in ``strong'' form, rather than as a supremum over a family of probability measures;
\item finally, in Lemma \ref{L:Fuhrman} we prove, roughly speaking, that given any $\alpha\in\Ac$ and $\eps>0$ there exist a probabilistic setting for the randomized control and a suitable $\hat\nu$ such that the ``distance'' under $\hat\P^{\hat\nu}$ between the pure jump process $\hat I$ and $\alpha$ is less than $\eps$. In order to do it, we need to introduce the following distance on $\hat\Ac$ (see Definition 3.2.3 in \cite{80Krylov}), for every fixed $\hat\nu\in\hat\Vc$:
\[
\hat d_{\textup{Kr}}^{\hat\nu}(\hat\alpha,\hat\beta) \ := \ \hat\E^{\hat\nu}\bigg[\int_0^T d_\Lambda(\hat\alpha_t,\hat\beta_t)\,dt\bigg],
\]
for all $\hat\alpha,\hat\beta\in\hat\Ac$.
\end{enumerate}

\begin{Lemma}\label{L:NewSetting}
Suppose that assumptions {\bf (A)}-{\bf (A${}_{\mathcal R}$)} hold. Consider a new probabilistic setting for the randomized control problem characterized by the tuple $(\bar\Omega,\bar\Fc,\bar\P;\bar x_0,\bar W,\bar\pi,\bar\theta;\bar I,$ $\bar X;\bar\Vc)$. Then
\[
\hat V_0^{\mathcal R} \ = \ \bar V_0^{\mathcal R}.
\]
\end{Lemma}
\textbf{Proof.}
The proof can be done proceeding along the same lines as in the proof of Proposition 3.1 in \cite{BaCoFuPh}. Here we just recall the main steps. Firstly we take $\hat\nu\in\hat\Vc$ which admits an explicit functional dependence on $(\hat x_0,\hat W,\hat\pi,\hat\theta)$. For such a $\hat\nu$ it is easy to find $\bar\nu\in\bar\Vc$ such that $(\hat\nu,\hat x_0,\hat W,\hat\pi,\hat\theta)$ has the same law as $(\bar\nu,\bar x_0,\bar W,\bar\pi,\bar\theta)$ (simply replacing $\hat x_0,\hat W,\hat\pi,\hat\theta$ by $\bar\nu,\bar x_0,\bar W,\bar\pi,\bar\theta$ in the expression of $\hat\nu$). So, in particular, $\hat J^{\mathcal R}(\hat\nu)=\bar J^{\mathcal R}(\bar\nu)$. By a monotone class argument, we deduce that the same equality holds true for every $\hat\nu\in\hat\Vc$, which implies $\hat V_0^{\mathcal R}\leq\bar V_0^{\mathcal R}$. Interchanging the role of $(\hat\Omega,\hat\Fc,\hat\P;\hat x_0,\hat W,\hat\pi,\hat\theta;\hat I,\hat X;\hat\Vc)$ and $(\bar\Omega,\bar\Fc,\bar\P;\bar x_0,\bar W,\bar\pi,\bar\theta;\bar I,\bar X;\bar\Vc)$, we obtain the other inequality, from which the claim follows.
\ep

\begin{Lemma}\label{L:hatJR=hatJ}
Suppose that assumptions {\bf (A)}-{\bf (A${}_{\mathcal R}$)} hold. Then, there exists a probabilistic setting for the randomized control problem $(\bar\Omega,\bar\Fc,\bar\P;\bar x_0,\bar W,\bar\pi,\bar\theta;\bar I,\bar X;\bar\Vc)$ and a $\sigma$-algebra $\bar\Gc\subset\bar\Fc$, independent of $\bar x_0$, $\bar W$, $\bar\pi$, with $\bar\Fc_\infty^{\bar\theta}\subset\bar\Gc$, such that: given any $\bar\nu\in\bar\Vc$ there exists $\bar\alpha^{\bar\nu}\in\bar\Ac^{\bar\Gc}$ satisfying
\begin{align}\label{Law}
\text{Law of $(\bar x_0,(\bar W_t)_{0\leq t\leq T},\bar\pi_{_{|[0,T]\times\Lambda}},(\bar I_t)_{0\leq t\leq T})$ under $\bar\P^{\bar\nu}$}& \notag \\
= \ \text{Law of $(\bar x_0,(\bar W_t)_{0\leq t\leq T},\bar\pi_{_{|[0,T]\times\Lambda}},\bar\alpha^{\bar\nu})$ under $\bar\P$.}\hspace{.7mm}&
\end{align}
So, in particular,
\[
\bar J^{\mathcal R}(\bar\nu) \ = \ \bar J(\bar\alpha^{\bar\nu}).
\]
\end{Lemma}
\begin{Remark}
{\rm
Recall that $\bar\Ac^{\bar\Gc}$ was defined just before Proposition \ref{P:V_0_G}, even though it was denoted $\hat\Ac^{\hat\Gc}$ since it was defined in the probabilistic setting $(\hat\Omega,\hat\Fc,\hat\P;\hat x_0,\hat W,\hat\pi,\hat\theta;\hat I,\hat X;\hat\Vc)$ instead of $(\bar\Omega,\bar\Fc,\bar\P;\bar x_0,\bar W,\bar\pi,\bar\theta;\bar I,\bar X;\bar\Vc)$.
\epR
}
\end{Remark}
\textbf{Proof (of Lemma \ref{L:hatJR=hatJ}).}
Let $(\Omega,\Fc,\P;x_0,W,\pi;X;\Ac)$ be the setting of the original stochastic control problem in Section \ref{SubS:Standard}.

Proceeding along the same lines as at the beginning of Section 4.1 in \cite{BaCoFuPh}, we construct an atomless finite measure $\lambda_0'$ on $(\R,\Bc(\R))$ and a surjective Borel-measurable map $\pi\colon\R\rightarrow\Lambda$ such that $\lambda_0=\lambda_0'\circ\pi^{-1}$. Let $(\Omega',\Fc',\P')$ be the completion of the canonical probability space of a Poisson random measure $\theta'=\sum_{n\geq1}\delta_{(T_n',\rho_n')}$ on $[0,\infty)\times\Lambda$ with intensity measure $\lambda_0'(dr)dt$, where $(T_n',\rho_n')_{n\geq1}$ is the marked point process associated with $\theta'$. Then, $\theta=\sum_{n\geq1}\delta_{(T_n',\pi(\rho_n'))}$ is a Poisson random measure on $[0,\infty)\times\Lambda$ with intensity measure $\lambda_0(dr)dt$.

Let $\bar\Omega=\Omega\times\Omega'$, $\bar\Fc$ the $\P\otimes\P'$-completion of $\Fc\otimes\Fc'$, and $\bar\P$ the extension of $\P\otimes\P'$ to $\bar\Fc$. Then, we consider the corresponding probabilistic setting for the randomized control problem $(\bar\Omega,\bar\Fc,\bar\P;\bar x_0,\bar W,\bar\pi,\bar\theta;\bar I,\bar X;\bar\Vc)$, where $\bar x_0$, $\bar W$, $\bar \pi$, $\bar\theta$ denote the canonical extensions of $x_0$, $W$, $\pi$, $\theta$ to $\bar\Omega$. We also denote by $\bar\theta'$ the canonical extension of $\theta'$ to $\bar\Omega$. Let $\bar\F^{\bar\theta'}=(\bar\Fc_t^{\bar\theta'})_{t\geq0}$ (resp. $\bar\F^{\bar\theta}=(\bar\Fc_t^{\bar\theta})_{t\geq0}$) the filtration generated by $\bar\theta'$ (resp. $\bar\theta$). We define $\bar\Gc:=\bar\Fc_\infty^{\bar\theta'}$. Notice that $\bar\Fc_\infty^{\bar\theta}\subset\bar\Gc$ and $\bar\Gc$ is independent of $\bar x_0$, $\bar W$, $\bar\pi$. Finally, we denote by $\bar\F^{\bar x_0,\bar W,\bar\pi,\bar\Gc}=(\bar\Fc_t^{\bar x_0,\bar W,\bar\pi,\bar\Gc})_{t\geq0}$ the $\bar\P$-completion of the filtration generated by $\bar x_0$, $\bar W$, $\bar\pi$, $\bar\Gc$ and satisfying $\bar\Gc\subset\bar\Fc_0^{\bar x_0,\bar W,\bar\pi,\bar\Gc}$.

Now, fix $\bar\nu\in\bar\Vc$. By an abuse of notation, we still denote by $\Fc$ the canonical extension of the $\sigma$-algebra $\Fc$ to $\bar\Omega$. Then, we notice that in the probabilistic setting $(\bar\Omega,\bar\Fc,\bar\P;\bar x_0,\bar W,\bar\pi,\bar\theta;$ $\bar I,\bar X;\bar\Vc)$ just introduced \eqref{Law} follows if we prove the following: \emph{there exists $\bar\alpha^{\bar\nu}\in\bar\Ac^{\bar\Gc}$ satisfying}
\begin{equation}\label{ConditionalLaw}
\text{\emph{Conditional law of $(\bar I_t)_{0\leq t\leq T}$ under $\bar\P^{\bar\nu}$ given $\Fc$}}= \text{\emph{Conditional law of $\bar\alpha^{\bar\nu}$ under $\bar\P$ given $\Fc$}}.
\end{equation}
It only remains to prove \eqref{ConditionalLaw}. To this end, we recall that the process $\bar I$ is defined as
\[
\bar I_t \ = \ \sum_{n \geq 0} \bar\eta_n\,1_{[\bar T_n , \bar T_{n+1})}(t), \qquad\qquad \text{for all }t \geq 0,
\]
where $(\bar T_0,\bar\eta_0):=(0,a_0)$, while $(\bar T_n,\bar\eta_n)$, $n\geq1$, denotes the canonical extension of $(T_n',\pi(\rho_n'))$ to $\bar\Omega$. Then, \eqref{ConditionalLaw} follows if we prove the following: \emph{there exists a sequence $(\bar T_n^{\bar\nu},\bar\eta_n^{\bar\nu})_{n\geq1}$ on $(\bar\Omega,\bar\Fc,\bar\P)$ such that:}
\begin{enumerate}[(i)]
\item \emph{$(\bar T_n^{\bar\nu},\bar\eta_n^{\bar\nu})\colon\bar\Omega\rightarrow(0,\infty)\times\Lambda$ and $\bar T_n^{\bar\nu}<\bar T_{n+1}^{\bar\nu}$;}
\item \emph{$\bar T_n^{\bar\nu}$ is a $\bar\F^{\bar x_0,\bar W,\bar\pi,\bar\Gc}$-stopping time and $\bar\eta_n^{\bar\nu}$ is $\Fc_{\bar T_n^{\bar\nu}}^{\bar x_0,\bar W,\bar\pi,\bar\Gc}$-measurable;}
\item \emph{$\lim_{n\rightarrow\infty}\bar T_n^{\bar\nu}=\infty$;}
\item \emph{the conditional law of the sequence $(\bar T_1,\bar\eta_1)\,1_{\{\bar T_1\leq T\}}$,\,$\ldots$, $(\bar T_n,\bar\eta_n)\,1_{\{\bar T_n\leq T\}}$,\,$\ldots$ under $\bar\P^{\bar\nu}$ given $\Fc$ is equal to the conditional law of the sequence $(\bar T_1^{\bar\nu},\bar\eta_1^{\bar\nu})\,1_{\{\bar T_1^{\bar\nu}\leq T\}}$, $\ldots$ , $(\bar T_n^{\bar\nu},\bar\eta_n^{\bar\nu})\,1_{\{\bar T_n^{\bar\nu}\leq T\}}$,\,$\ldots$ under $\bar\P$ given $\Fc$.}
\end{enumerate}
As a matter of fact, if there exists $(\bar T_n^{\bar\nu},\bar\eta_n^{\bar\nu})_{n\geq1}$ satisfying (i)-(ii)-(iii)-(iv), then the process $\bar\alpha^{\bar\nu}$, defined as
\[
\bar\alpha_t^{\bar\nu} \ := \ \sum_{n \geq 0} \bar\eta_n^{\bar\nu}\,1_{[\bar T_n^{\bar\nu} , \bar T_{n+1}^{\bar\nu})}(t), \qquad\qquad \text{for all $0\leq t\leq T$, with  $(\bar T_0^{\bar\nu},\bar\eta_0^{\bar\nu}):=(0,a_0)$},
\]
belongs to $\bar\Ac^{\bar\Gc}$ and \eqref{ConditionalLaw} holds.

Finally, concerning the existence of a sequence $(\bar T_n^{\bar\nu},\bar\eta_n^{\bar\nu})_{n\geq1}$ satisfying (i)-(ii)-(iii)-(iv), we do not report the proof of this result as it can be done proceeding along the same lines as in the proof of Lemma 4.3 in \cite{BaCoFuPh}, the only difference being that the filtration $\F^W$ in \cite{BaCoFuPh} (notice that in \cite{BaCoFuPh} $W$ denotes a finite dimensional Brownian motion) is now replaced by $\F^{x_0,W,\pi}$: this does not affect the proof of Lemma 4.3 in \cite{BaCoFuPh}.
\ep

\begin{Remark}\label{R:hatJR=hatJ}
{\rm
Let $(\bar\Omega,\bar\Fc,\bar\P;\bar x_0,\bar W,\bar\pi,\bar\theta;\bar I,\bar X;\bar\Vc)$ and $\bar\Gc$ be respectively the probabilistic setting for the randomized control problem and the $\sigma$-algebra mentioned in Lemma \ref{L:hatJR=hatJ}. We denote by $\bar\Ac^{\bar\Vc}$ the family of all controls $\bar\alpha\in\bar\Ac^{\bar\Gc}$ for which there exists some $\bar\nu\in\bar\Vc$ such that $\bar J(\bar\alpha)=\bar J^{\mathcal R}(\bar\nu)$. Then, by definition $\bar\Ac^{\bar\Vc}\subset\bar\Ac^{\bar\Gc}$. Moreover, by Lemma \ref{L:hatJR=hatJ} we have the following ``strong'' formulation of the randomized control problem:
\[
\bar V_0^{\mathcal R} \ = \ \sup_{\bar\alpha\in\bar\Ac^{\bar\Vc}}\bar J(\bar\alpha).\vspace{-5mm}
\]
\epR
}
\end{Remark}

\begin{Lemma}\label{L:Fuhrman}
Suppose that assumptions {\bf (A)}-{\bf (A${}_{\mathcal R}$)} hold. For any $\alpha\in\Ac$ and $\eps>0$ there exist:
\begin{enumerate}
\item[\textup{1)}] a probabilistic setting for the randomized control problem $(\bar\Omega,\bar\Fc,\bar\P^{\alpha,\eps};\bar x_0,\bar W,\bar\pi,\bar\theta^{\alpha,\eps};$ $\bar I^{\alpha,\eps},\bar X^{\alpha,\eps};\bar\Vc^{\alpha,\eps})$ (notice that $\bar\Omega,\bar\Fc,\bar x_0,\bar W,\bar\pi$ do not depend on $\alpha,\eps$);
\item[\textup{2)}] a probability measure $\bar\Q$ on $(\bar\Omega,\bar\Fc)$ equivalent to $\bar\P^{\alpha,\eps}$, which does not depend on $\alpha,\eps$;
\item[\textup{3)}] a stochastic process $\bar\alpha\colon[0,T]\times\bar\Omega\rightarrow\Lambda$, depending only on $\alpha$ but not on $\eps$, which is predictable with respect to the $\bar\P^{\alpha,\eps}$-completion (or, equivalently, $\bar\Q$-completion) of the filtration generated by $\bar x_0$, $\bar W$, $\bar\pi$;
\item[\textup{4)}] $\bar\nu^{\alpha,\eps}\in\bar\Vc^{\alpha,\eps}$,
\end{enumerate}
such that, denoting by $\bar\P^{\bar\nu^{\alpha,\eps}}$ the probability measure\footnote{Here $\bar\F^{\bar x_0,\bar W,\bar\pi,\bar\theta^{\alpha,\eps}}=(\bar\Fc_t^{\bar x_0,\bar W,\bar\pi,\bar\theta^{\alpha,\eps}})_{t\geq0}$ denotes the $\bar\P^{\alpha,\eps}$-completion of the filtration generated by $\bar x_0,\bar W,\bar\pi,\bar\theta^{\alpha,\eps}$, while $\bar\kappa^{\bar\nu^{\alpha,\eps}}$ is the Dol\'eans-Dade exponential given by \eqref{Doleans} with $\hat\nu$, $\hat\theta$ replaced respectively by $\bar\nu^{\alpha,\eps}$, $\bar\theta^{\alpha,\eps}$.} on $(\bar\Omega,\bar\Fc_T^{\bar x_0,\bar W,\bar\pi,\bar\theta^{\alpha,\eps}})$ defined as $d\bar\P^{\bar\nu^{\alpha,\eps}}=\bar\kappa_T^{\bar\nu^{\alpha,\eps}}d\bar\P^{\alpha,\eps}$, the following properties hold:
\begin{enumerate}
\item[\textup{(i)}] the restriction of $\bar\Q$ to $\bar\Fc_T^{\bar x_0,\bar W,\bar\pi,\bar\theta^{\alpha,\eps}}$ coincides with $\bar\P^{\bar\nu^{\alpha,\eps}}$;
\item[\textup{(ii)}] the following inequality holds:
\[
\bar\E^{\bar\Q}\bigg[\int_0^T d_\Lambda(\bar I_t^{\alpha,\eps},\bar\alpha_t)\,dt\bigg] \ \leq \ \eps;
\]
\item[\textup{(iii)}] the quadruple $(x_0,W, \pi,\alpha)$ under $\P$ has the same law as $(\bar x_0,\bar W, \bar\pi,\bar\alpha)$ under $\bar\P^{\alpha,\eps}$.
\end{enumerate}
\end{Lemma}
\textbf{Proof.}
Fix $\alpha\in\Ac$ and $\eps>0$. In order to construct the probabilistic setting of item 1), we apply Proposition A.1 in \cite{BaCoFuPh} (with filtration $\G=\F^{x_0,W,\pi}$ and $\delta=\eps$), from which we deduce the existence of a probability space $(\bar\Omega,\tilde\Fc,\tilde\Q)$ independent of $\alpha,\eps$ (corresponding to $(\hat\Omega,\hat\Fc,\Q)$ in the notation of Proposition A.1) and a marked point process $(\bar T_n^{\alpha,\eps},\bar\eta_n^{\alpha,\eps})_{n\geq1}$ with corresponding random measure $\bar\theta^{\alpha,\eps}=\sum_{n\geq1}\delta_{(\bar T_n^{\alpha,\eps},\bar\eta_n^{\alpha,\eps})}$ on $\bar\Omega$ (corresponding respectively to $(\hat S_n,\hat\eta_n)_{n\geq1}$ and $\hat\mu$ in Proposition A.1) with the following properties:
\begin{enumerate}[(a)]
\item there exists a probability space $(\Omega',\Fc',\P')$ such that $\bar\Omega=\Omega\times\Omega'$, $\tilde\Fc=\Fc\otimes\Fc'$, $\tilde\Q=\P\otimes\P'$; we denote by $\bar x_0$, $\bar W$, $\bar\pi$ the natural extensions of $x_0$, $W$, $\pi$ to $\bar\Omega$ (which obviously do not depend on $\alpha,\eps$); we also denote by $\tilde\F^{\bar x_0,\bar W,\bar\pi}$ the extension of $\F^{x_0,W,\pi}$ to $\bar\Omega$;
\item denoting $\tilde\E^{\tilde\Q}$ the expectation with respect to $\tilde\Q$, we have
\[
\tilde\E^{\tilde\Q}\bigg[\int_0^T d_\Lambda(\bar I_t^{\alpha,\eps},\bar\alpha_t)\,dt\bigg] \ \leq \ \eps,
\]
where $\bar\alpha$ is the natural extension of $\alpha$ to $\bar\Omega=\Omega\times\Omega'$ (which clearly depend only on $\alpha$, not on $\eps$), while $\bar I^{\alpha,\eps}$ is given by
\[
\bar I_t^{\alpha,\eps} \ = \ \sum_{n \geq 0} \bar \eta_n^{\alpha,\eps}\,1_{[\bar T_n^{\alpha,\eps},\bar T_{n+1}^{\alpha,\eps})}(t), \qquad\qquad \text{for all }t \geq 0,
\]
with $\bar T_0^{\alpha,\eps}=0$ and $\bar\eta_0^{\alpha,\eps}=a_0$;
\item let $\tilde\F^{\bar\theta^{\alpha,\eps}}=(\tilde\Fc_t^{\bar\theta^{\alpha,\eps}})_{t\geq0}$ denote the filtration generated by $\bar\theta^{\alpha,\eps}$; let also $\Pc(\tilde\Fc_t^{\bar x_0,\bar W,\bar\pi}\vee\bar\Fc_t^{\tilde\theta^{\alpha,\eps}})$ be the predictable $\sigma$-algebra on $[0,T]\times\bar\Omega$ associated with the filtration $(\tilde\Fc_t^{\bar x_0,\bar W,\bar\pi}\vee\tilde\Fc_t^{\bar\theta^{\alpha,\eps}})_{t\geq0}$; then, there exists a $\Pc(\tilde\Fc_t^{\bar x_0,\bar W,\bar\pi}\vee\bar\Fc_t^{\tilde\theta^{\alpha,\eps}})\otimes\Bc(\Lambda)$-measurable map $\bar\nu^{\alpha,\eps}\colon[0,T]\times\bar\Omega\times\Lambda\rightarrow(0,\infty)$, with $0<\inf_{[0,T]\times\bar\Omega\times\Lambda}\bar\nu^{\alpha,\eps}\leq\sup_{[0,T]\times\bar\Omega\times\Lambda}\bar\nu^{\alpha,\eps}<+\infty$, such that under $\tilde\Q$ the random measure $\bar\theta^{\alpha,\eps}$ has $(\tilde\Fc_t^{\bar x_0,\bar W,\bar\pi}\vee\tilde\Fc_t^{\bar\theta^{\alpha,\eps}})$-compensator on $[0,\,T] \times \Lambda$ given by $\bar\nu_t^{\alpha,\eps}(a)\lambda_0(da)dt$.
\end{enumerate}
Now, proceeding as in Section 4.2 of \cite{BaCoFuPh}, we consider the completion $(\bar\Omega,\bar\Fc,\bar\Q)$ of $(\bar\Omega,\tilde\Fc,\tilde\Q)$. Then, from item (b) above we immediately deduce item (ii).

Let $\bar\F^{\bar x_0,\bar W,\bar\pi,\bar\theta^{\alpha,\eps}}=(\bar\Fc_t^{\bar x_0,\bar W,\bar\pi,\bar\theta^{\alpha,\eps}})_{t\geq0}$ be the $\bar\Q$-completion of the filtration $(\tilde\Fc_t^{\bar x_0,\bar W,\bar\pi}\vee\tilde\Fc_t^{\bar\theta^{\alpha,\eps}})_{t\geq0}$. It easy to see that under $\bar\Q$ the $\bar\F^{\bar x_0,\bar W,\bar\pi,\bar\theta^{\alpha,\eps}}$-compensator of $\bar\theta^{\alpha,\eps}$ on $[0,\,T] \times \Lambda$ is still given by $\bar\nu_t^{\alpha,\eps}(a)\lambda_0(da)dt$. Denote by $\Pc(\bar\F^{\bar x_0, \bar W,\bar \pi,\bar\theta^{\alpha,\eps}})$ the predictable $\sigma$-algebra on $[0,T]\times\bar\Omega$ associated with $\bar\F^{\bar x_0, \bar W,\bar\pi, \bar\theta^{\alpha,\eps}}$. Then, we define $\bar\Vc^{\alpha,\eps}$ as the set of all $\Pc(\bar\F^{\bar x_0, \bar W,\bar\pi, \bar\theta^{\alpha,\eps}})\otimes\Bc(\Lambda)$-measurable functions $\bar \nu\colon[0,T]\times\bar\Omega\times \Lambda\rightarrow(0,\infty)$ which are bounded from above and bounded away from zero. Notice that $\bar\nu^{\alpha,\eps}\in\bar\Vc^{\alpha,\eps}$. Let $\bar\kappa^{\bar\nu^{\alpha,\eps}}$ be the Dol\'eans-Dade exponential given by \eqref{Doleans} with $\hat\nu$, $\hat\theta$ replaced respectively by $\bar\nu^{\alpha,\eps}$, $\bar\theta^{\alpha,\eps}$. Since $\inf_{[0,T]\times\bar\Omega\times\Lambda}\bar\nu^{\alpha,\eps}>0$, it follows that $\bar\nu^{\alpha,\eps}$ has bounded inverse, so that we can define the probability measure $\bar\P^{\alpha,\eps}$ on $(\bar\Omega,\bar\Fc)$, equivalent to $\bar\Q$, by $d\bar\P^{\alpha,\eps}=(\bar\kappa_T^{\bar\nu^{\alpha,\eps}})^{-1}d\bar\Q$. Notice that the restriction of $\bar\Q$ to $\bar\Fc_T^{\bar x_0,\bar W,\bar\pi,\bar\theta^{\alpha,\eps}}$ coincides with $\bar\P^{\bar\nu^{\alpha,\eps}}$, which is the probability measure on $(\bar\Omega,\bar\Fc_T^{\bar x_0,\bar W,\bar\pi,\bar\theta^{\alpha,\eps}})$ defined as $d\bar\P^{\bar\nu^{\alpha,\eps}}=\bar\kappa_T^{\bar\nu^{\alpha,\eps}}d\bar\P^{\alpha,\eps}$. This proves item (i).

By Girsanov's theorem, under $\bar\P^{\alpha,\eps}$ the random measure $\bar\theta^{\alpha,\eps}$ has $\bar\F^{\bar x_0,\bar W,\bar\pi,\bar\theta^{\alpha,\eps}}$-compensator on $[0,T]\times\Lambda$ given by $\lambda_0(da)dt$, so in particular it is a Poisson random measure. Moreover, under $\bar\P^{\alpha,\eps}$ the random variable $\bar x_0$ has still the same law, the process $\bar W$ is still a Brownian motion, and the random measure $\bar\pi$ is still a Poisson random measure with $\bar\F^{\bar x_0,\bar W,\bar\pi,\bar\theta^{\alpha,\eps}}$-compensator on $[0,T]\times U$ given by $\lambda_\pi(dz)dt$. In addition, $\bar x_0$, $\bar W$, $\bar\pi$, $\bar\theta$ are independent under $\bar\P^{\alpha,\eps}$. This shows the validity of item (iii) and concludes the proof.
\ep

\begin{Theorem}\label{Main_thm}
Under assumptions {\bf (A)}-{\bf (A${}_{\mathcal R}$)}, the following equality holds:
\[
V_0 \ = \ \hat V_0^{\mathcal R}.
\]
\end{Theorem}
\textbf{Proof.}
\emph{Proof of the inequality $V_0\geq\hat V_0^{\mathcal R}$.} Let $(\bar\Omega,\bar\Fc,\bar\P;\bar x_0,\bar W,\bar\pi,\bar\theta;\bar I,\bar X;\bar\Vc)$ and $\bar\Gc$ be respectively the probabilistic setting for the randomized control problem and the $\sigma$-algebra mentioned in Lemma \ref{L:hatJR=hatJ}. Recall from Proposition \ref{P:V_0_G} that
\[
V_0 \ = \ \sup_{\bar\alpha\in\bar\Ac^{\bar\Gc}}\bar J(\bar\alpha).
\]
Then, the inequality $V_0\geq\hat V_0^{\mathcal R}$ follows directly by Lemma \ref{L:NewSetting} and Remark \ref{R:hatJR=hatJ}, from which we have
\[
\hat V_0^{\mathcal R} \ = \ \bar V_0^{\mathcal R} \ = \ \sup_{\bar\alpha\in\bar\Ac^{\bar\Vc}}\bar J(\bar\alpha) \ \leq \ \sup_{\bar\alpha\in\bar\Ac^{\bar\Gc}}\bar J(\bar\alpha) \ = \ V_0.
\]

\vspace{2mm}

\noindent\emph{Proof of the inequality $V_0\leq\hat V_0^{\mathcal R}$.} Fix $\alpha\in\Ac$. Then, for every positive integer $k$, it follows from Lemma \ref{L:Fuhrman} with $\eps=1/k$ that there exist a probabilistic setting for the randomized control problem $(\bar\Omega,\bar\Fc,\bar\P^{\alpha,k};\bar x_0,\bar W,\bar\pi,\bar\theta^{\alpha,k};\bar I^{\alpha,k},\bar X^{\alpha,k};\bar\Vc^{\alpha,k})$, a probability measure $\bar\Q$  on $(\bar\Omega,\bar\Fc)$  equivalent to $\bar\P^{\alpha,k}$, $\bar\alpha\colon[0,T]\times\bar\Omega\rightarrow\Lambda$,      $\bar\nu^{\alpha,k}\in\bar\Vc^{\alpha,k}$ such that:
\begin{enumerate}
\item[\textup{(i)}] $\bar\Q\text{\tiny$|\bar\Fc_T^{\bar x_0,\bar W,\bar\pi,\bar\theta^{\alpha,k}}$}$ coincides with $\bar\P^{\bar\nu^{\alpha,k}}$;
\item[\textup{(ii)}] $\bar\E^{\bar\Q}\big[\int_0^T d_\Lambda(\bar I_t^{\alpha,k},\bar\alpha_t)\,dt\big]\leq1/k$, so, in particular,
\begin{equation}\label{Krylov_alpha_k}
\bar\E^{\bar\Q}\bigg[\int_0^T d_\Lambda(\bar I_t^{\alpha,k},\bar\alpha_t)\,dt\bigg] \ \overset{k\rightarrow+\infty}{\longrightarrow} \ 0;
\end{equation}
\item[\textup{(iii)}] $(x_0,W, \pi,\alpha)$ under $\P$ has the same law as $(\bar x_0,\bar W, \bar\pi,\bar\alpha)$ under $\bar\P^{\alpha,k}$.
\end{enumerate}
The claim follows if we prove that
\begin{equation}\label{FirstIneq}
\lim_{k\rightarrow+\infty} \bar J^{\mathcal R,\alpha,k}(\bar\nu^{\alpha,k}) \ = \ J(\alpha),
\end{equation}
where $\bar J^{\mathcal R,\alpha,k}$ denotes the gain functional for the randomized control problem $(\bar\Omega,\bar\Fc,\bar\P^{\alpha,k};\bar x_0,$ $\bar W,\bar\pi,$ $\bar\theta^{\alpha,k};\bar I^{\alpha,k},\bar X^{\alpha,k};\bar\Vc^{\alpha,k})$, which is given by
\[
\bar J^{\mathcal R,\alpha,k}(\bar\nu^{\alpha,k}) \ = \ \bar\E^{\bar \nu^{\alpha,k}}
\bigg[\int_0^T f_t(\bar X^{\alpha,k},\bar I_t^{\alpha,k})\,dt + g(\bar X^{\alpha,k})\bigg],
\]
with
\[
\begin{cases}
\vspace{1mm}\displaystyle d\bar X_t^{\alpha,k} = A\bar X_t^{\alpha,k}dt \!+\! b_t(\bar X^{\alpha,k},\bar I_t^{\alpha,k})dt \!+ \!\sigma_t(\bar X^{\alpha,k},\bar I_t^{\alpha,k})d\bar W_t \!+ \!\!\!\int_{U \setminus \{0\}}\!\!\!\!\!\!\!\gamma_t(\bar X^{\alpha,k},\bar I_{t-}^{\alpha,k}, z)\big(\bar\pi(dt\,dz)\!-\!\lambda_\pi(dz)dt\big), \\
\displaystyle \bar X_0^{\alpha,k} = \bar x_0.
\end{cases}
\]
As a matter of fact, if \eqref{FirstIneq} holds true then for every $\eps>0$ there exists $k_\eps$ such that $J(\alpha)\leq\bar J^{\mathcal R,\alpha,k}(\bar\nu^{\alpha,k})+\eps\leq\bar V_0^{\mathcal R,\alpha,k}+\eps$, for all $k\geq k_\eps$. By Lemma \ref{L:NewSetting} we know that $\bar V_0^{\mathcal R,\alpha,k}=\hat V_0^{\mathcal R}$, so the claim follows.

It remains to prove \eqref{FirstIneq}. By item (i) above we notice that $\bar J^{\mathcal R,\alpha,k}(\bar\nu^{\alpha,k})$ can be equivalently written in terms of $\bar\E^{\bar\Q}$:
\[
\bar J^{\mathcal R,\alpha,k}(\bar\nu^{\alpha,k}) \ = \ \bar\E^{\bar\Q}
\bigg[\int_0^T f_t(\bar X^{\alpha,k},\bar I_t^{\alpha,k})\,dt + g(\bar X^{\alpha,k})\bigg].
\]
On the other hand, by item (iii) above, $J(\alpha)$ is also given by
\[
J(\alpha) \ = \ \bar\E^{\bar\Q}
\bigg[\int_0^T f_t(\bar X^{\bar\alpha},\bar\alpha_t)\,dt + g(\bar X^{\bar\alpha})\bigg],
\]
with
\[
\begin{cases}
\vspace{1mm}\displaystyle d\bar X_t^{\bar\alpha} = A\bar X_t^{\bar\alpha}dt + b_t(\bar X^{\bar\alpha},\bar\alpha_t)dt + \sigma_t(\bar X^{\bar\alpha},\bar\alpha_t)d\bar W_t + \!\int_{U \setminus \{0\}}\!\!\!\gamma_t(\bar X^{\bar\alpha},\bar\alpha_t, z)\big(\bar\pi(dt\,dz)-\lambda_\pi(dz)dt\big), \\
\displaystyle \bar X_0^{\bar\alpha} = \bar x_0.
\end{cases}
\]
Hence, \eqref{FirstIneq} can be equivalently rewritten as follows:
\begin{equation}\label{FirstIneq2}
\bar\E^{\bar\Q}
\bigg[\int_0^T f_t(\bar X^{\alpha,k},\bar I_t^{\alpha,k})\,dt + g(\bar X^{\alpha,k})\bigg] \ \overset{k\rightarrow+\infty}{\longrightarrow} \ \bar\E^{\bar\Q}
\bigg[\int_0^T f_t(\bar X^{\bar\alpha},\bar\alpha_t)\,dt + g(\bar X^{\bar\alpha})\bigg].
\end{equation}
Now, we notice that, under assumptions {\bf (A)}-{\bf (A${}_{\mathcal R}$)}, proceeding along the same lines as in the proof of Proposition \ref{P:SDE}, we can prove the following result: for every $1\leq p\leq p_0$,
\begin{equation}\label{Est_X-X_rand}
\bar\E^{\bar\Q}\Big[\sup_{t\in[0,T]}\big|\bar X_t^{\alpha,k} - \bar X_t^{\bar\alpha}\big|^p\Big] \ \overset{k\rightarrow+\infty}{\longrightarrow} \ 0.
\end{equation}
It is then easy to see that, from the continuity and polynomial growth assumptions on $f$ and $g$ in {\bf (A)}-(v) and {\bf (A)}-(vi), convergence \eqref{FirstIneq2} follows directly from \eqref{Krylov_alpha_k} and \eqref{Est_X-X_rand}. This concludes the proof of the inequality $V_0\leq\hat V_0^{\mathcal R}$.
\ep

\section{BSDE with non-positive jumps}
\label{S:BSDE}

Let $(\hat \Omega, \hat{\mathcal F}, \hat \P)$ be the
complete probability space on which are defined $\hat x_0$, $\hat W$, $\hat \pi$, $\hat\theta$ as in Section \ref{SubS:RandomizedSetting}. $\hat\F^{\hat x_0,\hat W,\hat\pi,\hat\theta}= (\hat{\mathcal F}_t^{\hat x_0,\hat W,\hat \pi,\hat\theta})_{t \geq 0}$ still denotes the $\hat\P$-completion of the filtration generated by $\hat x_0$, $\hat W$, $\hat\pi$, $\hat\theta$; we also recall that $\Pc(\hat\F^{\hat x_0, \hat W,\hat\pi,\hat\theta})$ is the predictable $\sigma$-algebra on $[0,T]\times\hat\Omega$ corresponding to $\hat\F^{\hat x_0,\hat W,\hat\pi,\hat\theta}$. We begin introducing the following notations.
\begin{itemize}
\item ${\bf S^2}$ denotes the set of c\`adl\`ag $\hat\F^{\hat x_0,\hat W,\hat\pi,\hat\theta}$-adapted processes $Y\colon[0,T]\times\hat\Omega\rightarrow\R$ satisfying
\[
\|Y\|_{_{{\bf S^2}}}^2 \ := \ \hat\E\Big[ \sup_{0\leq t\leq T} |Y_t|^2 \Big] \ < \ \infty.
\]
\item ${\bf L^p(0,T)}$, $p$ $\geq$ $1$, denotes the set of $\hat\F^{\hat x_0, \hat W,\hat\pi,\hat\theta}$-adapted processes $\phi\colon[0,T]\times\hat\Omega\rightarrow\R$ satisfying
\[
\|\phi\|_{_{{\bf L^p(0,T)}}}^p \ := \ \hat\E\bigg[\int_0^T |\phi_t|^p\,dt\bigg] \ < \ \infty.
\]
\item ${\bf L^p(\hat W)}$, $p$ $\geq$ $1$, denotes the set of $\Pc(\hat\F^{\hat x_0, \hat W,\hat\pi,\hat\theta})$-measurable processes
$Z\colon[0,T]\times\hat\Omega\rightarrow\Xi$ satisfying
\[
\|Z\|_{_{\bf L^p(\hat W)}}^p \ := \ \hat\E\bigg[\bigg(\int_0^T |Z_t|_\Xi^2\,dt\bigg)^{\frac{p}{2}}\bigg] \ < \ \infty.
\]
We shall identify $\Xi$ with its dual $\Xi^*$. Notice also that $\Xi^*=L_2(\Xi,\R)$, the space of Hilbert-Schmidt operators from $\Xi$ into $\R$ endowed with the usual scalar product.
\item ${\bf L^p(\hat\pi)}$, $p$ $\geq$ $1$, denotes the set of $\Pc(\hat\F^{\hat x_0,\hat W,\hat\pi,\hat\theta})\otimes\Bc(U)$-measurable maps $L\colon[0,T]\times\hat\Omega\times U\rightarrow\R$ satisfying
\[
\| L \|_{_{{\bf L^p(\hat\pi)}}}^p \ := \ \hat\E\bigg[ \bigg(\int_0^T\int_U |L_t(z)|^2\, \lambda_\pi(dz)\,dt\bigg)^{\frac{p}{2}} \bigg] \ < \ \infty.
\]
\item ${\bf L^p(\hat\theta)}$, $p$ $\geq$ $1$, denotes the set of $\Pc(\hat\F^{\hat x_0,\hat W,\hat\pi,\hat\theta})\otimes\Bc(\Lambda)$-measurable maps $R\colon[0,T]\times\hat\Omega\times \Lambda\rightarrow\R$ satisfying
\[
\| R \|_{_{{\bf L^p(\hat\theta)}}}^p \ := \ \hat\E\bigg[ \bigg(\int_0^T\int_\Lambda  |R_t(b)|^2\, \lambda_0(db)\,dt\bigg)^{\frac{p}{2}} \bigg] \ < \ \infty.
\]
\item ${\bf L^p(\lambda_0)}$, $p$ $\geq$ $1$, denotes the set of $\Bc(\Lambda)$-measurable maps $r\colon\Lambda\rightarrow\R$ satisfying
\[
\| r \|_{_{{\bf L^p(\lambda_0)}}}^p \ := \ \int_\Lambda  |r(b)|^p\, \lambda_0(db) \ < \ \infty.
\]
\item ${\bf K^2}$ denotes the set of non-decreasing $\Pc(\hat\F^{\hat x_0,\hat W,\hat\pi,\hat\theta})$-measurable processes $K\in{\bf S^2}$ satisfying $K_0=0$, so that
\[
\|K\|_{_{\bf S^2}}^2 \ = \ \hat\E|K_T|^2.
\]
\end{itemize}
Consider the following backward stochastic differential equation with non-positive jumps:\begin{align}
Y_t \ &= \ g(\hat X) + \int_t^T f(\hat X,\hat I_s)\,ds + K_T - K_t - \int_t^T\int_\Lambda R_s(b)\,\hat \theta(ds\,db) \label{BSDE} \\
&\quad \ -\int_t^T Z_s\,d\hat W_s - \int_t^T\int_U L_s(z)\,(\hat\pi(ds\,dz)-\lambda_\pi(dz)\,ds),\qquad 0 \leq t \leq T,\,\,\hat\P\text{-a.s.} \notag \\
R_t(b) \ &\leq \ 0, \qquad dt\otimes d\hat\P\otimes\lambda_0(db)\text{-a.e. on }\hat\Omega\times[0,T]\times  \Lambda. \label{Constraint}
\end{align}

\begin{Definition}
A \textbf{minimal solution} to equation \eqref{BSDE}-\eqref{Constraint} is a quintuple $(Y,Z,L,$ $R,K)\in{\bf S^2}\times{\bf L^2(\hat W)}\times{\bf L^2(\hat\pi)}\times{\bf L^2(\hat\theta)}\times{\bf K^2}$
satisfying \eqref{BSDE}-\eqref{Constraint} such that for any other quintuple $(\tilde Y,\tilde Z,\tilde L,\tilde R,\tilde K)\in{\bf S^2}\times{\bf L^2(\hat W)}\times{\bf L^2(\hat\pi)}\times{\bf L^2(\hat\theta)}\times{\bf K^2}$ satisfying \eqref{BSDE}-\eqref{Constraint}, we have
\[
Y_t \ \leq \ \tilde Y_t,\quad  0 \leq t \leq T, \,\,\hat\P\textup{-a.s.}
\]

\end{Definition}

\begin{Lemma}
\label{L:Uniqueness}
Under assumptions {\bf (A)}-{\bf (A${}_{\mathcal R}$)}, there exists at most one minimal solution to equation \eqref{BSDE}-\eqref{Constraint}.
\end{Lemma}
\textbf{Proof.}
The uniqueness of $Y$ follows from the definition of minimal solution.
Now, let $(Y,Z,L,R,K)$, $(Y,\tilde Z,\tilde L,\tilde R,\tilde K)\in{\bf S^2}\times{\bf L^2(\hat W)}\times{\bf L^2(\hat\pi)}\times{\bf L^2(\hat\theta)}\times{\bf K^2}$ be two minimal solutions. Then
\begin{align}\label{BSDEdifference}
&K_t - \tilde K_t - \int_0^t \big(Z_s - \tilde Z_s\big)\,d\hat W_s + \int_0^t\int_U \big(L_s(z) - \tilde L_s(z)\big)\,\lambda_\pi(dz)ds\notag\\
&= \int_0^t\int_U \big(L_s(z) - \tilde L_s(z)\big)\,\hat\pi(ds\,dz) +  \int_0^t\int_\Lambda \big(R_s(b) - \tilde R_s(b)\big)\,\hat\theta(ds\,db),
\end{align}
for all $0\leq t\leq T$, $\hat\P$-a.s..
Observe that on the left-hand side of \eqref{BSDEdifference} there is a predictable process, which has therefore no totally inaccessible jumps, while on the right-hand side in \eqref{BSDEdifference} there is a pure jump process which has only totally inaccessible jumps. We deduce that both sides must be equal to zero. Therefore, we obtain the two following equalities: for all $0\leq t\leq T$, $\hat\P$-a.s.,
\begin{align*}
K_t - \tilde K_t + \int_0^t\int_U \big(L_s(z) - \tilde L_s(z)\big)\,\lambda_\pi(dz)ds &= \int_0^t \big(Z_s - \tilde Z_s\big)\,d\hat W_s,\\
 \int_0^t\int_U \big(L_s(z) - \tilde L_s(z)\big)\,\hat\pi(ds\,dz)
&=   \int_0^t\int_\Lambda \big(R_s(b) - \tilde R_s(b)\big)\,\hat\theta(ds\,db).
\end{align*}
Concerning the first equation,
 the left-hand side is a finite variation process, while the process on the right-hand side
has not finite variation, unless $Z = \tilde Z$ and $K - \tilde K + \int_0^\cdot\int_U (L_s(z) - \tilde L_s(z))\lambda_\pi(dz)ds=0$. On the other hand,  since $\hat \pi$ and $\hat \theta$ are independent, they have disjoint jump times, therefore from the second equation above
we find $L=\tilde L$ and  $R=\tilde R$, from which we also obtain $K=\tilde K$.
\ep

\vspace{3mm}

We now prove that focus on the existence of a minimal solution to \eqref{BSDE}-\eqref{Constraint}. To this end, we introduce, for every integer $n\geq1$, the following penalized backward stochastic differential equation:
\begin{align}
\label{BSDE_pen}
Y_t^n \ &= \ g( \hat X) + \int_t^T f(\hat X,\hat I_s)\,ds + K_T^n - K_t^n - \int_t^T\int_\Lambda R_s^n(b)\,\hat\theta(ds\,db)\\
&\quad \ - \int_t^T Z_s^n\,d\hat W_s - \int_t^T\int_U L_s^n(z)\,(\hat\pi(ds\,dz)-\lambda_\pi(dz) ds), \qquad 0 \leq t \leq T,\;\hat\P\text{-a.s.} \notag
\end{align}
where
\[
K_t^n \ = \ n\int_0^t\int_\Lambda \big(R_s^n(b)\big)^+\, \lambda_0(db)ds, \qquad 0 \leq t \leq T,\;\hat\P\text{-a.s.}
\]
with $f^+= \max (f,0)$ denoting the positive part of the function $f$.

\begin{Lemma}[Martingale representation]
\label{L:MartRepr}
Suppose that assumptions {\bf (A)}\textup{-(iii)} and {\bf (A${}_{\mathcal R}$)}\textup{-(i)} hold. Given any $\xi\in{\bf L^2}(\hat\Omega,\hat\Fc_T^{\hat x_0,\hat W,\hat\pi,\hat\theta},\hat\P)$, there exist $Z\in{\bf L^2(\hat W)}$, $L\in{\bf L^2(\hat\pi)}$, $R\in{\bf L^2(\hat\theta)}$ such that
\begin{equation}\label{Repr_xi}
\xi \ = \ \hat\E[\xi|\hat x_0] + \int_0^T Z_t\,d\hat W_t + \int_0^T \int_U L_t(z)\,\hat\pi(dt\,dz) + \int_0^T \int_\Lambda R_t(b)\,\hat\theta(dt\,db), \qquad \hat\P\text{-a.s.}
\end{equation}
\end{Lemma}
\textbf{Proof.}
We begin noting that, when $\hat W$ is a finite-dimensional Brownian motion, representation \eqref{Repr_xi} for $\xi$ can be easily proved using for instance Lemma 2.3 in \cite{tangli94}. As a matter of fact, let $\hat\F^{\hat x_0}=(\hat\Fc_t^{\hat x_0})_{t\geq0}$, $\hat\F^{\hat W_0}=(\hat\Fc_t^{\hat W_0})_{t\geq0}$, $\hat\F^{\hat\pi_0}=(\hat\Fc_t^{\hat\pi_0})_{t\geq0}$, $\hat\F^{\hat\theta_0}=(\hat\Fc_t^{\hat\theta_0})_{t\geq0}$ be the $\hat\P$-completion of the filtration generated respectively by $\hat x_0$, $\hat W$, $\hat\pi$, $\hat\theta$. When $\xi=1_{\text{\tiny$E_{\hat x_0}$}}1_{\text{\tiny$E_{\hat W_0}$}}1_{\text{\tiny$E_{\hat\pi_0}$}}1_{\text{\tiny$E_{\hat\theta_0}$}}$, with $E_{\hat x_0}\in\hat\Fc_T^{\hat x_0}$, $E_{\hat W_0}\in\hat\Fc_T^{\hat W_0}$, $E_{\hat\pi_0}\in\hat\Fc_T^{\hat\pi_0}$, $E_{\hat\theta_0}\in\hat\Fc_T^{\hat\theta_0}$, then representation \eqref{Repr_xi} for $\xi$ follows easily by Lemma 2.3 in \cite{tangli94}. Since the linear span of the random variables of the form $1_{\text{\tiny$E_{\hat x_0}$}}1_{\text{\tiny$E_{\hat W_0}$}}1_{\text{\tiny$E_{\hat\pi_0}$}}1_{\text{\tiny$E_{\hat\theta_0}$}}$ is dense in ${\bf L^2}(\hat\Omega,\hat\Fc_T^{\hat x_0,\hat W,\hat\pi,\hat\theta},\hat\P)$, we deduce the validity of \eqref{Repr_xi} for a general $\xi\in{\bf L^2}(\hat\Omega,\hat\Fc_T^{\hat x_0,\hat W,\hat\pi,\hat\theta},\hat\P)$.

In the infinite-dimensional case, let $(e_k)_{k\geq1}$ be an orthonormal basis of $\Xi$ and define $\hat W_t^{(k)} = \langle \hat W_t,e_k\rangle_\Xi$, for $t\geq0$. The processes $W^{(k)}$ are independent standard real Brownian motions. For any positive integer $n$, let $\hat\F^{(n)}=(\hat\Fc_t^{(n)})_{t\geq0}$ denote the $\hat\P$-completion of the filtration generated by $\hat x_0$, $\hat W^{(1)},\ldots,\hat W^{(n)}$, $\hat\pi$, $\hat\theta$. Notice that $\hat\F^{(n)}$ satisfies the usual conditions. Denote $\xi^{(n)} = \hat\E[\xi|\hat\Fc_T^{(n)}]$. By the previously mentioned finite-dimensional version of representation \eqref{Repr_xi}, we have a martingale representation for $\xi^{(n)}$. It is then easy to see that, letting $n\rightarrow+\infty$ in such a martingale representation, \eqref{Repr_xi} follows.
\ep

\begin{Proposition}
\label{P:BSDEpen}
Under assumptions {\bf (A)}-{\bf (A${}_{\mathcal R}$)}, for every integer $n\geq1$ there exists a unique solution $(Y^n,Z^n,L^n, R^n)\in{\bf S^2}\times{\bf L^2(\hat W)}\times{\bf L^2(\hat\pi)} \times{\bf L^2(\hat\theta)}$ to equation \eqref{BSDE_pen}. In addition, the following estimate holds:
\begin{align}\label{EstimateZLRK_n}
\|Z^n\|_{_{\bf L^2(\hat W)}}^2 \!\!+ \|L^n\|_{_{\bf L^2(\hat\pi)}}^2 \!\!+ \|R^n\|_{_{\bf L^2(\hat\theta)}}^2 \!\!+ \|K^n\|_{_{\bf S^2}}^2 \leq \hat C\bigg(\|Y^n\|_{_{\bf S^2}}^2 + \hat\E\bigg[\int_0^T \!\!|f(\hat X,\hat I_t)|^2dt\bigg]\bigg),
\end{align}
for some constant $\hat C\geq0$, depending only on $T$ and on the constant $L$ in assumption {\bf (A)}\textup{-(vi)}, independent of $n$.
\end{Proposition}
\textbf{Proof.}
The existence and uniqueness result can be proved as in the finite-dimensional case $\text{dim}\,\Xi<\infty$, see Lemma 2.4 in \cite{tangli94}. We simply recall that, as usual, it is based on a fixed point argument and on the martingale representation (concerning this latter result, since we did not find a reference for it suitable for our setting, we proved it in Lemma \ref{L:MartRepr}).

Similarly, estimate \eqref{EstimateZLRK_n} can be proved proceeding along the same lines as in the finite-dimensional case $\text{dim}\,\Xi<\infty$, for which we refer to Lemma 2.3 in \cite{khapha12}; we just recall that its proof is based on the application of It\^o's formula to $|Y^n|^2$, as well as on Gronwall's lemma and the Burkholder-Davis-Gundy inequality.
\ep

\vspace{3mm}

For every integer $n\geq1$, we provide the following representation of $Y^n$ in terms of a suitable penalized randomized control problem. To this end, we define $\hat\Vc_n$ as the subset of $\hat\Vc$ of all maps $\hat\nu$ bounded from above by $n$.

We recall that, for every $\hat\nu \in \hat \Vc$, $\hat\E^{\hat\nu}$ denotes the expectation with respect to the probability measure on $(\hat\Omega,\hat\Fc_T^{\hat x_0, \hat W,\hat\pi, \hat \theta})$ given by $d\hat\P^{\hat \nu}=\hat\kappa_T^{\hat \nu}\,d\hat\P$, where $(\hat\kappa_t^{\hat\nu})_{t\in[0,T]}$ denotes the Dol\'eans-Dade exponential defined in (\ref{Doleans}).

\begin{Lemma}\label{L:Y_n_formula}
Under assumptions {\bf (A)}-{\bf (A${}_{\mathcal R}$)}, for every integer $n\geq1$ the following equalities hold:
\begin{equation}\label{Y_pen_formula}
Y_t^n \ = \ \esssup_{\hat\nu\in\hat\Vc_n} \hat\E^{\hat\nu}
\bigg[\int_t^T f(\hat X,\hat I_s)\,ds + g(\hat X)\bigg|\hat\Fc_t^{\hat x_0, \hat W,\hat\pi, \hat \theta}\bigg], \qquad \hat\P\text{-a.s.},\,0\leq t\leq T
\end{equation}
and
\begin{equation}\label{Y_pen_formula_0}
\hat\E[Y_0^n] \ = \ \sup_{\hat\nu\in\hat\Vc_n} \hat\E^{\hat\nu}
\bigg[\int_0^T f(\hat X,\hat I_s)\,ds + g(\hat X)\bigg],
\end{equation}
with $\hat\E[Y_0^n]=Y_0^n$, $\hat\P$-a.s., when $\hat x_0$ is deterministic. In addition, we have:
\begin{itemize}
\item for every $0\leq t\leq T$, the sequence $(Y_t^n)_n$ is non-decreasing;
\item there exists a constant $\bar C\geq0$, depending only on $T$, $\bar p$, and on the constant $L$ in assumption {\bf (A)}\textup{-(vi)}, independent of $n$, such that
\begin{align}
\label{EstimateYn}
\sup_{s \in [0,\,T]} |Y_s^n| \ \leq \ \bar C\Big(1 + \sup_{s \in [0,\,T]} |\hat X_s|^{\bar p}\Big), \qquad \hat\P\text{-a.s.}
\end{align}
\end{itemize}
\end{Lemma}
\textbf{Proof.}
\emph{Proof of formulae \eqref{Y_pen_formula} and \eqref{Y_pen_formula_0}.} We report the proof of formula \eqref{Y_pen_formula}, as \eqref{Y_pen_formula_0} can be proved proceeding along the same lines (simply replacing all the $\hat\P^{\hat\nu}$-conditional expectations with normal $\hat\P^{\hat\nu}$-expectations, and also noting that $\hat\P^{\hat\nu}$ coincides with $\hat\P$ on $\hat\Fc_0^{\hat x_0,\hat W,\hat\pi,\hat\theta}$, which is the $\hat\P$-completion of the $\sigma$-algebra generated by $\hat x_0$). Fix an integer $n \geq 1$ and let $(Y^n,Z^n,L^n,R^n)$ be the solution to \eqref{BSDE_pen}, whose existence follows from Proposition \ref{P:BSDEpen}. As consequence of the Girsanov Theorem, the two following processes
$$\int_0^t Z_s^n\,d\hat W_s, \qquad\qquad  \int_0^t\int_U L_s^n(z)\,\big(\hat\pi(ds\,dz)-\lambda_\pi(dz) ds\big),$$
are $\hat \P^{\hat\nu}$-martingales (see e.g. Theorem 15.3.10 in \cite{Coh-Ell} or Theorem 12.31 in \cite{HeWangYan}). Moreover
$$ \hat\E^{\hat\nu} \bigg[\int_t^T\int_\Lambda R_s^n(b)\,\hat\theta(ds\,db)\bigg| {\hat\Fc_t^{\hat x_0, \hat W,\hat\pi, \hat \theta}}\bigg]=\hat\E^{\hat\nu} \bigg[\int_t^T\int_\Lambda R_s^n(b)\,{\hat\nu}_s(b)\lambda_0(db)ds\bigg| {\hat\Fc_t^{\hat x_0, \hat W,\hat\pi, \hat \theta}}\bigg].$$
Therefore, taking the $\hat \P^{\hat\nu}$-conditional expectation given $\hat\Fc_t^{\hat x_0, \hat W,\hat\pi, \hat \theta}$ in \eqref{BSDE_pen},
we
obtain
\begin{align}\label{eqYncond}
Y_t^n \ &= \ \hat\E^{\hat\nu} \bigg[g( \hat X_T) + \int_t^T f(\hat X_s,\hat I_s)\,ds \bigg| {\hat\Fc_t^{\hat x_0, \hat W,\hat\pi, \hat \theta}}\bigg] \\
&\quad \ + \hat\E^{\hat\nu} \bigg[\int_t^T\int_\Lambda [n(R_s^n(b))^+-{\hat\nu}_s(b)\,R_s^n(b)]\,\lambda_0(db)ds \bigg| {\hat\Fc_t^{\hat x_0, \hat W,\hat\pi, \hat \theta}}\bigg], \qquad \hat\P\text{-a.s.},\,0 \leq t \leq T.\notag
\end{align}
Firstly, we notice that $n u^+ - \nu u \geq 0$ for all $u \in \R$, $\nu \in (0, n]$, so that \eqref{eqYncond} gives
\begin{align}\label{Ynineq1}
Y_t^n \ &\geq \ \esssup_{\hat\nu\in\hat\Vc_n}\hat\E^{\hat\nu} \bigg[g( \hat X_T) + \int_t^T f(\hat X_s,\hat I_s)\,ds \bigg| {\hat\Fc_t^{\hat x_0, \hat W,\hat\pi, \hat \theta}}\bigg] \qquad \hat\P\text{-a.s.},\,0 \leq t \leq T.
\end{align}
On the other hand, since $R^n \in {\bf L^2(\hat\theta)}$, by Lebesgue's dominated convergence theorem for conditional expectation, we obtain
$$
\lim_{N \rightarrow \infty}\hat\E \bigg[\int_t^T\int_\Lambda |R_s^n(b)|^2\,1_{\{R_s^n(b) \leq -N\}}\,\lambda_0(db)ds \bigg| {\hat\Fc_t^{\hat x_0, \hat W,\hat\pi, \hat \theta}}\bigg] \ = \ 0.
$$
So, in particular, for every $n\geq1$ there exists a positive integer $N_n$ such that
\begin{equation}\label{Lebesgue}
\hat\E \bigg[\int_t^T\int_\Lambda |R_s^n(b)|^2\,1_{\{R_s^n(b) \leq -N_n\}}\,\lambda_0(db)ds \bigg] \ \leq \ e^{-(n-1)\lambda_0(\Lambda)(T-t)}.
\end{equation}
Now, let us define
$$
\hat\nu_s^{n,\varepsilon}(b) \ := \ n 1_{\{R_s^n(b) \geq 0\}} + \varepsilon 1_{\{-1 < R_s^n(b) <0\}} - \varepsilon R_s^n(b)^{-1}1_{\{- N_n < R_s^n(b) \leq -1\}} + \varepsilon 1_{\{R_s^n(b) \leq - N_n\}}.
$$
It is easy to see  that $\hat\nu^{n,\varepsilon} \in\hat\Vc_n$. Moreover, we have
\begin{align}\label{ineqfin2}
&\hat\E^{\hat\nu^{n,\varepsilon}}\bigg[\int_t^T\int_\Lambda [n(R_s^n(b))^+-\hat\nu_s^{n,\varepsilon}(b)\,R_s^n(b)]\,\lambda_0(db)ds \bigg| {\hat\Fc_t^{\hat x_0,\hat W,\hat\pi, \hat \theta}}\bigg] \ \leq \ \varepsilon \sqrt{(T-t) \lambda_0(\Lambda)}\Bigg\{\sqrt{(T-t) \lambda_0(\Lambda)} \notag \\
&+\sqrt{\hat\E\bigg[\Big|\frac{\hat\kappa^{\hat\nu^{n,\varepsilon}}_T}{\hat\kappa^{\hat\nu^{n,\varepsilon}}_t}\Big|^2\bigg|{\hat\Fc_t^{\hat x_0, \hat W,\hat\pi, \hat \theta}}\bigg]} \sqrt{\hat\E \bigg[\int_t^T\int_\Lambda |R_s^n(b)|^2\,1_{\{R_s^n(b) \leq -N_n\}}\,\lambda_0(db)ds \bigg| {\hat\Fc_t^{\hat x_0, \hat W,\hat\pi, \hat \theta}}\bigg] }\Bigg\}.
\end{align}
Recalling that, for every $\hat\nu\in\hat\Vc$, it holds that $|\hat\kappa^{\hat\nu}_s|^2 = \hat\kappa_s^{\hat\nu^2} e^{\int_0^s \int_{\Lambda} (\hat\nu_r(b)-1)\lambda_0(db) dr}$, $s \in [0,\,T]$ (see e.g. the proof of Lemma 4.1 in \cite{khapha12}), we obtain
\begin{align}\label{kappa/kappa}
\hat\E\bigg[\Big|\frac{\hat\kappa^{\hat\nu^{n,\varepsilon}}_T}{\hat\kappa^{\hat\nu^{n,\varepsilon}}_t}\Big|^2\bigg|{\hat\Fc_t^{\hat x_0, \hat W,\hat\pi, \hat \theta}}\bigg] \ &= \ \hat\E\bigg[\frac{\hat\kappa^{|\hat\nu^{n,\varepsilon}|^2}_T}{\hat\kappa^{|\hat\nu^{n,\varepsilon}|^2}_t}e^{\int_t^T \int_{\Lambda} (\hat\nu_r^{n,\varepsilon}(b)-1)\lambda_0(db) dr}\bigg|{\hat\Fc_t^{\hat x_0, \hat W,\hat\pi, \hat \theta}}\bigg] \\
&\leq \ \hat\E\bigg[\frac{\hat\kappa^{|\hat\nu^{n,\varepsilon}|^2}_T}{\hat\kappa^{|\hat\nu^{n,\varepsilon}|^2}_t}e^{(n-1)\lambda_0(\Lambda)(T-t)}\bigg|{\hat\Fc_t^{\hat x_0, \hat W,\hat\pi, \hat \theta}}\bigg] \ = \ e^{(n-1)\lambda_0(\Lambda)(T-t)}, \notag
\end{align}
where the last equality follows from the fact that, for every $\hat\nu\in\hat\Vc$, we have $\hat\nu^2\in\hat\Vc$, so that $\hat\kappa^{\hat\nu^2}$ is a martingale. Plugging \eqref{Lebesgue} and \eqref{kappa/kappa} into \eqref{ineqfin2}, we end up with
\begin{align}\label{ineqfin3}
\hat\E^{\hat\nu^{n,\varepsilon}} \bigg[\int_t^T\int_\Lambda [n(R_s^n(b))^+-\hat\nu^{n,\varepsilon}_s(b)\,R_s^n(b)]\,\lambda_0(db)ds \bigg| {\hat\Fc_t^{\hat x_0,\hat W,\hat\pi, \hat \theta}}\bigg]&\\
\leq \ \varepsilon \sqrt{(T-t) \lambda_0(\Lambda)} \Big\{ \sqrt{(T-t) \lambda_0(\Lambda)}+1\Big\}& \ = \ \varepsilon\,\tilde C,\notag
\end{align}
with $\tilde C:=\sqrt{(T-t) \lambda_0(\Lambda)}\{ \sqrt{(T-t) \lambda_0(\Lambda)}+1\}$. Plugging \eqref{ineqfin3} into \eqref{eqYncond} we get
\begin{align*}
Y_t^n \ &\leq \ \hat\E^{\hat\nu^{n,\varepsilon}}\bigg[g( \hat X_T) + \int_t^T f(\hat X_s,\hat I_s)\,ds \bigg| {\hat\Fc_t^{\hat x_0, \hat W,\hat\pi, \hat \theta}}\bigg] +\varepsilon\,\tilde C\notag\\
&\leq \ \esssup_{\hat\nu\in\hat\Vc_n} \hat\E^{\hat\nu}\bigg[g( \hat X_T) + \int_t^T f(\hat X_s,\hat I_s)\,ds \bigg| {\hat\Fc_t^{\hat x_0, \hat W,\hat\pi, \hat \theta}}\bigg] +\varepsilon\,\tilde C, \quad \hat\P\text{-a.s.},\,0 \leq t \leq T.
\end{align*}
From the arbitrariness of $\eps$, we find the reverse inequality of \eqref{Ynineq1}, from which \eqref{Y_pen_formula} follows.

\medskip

\noindent \emph{Proof of the monotonicity of $(Y^n)_n$.} By definition $\hat\Vc_n \subset
\hat\Vc_{n+1}$. Then inequality $Y^n_t \leq  Y^{n+1}_t$, $\hat \P$-a.s. for all $t \in [0, T]$, follows directly from \eqref{Y_pen_formula}.

\medskip

\noindent \emph{Proof of formula \eqref{EstimateYn}.} In the sequel we denote by $\bar C$ a non-negative constant, depending only on $T$, $\bar p$, and on the constant $L$ in assumption {\bf (A)}\textup{-(vi)}, independent of $n$, which may change from line to line.

Recalling the polynomial growth condition \eqref{PolynomialGrowth_f_g} on $f$ and $g$ in assumption {\bf (A)}-(vi), it follows from formula \eqref{Y_pen_formula} that
$$
|Y^n_t | \ \leq \ \bar C \,\esssup_{\hat\nu\in\hat\Vc_n} \hat\E^{\hat\nu}\Big[1 + \sup_{s \in[0,T]}|\hat X_s|^{\bar p} \Big| {\hat\Fc_t^{\hat x_0, \hat W,\hat\pi, \hat \theta}}\Big], \qquad \hat\P\textup{-a.s.},\,0 \leq t \leq T.
$$
Finally, by estimate \eqref{Est_X_rand_conditional}, together with the fact that $Y^n$ is a c\`{a}dl\`{a}g process, we see that \eqref{EstimateYn} follows.
\ep

\vspace{3mm}

We can now prove the main result of this section.

\begin{Theorem}
\label{ThmExistence}
Under assumptions {\bf (A)}-{\bf (A${}_{\mathcal R}$)}, there exists a unique minimal solution $(Y,Z,$ $L,R, K)\in{\bf S^2}\times{\bf L^2(\hat W)}\times{\bf L^2(\hat\pi)}\times{\bf L^2(\hat\theta)}\times{\bf K^2}$ to \eqref{BSDE}-\eqref{Constraint}, satisfying
\begin{equation}\label{Y_formula}
Y_t \ = \ \esssup_{\hat\nu\in\hat\Vc} \hat\E^{\hat\nu}
\bigg[\int_t^T f_s(\hat X,\hat I_s)\,ds + g(\hat X)\bigg|\hat\Fc_t^{\hat x_0, \hat W,\hat\pi, \hat \theta}\bigg], \qquad \hat\P\text{-a.s.},\,0\leq t\leq T
\end{equation}
and
\begin{equation}\label{Y_formula_0}
\hat\E[Y_0] \ = \ \sup_{\hat\nu\in\hat\Vc} \hat\E^{\hat\nu}
\bigg[\int_0^T f_s(\hat X,\hat I_s)\,ds + g(\hat X)\bigg] \ = \ \hat V_0^{\mathcal R},
\end{equation}
with $\hat\E[Y_0]=Y_0$, $\hat\P$-a.s., when $\hat x_0$ is deterministic. In addition, we have:
\begin{enumerate}
\item[\textup{(i)}] for every $0\leq t\leq T$, the sequence $(Y_t^n)_n$ increasingly converges to $Y_t$; moreover, $Y^n\rightarrow Y$ in ${\bf L^2(0,T)}$;
\item[\textup{(ii)}] the following estimate holds:
\begin{align}
\label{EstimateY}
\sup_{s \in [0,\,T]} |Y_s| \ \leq \ \bar C\Big(1 + \sup_{s \in [0,\,T]} |\hat X_s|^{\bar p}\Big), \qquad \hat\P\text{-a.s.},
\end{align}
with the same constant $\bar C$ as in \eqref{EstimateYn};
\item[\textup{(iii)}] the sequence $(Z^n,L^n,R^n)_n$ weakly converges to $(Z,L,R)$ in ${\bf L^2(\hat W)}\times{\bf L^2(\hat \pi)}\times{\bf L^2(\hat\theta)}$;
\item[\textup{(iv)}] for every $0\leq t\leq T$, the sequence $(K_t^n)_n$ weakly converges to $K_t$ in ${\bf L^2}(\hat\Omega,\hat \Fc_t^{\hat x_0,\hat W,\hat\pi,\hat\theta},\hat\P)$.
\end{enumerate}
Finally, the so-called randomized dynamic programming principle holds: for every $t\in[0,T]$ and any $\hat\F^{\hat x_0,\hat W,\hat\pi,\hat\theta}$-stopping time $\hat\tau$ taking values in $[t,T]$, we have
\begin{equation}\label{RandDynProgPr}
Y_t \ = \ \esssup_{\hat\nu\in\hat\Vc} \hat\E^{\hat\nu}
\bigg[\int_t^{\hat\tau} f_s(\hat X,\hat I_s)\,ds + Y_{\hat\tau}\bigg|\hat\Fc_t^{\hat x_0, \hat W,\hat\pi, \hat \theta}\bigg], \qquad \hat\P\text{-a.s.}
\end{equation}
\end{Theorem}
\textbf{Proof.} \emph{Construction of $(Y,Z,L,R,K)$ in ${\bf S^2}\times{\bf L^2(\hat W)}\times{\bf L^2(\hat\pi)}\times{\bf L^2(\hat\theta)}\times{\bf K^2}$ solution to \eqref{BSDE}.} By Lemma \ref{L:Y_n_formula} we know that, for every $0\leq t\leq T$, the sequence $(Y_t^n)_n$ is non-decreasing. Since $Y^n$ is c\`{a}dl\`{a}g, it follows that there exists a $\hat\P$-null set $\hat N$ such that, for every integer $n\geq1$,
\[
Y_t^n(\hat\omega) \ \leq \ Y_t^{n+1}(\hat\omega), \qquad 0\leq t\leq T,\;\hat\omega\notin\hat N.
\]
This property, together with estimate \eqref{EstimateYn}, shows that there exists a measurable $\hat\F^{\hat x_0,\hat W,\hat\pi,\hat\theta}$-adapted process $Y=(Y_t)_{t\geq0}$ such that $Y_t^n(\hat\omega)$ increasingly converges to $Y_t(\hat\omega)$, $0\leq t\leq T$, $\hat\omega\notin\hat N$. Moreover, estimate \eqref{EstimateY} holds, from which we also deduce that $Y^n\rightarrow Y$ in ${\bf L^2(0,T)}$. In addition, noting that $\hat\Vc_n\subset\hat\Vc_{n+1}$ and $\cup_n\hat\Vc_n=\hat\Vc$, letting $n\rightarrow\infty$ in equalities \eqref{Y_pen_formula} and \eqref{Y_pen_formula_0}, we obtain formulae \eqref{Y_formula} and \eqref{Y_formula_0}, respectively.

By estimate \eqref{EstimateZLRK_n}, we see that the sequence $(Z^n,L^n,R^n)_n$ is bounded in the Hilbert space ${\bf L^2(\hat W)}\times{\bf L^2(\hat\pi)}\times{\bf L^2(\hat\theta)}$. So, in particular, $(Z^n,L^n,R^n)_n$ admits a weakly convergent subsequence $(Z^{n_k},L^{n_k},R^{n_k})_k$ going towards some $(Z,L,R)\in{\bf L^2(\hat W)}\times{\bf L^2(\hat\pi)}\times{\bf L^2(\hat\theta)}$. Then, for any $\hat\F^{\hat x_0,\hat W,\hat\pi,\hat\theta}$-stopping time $\hat\tau$ taking values in $[0,T]$, we obtain
\begin{align*}
&\int_0^{\hat\tau} Z_s^{n_k}\,d\hat W_s \ \rightharpoonup \ \int_0^{\hat\tau} Z_s\,d\hat W_s, \qquad\qquad \int_0^{\hat\tau}\int_\Lambda R_s^{n_k}(b)\,\hat\theta(ds\,db) \ \rightharpoonup \ \int_0^{\hat\tau}\int_\Lambda R_s(b)\,\hat\theta(ds\,db), \\
&\int_0^{\hat\tau}\int_U L_s^{n_k}(z)\,(\hat\pi(ds\,dz)-\lambda_\pi(dz) ds) \ \rightharpoonup \ \int_0^{\hat\tau}\int_U L_s(z)\,(\hat\pi(ds\,dz)-\lambda_\pi(dz) ds).
\end{align*}
By equation \eqref{BSDE_pen}, we have
\begin{align*}
K_{\hat\tau}^n \ &= \ Y_{\hat\tau}^n - Y_0^n - g(\hat X) - \int_0^{\hat\tau} f_s(\hat X,\hat I_s)\,ds + \int_0^{\hat\tau}\int_\Lambda R_s^n(b)\,\hat\theta(ds\,db) \\
&\quad \ + \int_0^{\hat\tau} Z_s^n\,d\hat W_s + \int_0^{\hat\tau}\int_U L_s^n(z)\,(\hat\pi(ds\,dz)-\lambda_\pi(dz) ds).
\end{align*}
Noting that $Y_{\hat\tau}^n\rightarrow Y_{\hat\tau}$ strongly in ${\bf L^2}(\hat\Omega,\hat \Fc_{\hat\tau}^{\hat x_0,\hat W,\hat\pi,\hat\theta},\hat\P)$, we get
\[
K_{\hat\tau}^{n_k} \ \rightharpoonup \ K_{\hat\tau},
\]
where
\begin{align*}
K_t \ &:= \ Y_t - Y_0 - g(\hat X) - \int_0^t f_s(\hat X,\hat I_s)\,ds + \int_0^t\int_\Lambda R_s(b)\,\hat\theta(ds\,db) \\
&\quad \ + \int_0^t Z_s\,d\hat W_s + \int_0^t\int_U L_s(z)\,(\hat\pi(ds\,dz)-\lambda_\pi(dz) ds), \qquad\qquad 0\leq t\leq T.
\end{align*}
Since $K_T^{n_k}\rightharpoonup K_T$, from the lower semicontinuity of the norm with respect to the weak topology on ${\bf L^2}(\hat\Omega,\hat \Fc_T^{\hat x_0,\hat W,\hat\pi,\hat\theta},\hat\P)$, we deduce that $\hat\E|K_T|^2<\infty$. It is also easy to see that $K^{n_k}$ weakly converges to $K$ in ${\bf L^2(0,T)}$. Since the set of $\hat\F^{\hat x_0,\hat W,\hat\pi,\hat\theta}$-predictable processes is convex and strongly closed in ${\bf L^2(0,T)}$, it is also weakly closed, so that $K$ is $\hat\F^{\hat x_0,\hat W,\hat\pi,\hat\theta}$-predictable.

Now, given any $\hat\F^{\hat x_0,\hat W,\hat\pi,\hat\theta}$-stopping times $\hat\tau$ and $\hat\tau'$, with $0\leq\hat\tau\leq\hat\tau'\leq T$, since $K_{\hat\tau}^n\leq K_{\hat\tau'}^n$, $\hat\P$-a.s., we deduce that $K_{\hat\tau}\leq K_{\hat\tau'}$, $\hat\P$-a.s.. This implies that $K$ is a non-decreasing process. As a matter fact, $K$ is non-decreasing if and only if the two processes $K$ and $\sup_{0\leq s\leq \cdot}K_s$ are $\hat\P$-indistinguishable. Since $K$ is predictable, we notice that $\sup_{0\leq s\leq \cdot}K_s$ is also predictable (by the proof of item (a) of Theorem IV.33 in \cite{DellacherieMeyer} we know that $\sup_{0\leq s< \cdot}K_s$ is progressively measurable and left-continuous, hence it is predictable; since $K$ is predictable and $\sup_{0\leq s\leq \cdot}K_s=K_\cdot\vee\sup_{0\leq s< \cdot}K_s$, we deduce that $\sup_{0\leq s\leq \cdot}K_s$ is predictable). Let
\[
\hat\tau \ = \ \inf\Big\{t\geq0\colon K_t\,<\,\sup_{0\leq s\leq t}K_s\Big\}, \qquad \hat\tau' \ = \ \inf\Big\{t\geq\hat\tau\colon K_t\,=\,\sup_{0\leq s\leq t}K_s\Big\},
\]
with $\inf\emptyset=\infty$. The claim follows if we prove that $\hat\P(\hat\tau<\infty)=0$. We proceed by contradiction, assuming that $E:=\{\hat\tau<\infty\}$ is such that $\hat\P(E)>0$. We begin noting that $\hat\tau<\hat\tau'$ on $E$. Now, for every $\hat\omega\in E$ and any $t$ satisfying $\hat\tau(\hat\omega)<t<\hat\tau'(\hat\omega)$, we obtain
\begin{equation}\label{contradiction}
K_t(\hat\omega) \ < \ \sup_{0\leq s\leq t} K_s(\hat\omega) \ = \ \sup_{0\leq s\leq\hat\tau(\hat\omega)} K_s(\hat\omega).
\end{equation}
Since $K$ and $\sup_{0\leq s\leq \cdot}K_s$ are predictable, $\hat\tau$ (resp. $\hat\tau'$) is a predictable time, so, in particular, there exists a sequence of stopping times $\hat\tau_m\uparrow\hat\tau$, with $\hat\tau_m<\hat\tau_{m+1}<\hat\tau$ whenever $\hat\tau\neq0$ (resp. $\hat\tau_m'\uparrow\hat\tau'$, with $\hat\tau_m'<\hat\tau_{m+1}'<\hat\tau'$ whenever $\hat\tau'\neq0$). It is then easy to prove (using that $\hat\tau<\hat\tau'$ on $E$ and $\hat\tau'$ is announceable) the existence of a stopping time $\bar\tau$ satisfying $\hat\tau<\bar\tau<\hat\tau'$ on $E$. Moreover, using that $\hat\tau$ is announceable, we obtain $K_{\hat\tau}=\sup_{0\leq s\leq\hat\tau}K_s$, arguing as follows. Let $F:=\{K_{\hat\tau}<\sup_{0\leq s\leq\hat\tau}K_s\}\cap E$. On $F$ it holds that $\sup_{0\leq s\leq\hat\tau}K_s=\sup_{0\leq s<\hat\tau}K_s$. Since $\hat\tau_m\uparrow\hat\tau$ and the stochastic process $\sup_{0\leq s<\cdot}K_s$ is left-continuous, we have $\sup_{0\leq s<\hat\tau_m}K_s\uparrow\sup_{0\leq s<\hat\tau}K_s$. As $\hat\tau_m<\hat\tau_{m+1}$ on $E$, it follows that $\sup_{0\leq s<\hat\tau_m}K_s\leq\sup_{0\leq s\leq\hat\tau_m}K_s\leq\sup_{0\leq s<\hat\tau_{m+1}}K_s$ on $E$, therefore $K_{\hat\tau_m}=\sup_{0\leq s\leq\hat\tau_m}K_s\uparrow\sup_{0\leq s<\hat\tau}K_s$ on $E$. Recalling that $\sup_{0\leq s<\hat\tau}K_s>K_{\hat\tau}$ on $F$, we get a contradiction with $K_{\hat\tau_m}\leq K_{\hat\tau}$, unless $F$ is a $\hat\P$-null set. Finally, from \eqref{contradiction} with $t=\bar\tau(\hat\omega)$, we obtain
\[
K_{\bar\tau(\hat\omega)} \ < \ K_{\hat\tau(\hat\omega)}, \qquad \text{for every }\hat\omega\in E\backslash F,
\]
which is in contradiction with $K_{\bar\tau}\geq K_{\hat\tau}$, unless $E$ is a $\hat\P$-null set. This shows that $\hat\P(\hat\tau<\infty)=\hat\P(E)=0$ and proves that $K$ is a non-decreasing process. Finally, by Lemma 2.2 in \cite{Peng99} it follows that both $Y$ and $K$ are c\`{a}dl\`{a}g, so, in particular, they belong to ${\bf S^2}$. We conclude that $(Y,Z,L,R, K)\in{\bf S^2}\times{\bf L^2(\hat W)}\times{\bf L^2(\hat\pi)}\times{\bf L^2(\hat\theta)}\times{\bf K^2}$ is a solution to equation \eqref{BSDE}.

Proceeding along the same lines as in the proof of Lemma \ref{L:Uniqueness}, we deduce that given $Y$ there exists a unique quadruple $(Z,L,R,K)$ in ${\bf L^2(\hat W)}\times{\bf L^2(\hat\pi)}\times{\bf L^2(\hat\theta)}\times{\bf K^2}$ satisfying equation \eqref{BSDE}. It follows that the entire sequence $(Z^n,L^n,R^n)_n$ weakly converges to $(Z,L,R)$ in ${\bf L^2(\hat W)}\times{\bf L^2(\hat\pi)}\times{\bf L^2(\hat\theta)}$, so that item (iii) holds. Similarly, item (iv) holds.

\vspace{1mm}

\noindent\emph{Jump constraint \eqref{Constraint}.} Let $\Phi\colon{\bf L^2(\hat\theta)}\rightarrow\R$ be given by
\[
\Phi(\tilde R) \ = \ \hat \E\bigg[\int_0^T\int_\Lambda (\tilde R_t(a))^+\,\lambda_0(db)dt\bigg]^2, \qquad \forall\,\tilde R\in{\bf L^2(\hat\theta)}.
\]
Since $\Phi$ is convex and strongly continuous, it is also weakly lower-semicontinuous, therefore
\[
\Phi(R) \ \leq \ \liminf_{n\rightarrow\infty} \Phi(R^n) \ = \ \liminf_{n\rightarrow\infty}\frac{\hat\E|K_T^n|^2}{n^2} \ = \ 0,
\]
where the last equality follows from estimates \eqref{EstimateZLRK_n} and \eqref{EstimateYn}. This implies that $\Phi(R)=0$, that is
\[
\hat \E\bigg[\int_0^T\int_\Lambda (R_t(a))^+ \lambda_0(db)dt\bigg]^2 \ = \ 0,
\]
which means that the jump constraint \eqref{Constraint} is satisfied. In conclusion, $(Y,Z,L,R,K)$ is a solution to \eqref{BSDE}-\eqref{Constraint}.

\vspace{1mm}

\noindent\emph{Proof of the minimality of $(Y,Z,L,R,K)$.} The minimality follows from $Y=\lim_n Y_n$. In fact, let $(\tilde Y,\tilde Z,\tilde L,\tilde R, \tilde K)\in{\bf S^2}\times{\bf L^2(W)}\times{\bf L^2(\tilde\pi)}\times{\bf L^2(\tilde\theta)}\times{\bf K^2}$ be another solution to \eqref{BSDE}-\eqref{Constraint}. Proceeding as in the proof of formula \eqref{Y_pen_formula} (see the beginning of the proof of Lemma \ref{L:Y_n_formula}), given any $t\in[0,T]$ and $\hat\nu\in\hat\Vc$, taking the $\hat \P^{\hat\nu}$-conditional expectation with respect to $\hat\Fc_t^{\hat x_0, \hat W,\hat\pi, \hat \theta}$ in \eqref{BSDE}, we obtain, $\hat\P$-a.s.,
\begin{align*}
\tilde Y_t \ &= \ \hat\E^{\hat\nu} \bigg[g( \hat X_T) + \int_t^T f_s(\hat X,\hat I_s)\,ds \bigg| {\hat\Fc_t^{\hat x_0, \hat W,\hat\pi, \hat \theta}}\bigg] + \hat\E^{\hat\nu}\big[\tilde K_T - \tilde K_t\big|{\hat\Fc_t^{\hat x_0, \hat W,\hat\pi, \hat \theta}}\big] \\
&\quad \ - \hat\E^{\hat\nu} \bigg[\int_t^T\int_\Lambda {\hat\nu}_s(b)\,R_s(b)\,\lambda_0(db)ds \bigg| {\hat\Fc_t^{\hat x_0, \hat W,\hat\pi, \hat \theta}}\bigg] \ \geq \ \hat\E^{\hat\nu} \bigg[g( \hat X_T) + \int_t^T f_s(\hat X,\hat I_s)\,ds \bigg| {\hat\Fc_t^{\hat x_0, \hat W,\hat\pi, \hat \theta}}\bigg].
\end{align*}
From the arbitrariness of $\hat\nu$, we get
\[
\tilde Y_t \ \geq \ \esssup_{\hat\nu\in\hat\Vc}\hat\E^{\hat\nu} \bigg[g( \hat X_T) + \int_t^T f_s(\hat X,\hat I_s)\,ds \bigg| {\hat\Fc_t^{\hat x_0, \hat W,\hat\pi, \hat \theta}}\bigg] \qquad \hat\P\text{-a.s.},\,0 \leq t \leq T.
\]
By formula \eqref{Y_pen_formula}, recalling that $\hat\Vc_n\subset\hat\Vc$, we conclude that $Y_t^n \leq \tilde Y_t$, $0 \leq t \leq T$, $\hat\P$-a.s.. Letting $n\rightarrow\infty$, we obtain $Y_t \leq \tilde Y_t$, $0 \leq t \leq T$, $\P$-a.s., which proves the minimality of $(Y,Z,L,R,K)$. Finally, by Proposition \ref{L:Uniqueness} we know that $(Y,Z,L,R,K)$ is unique.

\vspace{1mm}

\noindent\emph{Proof of the randomized dynamic programming principle \eqref{RandDynProgPr}.} Fix $t\in[0,T]$ and let $\hat\tau$ be a $\hat\F^{\hat x_0,\hat W,\hat\pi,\hat\theta}$-stopping time taking values in $[t,T]$. Given any integer $n\geq1$, consider the penalized equation \eqref{BSDE_pen} between $0$ and $\hat\tau$ with terminal condition $Y_{\hat\tau}^n$. Then, proceeding along the same lines as in the proof of formula \eqref{Y_pen_formula}, we obtain
\[
Y_t^n \ = \ \esssup_{\hat\nu\in\hat\Vc_n} \hat\E^{\hat\nu}
\bigg[\int_t^{\hat\tau} f_s(\hat X,\hat I_s)\,ds + Y_{\hat\tau}^n\bigg|\hat\Fc_t^{\hat x_0, \hat W,\hat\pi, \hat \theta}\bigg], \qquad \hat\P\text{-a.s.}
\]
Recalling that $\hat\Vc_n\subset\hat\Vc$ and $Y^n\leq Y$, we find $Y_t^n$ $\leq$ $\esssup_{\hat\nu\in\hat\Vc} \hat\E^{\hat\nu}
[\int_t^{\hat\tau} f_s(\hat X,\hat I_s)\,ds + Y_{\hat\tau}|\hat\Fc_t^{\hat x_0, \hat W,\hat\pi, \hat \theta}]$. Letting $n\rightarrow\infty$, we conclude that
\[
Y_t \ \leq \ \esssup_{\hat\nu\in\hat\Vc} \hat\E^{\hat\nu}
\bigg[\int_t^{\hat\tau} f_s(\hat X,\hat I_s)\,ds + Y_{\hat\tau}\bigg|\hat\Fc_t^{\hat x_0, \hat W,\hat\pi, \hat \theta}\bigg], \qquad \hat\P\text{-a.s.}
\]
In order to prove the reverse inequality, take a positive integer $m$, then, for every $n\geq m$,
\begin{align*}
Y_t &\geq \esssup_{\hat\nu\in\hat\Vc_n} \hat\E^{\hat\nu}
\bigg[\int_t^{\hat\tau} f_s(\hat X,\hat I_s)\,ds + Y_{\hat\tau}^n\bigg|\hat\Fc_t^{\hat x_0, \hat W,\hat\pi, \hat \theta}\bigg] \geq \esssup_{\hat\nu\in\hat\Vc_n} \hat\E^{\hat\nu}
\bigg[\int_t^{\hat\tau} f_s(\hat X,\hat I_s)\,ds + Y_{\hat\tau}^m\bigg|\hat\Fc_t^{\hat x_0, \hat W,\hat\pi, \hat \theta}\bigg],
\end{align*}
where we have used that $Y_t\geq Y_t^n$ and $Y_{\hat\tau}^n\geq Y_{\hat\tau}^m$. From the arbitrariness of $n$, we end up with $Y_t$ $\geq$ $\hat\E^{\hat\nu}[\int_t^{\hat\tau} f_s(\hat X,\hat I_s)\,ds + Y_{\hat\tau}^m|\hat\Fc_t^{\hat x_0, \hat W,\hat\pi, \hat \theta}]$, for any $\hat\nu\in\hat\Vc$ and $m\geq1$. Letting $m\rightarrow\infty$ and taking the essential supremum over $\hat\Vc$, we see that the claim follows.
\ep

\section{HJB equation in Hilbert spaces: the Markovian case}
\label{S:HJB}

\setcounter{equation}{0} \setcounter{Assumption}{0}
\setcounter{Theorem}{0} \setcounter{Proposition}{0}
\setcounter{Corollary}{0} \setcounter{Lemma}{0}
\setcounter{Definition}{0} \setcounter{Remark}{0}

In the present section, we replace assumptions {\bf (A)} by the set of assumptions {\bf (A$_{\text{\textnormal{\textsc M}}}$)} reported below. Before stating {\bf (A$_{\text{\textnormal{\textsc M}}}$)}, we notice that in this section, $A$ still denotes a linear operator from $\Dc(A)\subset H$ into $H$, while the coefficients $b$, $\sigma$, $\gamma$, $f$, $g$ are \emph{non}-path-depedent, namely $b\colon[0,T]\times H\times \Lambda\rightarrow H$, $\sigma\colon[0,T]\times H\times \Lambda\rightarrow L(\Xi;H)$, $\gamma: [0,T] \times H \times \Lambda \times U \rightarrow H$, $f\colon[0,T]\times H\times \Lambda\rightarrow\R$, $g\colon H\rightarrow\R$. In what follows, we shall impose the following assumptions on $A$, $b$, $\sigma$, $\gamma$, $f$, $g$.

\vspace{3mm}

\noindent {\bf (A$_{\text{\textnormal{\textsc M}}}$)}
\begin{itemize}
\item [(i)] $A$ is a linear, densely defined, maximal dissipative operator in $H$. In particular, $A$ is the generator of a strongly continuous semigroup $\{e^{tA},\ t\geq 0\}$ of contractions.
Moreover, there exists (see e.g. Theorem 3.11 in \cite{fabbrigozziswiech17}) an operator $B\colon H\rightarrow H$, which is linear, bounded, strictly positive, self-adjoint, with $A^* B$ bounded on $H$, such that the \emph{weak $B$-condition} for $A$ holds
\[
\langle(-A^*B + c_0 B)x,x\rangle \geq 0, \qquad \text{for all }x\in H,
\]
for some constant $c_0 \geq 0$.\\
\emph{We define on $H$ the norm $|\cdot|_{-1}$, defined as
$|x|_{-1} := \big|B^{1/2}x\big|$, for every $x\in H$. In addition, we define the space $H_{-1}$ to be the completion of $H$ under the norm $|\cdot|_{-1}$. $H_{-1}$ is a Hilbert space equipped with the scalar product
\[
\langle x,y\rangle_{-1} := \big\langle B^{1/2}x,B^{1/2}y\big\rangle.
\]}
\item[(ii)] There exists a Borel measurable function $\rho\colon U\rightarrow\R$, bounded on bounded subsets of $U$, such that
\[
\inf_{|z|_{U}> R} \rho(z) \ > \ 0, \quad \text{for every }R \ > \ 0 \qquad\quad \text{ and } \qquad\quad
\int_U |\rho(z)|^2 \lambda_\pi(dz) \ < \ \infty.
\]
\item [(iii)] The maps $b$, $\gamma$, $f$, $g$ are Borel measurable. For every $v\in H$, the map $\sigma(\cdot,\cdot,\cdot)v\colon[0,T]\times H\times \Lambda\rightarrow H$ is Borel measurable.
\item [(iv)] The map $g$ is continuous on $H$ with respect to the supremum norm. For every $t\in[0,T]$, the maps $b(t,\cdot,\cdot)$ and $f(t,\cdot,\cdot)$ are continuous on $H\times \Lambda$. For every $(t,z)\in[0,T]\times U$, the map $\gamma(t,\cdot,\cdot,z)$ is continuous on $H\times \Lambda$. For every $t\in[0,T]$ and any $s\in(0,T]$, we have $e^{sA}\sigma(t,x,a)\in L_2(\Xi;H)$, for all $(x,a)\in H\times \Lambda$, and the map $e^{sA}\sigma(t,\cdot,\cdot)\colon H\times \Lambda\rightarrow L_2(\Xi;H)$ is continuous.
\item [(v)] For all $t\in[0,T]$, $s\in(0,T]$, $x,x'\in H$, $a\in \Lambda$, $z\in U$,
\begin{align*}
|b(t,x,a) - b(t,x',a)| + |e^{sA}\sigma(t,x,a)-e^{sA}\sigma(t,x',a)|_{L_2(\Xi;H)}
\ & \leq \  L   |x-x'|_{-1}, \notag \\
|\gamma(t,x,a,z)-\gamma(t,x',a,z)| \ &\leq \ L\,\rho(z) |x-x'|_{-1}, \notag \\
|b(t,0,a)| + |\sigma(t,0,a)|_{L_2(\Xi;H)}
 \ & \leq \ L, \notag
\\
|\gamma(t,0,a,z)| \ & \leq \ L\,\rho(z), \notag
\\
|f(t,x,a) - f(t,x',a)| + |g(x) - g(x')| \ & \leq \ \omega(|x-x'|_{-1}), \\
|f(t,0,a)| \ &\leq \ L,
\end{align*}
for some constant $L\geq0$ and some \emph{modulus of continuity} $\omega$, i.e. a continuous, non-decreasing, subadditive map $\omega\colon[0,\infty)\rightarrow[0,\infty)$ satisfying $\omega(0)=0$ and $\omega(r)>0$, for any $r>0$.
\end{itemize}

\paragraph{Stochastic optimal control problem.}

We now formulate the stochastic optimal control problem in such a setting. Since the formulation can be done proceeding along the same lines as in subsection \ref{SubS:Standard}, we focus on the main steps. We consider a complete probability space $(\Omega,\Fc,\P)$ on which are defined a cylindrical Brownian motion $W=(W_t)_{t\geq0}$, with values in $\Xi$, and an independent Poisson random measure  $\pi(dt \,dz)$ on $[0,\,\infty) \times U$ with compensator $\lambda_{\pi}(dz)\,dt$. For every $t\geq0$, we denote by $\F^{t,W,\pi}= (\mathcal F_s^{t,W,\pi})_{s \geq t}$ the $\P$-completion of the filtration generated by $(W_s-W_t)_{s\geq t}$ and the restriction of $\pi(dt\,dz)$ to $[t,\infty)\times U$.

For every $t\in[0,T]$, an admissible control process at time $t$ will be any $\F^{t,W,\pi}$-predictable process $\alpha\colon[t,T]\times\Omega\rightarrow\Lambda$. For every $t\in[0,T]$, the set of all admissible control processes at time $t$ is denoted by $\Ac_t$. For every $(t,x)\in[0,T]\times H$ and any $\alpha\in\Ac_t$, the controlled equation has the form
\begin{equation}\label{SDE_Markov}
\begin{cases}
\vspace{1mm}\displaystyle dX_s \ = \ AX_s\,ds + b(s,X_s,\alpha_s)\,ds + \sigma(s,X_s,\alpha_s)\,dW_s  \\
\vspace{1mm}\displaystyle\qquad\quad\;\; + \, \int_{U \setminus \{0\}} \gamma(s,X_s, \alpha_s, z)\,\big(\pi(ds\,dz)-\lambda_\pi(dz)\,ds\big), \qquad\qquad t\leq s\leq T,\\
\displaystyle X_t \ = \ x.
\end{cases}
\end{equation}
We have the following result.
\begin{Proposition}
Under assumption {\bf (A$_{\text{\textnormal{\textsc M}}}$)}, for every $(t,x)\in[0,T]\times H$ and any $\alpha \in \mathcal A_t$, there exists a unique mild solution $X^{t,x,\alpha}=(X_s^{t,x,\alpha})_{s\in[t,T]}$ to equation \eqref{SDE_Markov}. Moreover, for every $p\geq1$,
\begin{equation}\label{EstXalpha_Markov}
\E\Big[\sup_{s\in[t,T]}|X_s^{t,x,\alpha}|^p\Big] \ \leq \ C_p\,\big(1 + |x|^p\big),
\end{equation}
for some positive constant $C_p$, independent of $t$, $x$, $\alpha$.
\end{Proposition}
\textbf{Proof.}
The proof can be done proceeding along the same lines as in the proof of Proposition \ref{SDE}.
\ep

\vspace{3mm}

\noindent The controller aims at maximizing over all $\alpha \in \mathcal A_t$ the gain functional
\[
J(t,x,\alpha) \ = \ \E\bigg[\int_t^T f(s,X_s^{t,x,\alpha},\alpha_s)\,ds + g(X_T^{t,x,\alpha})\bigg].
\]
Finally, the value function of the stochastic control problem is given by
\begin{equation}\label{value}
v(t,x) \ = \ \sup_{\alpha\in\mathcal A_t} J(t,x,\alpha), \qquad (t,x)\in[0,T]\times H.
\end{equation}

\begin{Lemma}
Let assumption {\bf (A$_{\text{\textnormal{\textsc M}}}$)} hold. There exist a modulus of continuity $\omega_v$ and a constant $C\geq0$ such that
\begin{align}
|J(t,x,\alpha) - J(t,x',\alpha)| \ &\leq \ \omega_v(|x-x'|_{-1}), \label{J_unif} \\
|J(t,x,\alpha)| \ &\leq \ C\,\big(1+|x|_{-1}\big), \label{J_pol}
\end{align}
for all $t\in[0,T]$, $x,x'\in H$, $\alpha\in\Ac_t$. In particular,
\begin{align}
|v(t,x) - v(t,x')| \ &\leq \ \omega_v(|x-x'|_{-1}), \label{v_unif} \\
|v(t,x)| \ &\leq \ C\,\big(1+|x|_{-1}\big), \notag
\end{align}
for all $t\in[0,T]$, $x,x'\in H$.
\end{Lemma}
\textbf{Proof.}
We begin noting that, proceeding along the same lines as in the proof of estimate (3.12) of Theorem 3.4 in \cite{swiechzabczyk16}, we can prove that the following estimate holds:
\begin{equation}\label{Estimate-1}
\sup_{t\leq s\leq T}\E\big[|X_s^{t,x,\alpha} - X_s^{t,x',\alpha}|_{-1}^2\big] \ \leq \ \bar C\,|x - x'|_{-1},
\end{equation}
for some constant $\bar C\geq0$, independent of $t$, $x$, $x'$, $\alpha$. Then, \eqref{J_unif} follows directly from estimate \eqref{Estimate-1} and the assumptions on $f$ and $g$ in {\bf (A$_{\text{\textnormal{\textsc M}}}$)}-(v). On the other hand, \eqref{J_pol} follows from estimate \eqref{EstXalpha_Markov}, using again the assumptions on $f$ and $g$ in {\bf (A$_{\text{\textnormal{\textsc M}}}$)}-(v).
\ep

\paragraph{Randomized setting.}

We now consider, following Section \ref{S:Randomized}, the randomized setting. We focus on the main steps. We consider a complete probability space $(\hat\Omega,\hat\Fc,\hat\P)$ on which are defined a cylindrical Brownian motion $\hat W=(\hat W_t)_{t\geq0}$ with values in $\Xi$, a Poisson random measure  $\hat\pi(dt \,dz)$ on $[0,\,\infty) \times U$ with compensator $\lambda_{\pi}(dz)\,dt$, and a Poisson random measure  $\hat\theta(dt \,da)$ on $[0,\,\infty) \times\Lambda$ with compensator $\lambda_0(da)\,dt$ (satisfying assumption {\bf (A$_{\mathcal R}$)}\textup{-(i)}). For every $t\geq0$, we denote by $\hat\F^{t,\hat W,\hat\pi,\hat\theta}= (\hat{\mathcal F}_s^{t,\hat W,\hat\pi,\hat\theta})_{s \geq t}$ the $\hat\P$-completion of the filtration generated by $(\hat W_s-\hat W_t)_{s\geq t}$, the restriction of $\hat\pi(dt\,dz)$ to $[t,\infty)\times U$, the restriction of $\hat\theta(dt\,da)$ to $[t,\infty)\times\Lambda$. Finally, we denote by $\Pc(\hat\F^{t,\hat W,\hat\pi,\hat\theta})$ the predictable $\sigma$-algebra on $[t,T]\times\hat\Omega$ associated with $\hat\F^{t,\hat W,\hat\pi,\hat\theta}$.

For every $t\in[0,T]$, we denote by $\hat\Vc_t$ the set of all $\Pc(\hat\F^{t,\hat W,\hat\pi,\hat\theta})\otimes\Bc(\Lambda)$-measurable functions $\hat\nu\colon[t,T]\times\hat\Omega\times\Lambda\rightarrow(0,\infty)$ which are bounded from above and bounded away from zero. Given $\hat \nu\in\hat \Vc_t$, as in Section \ref{S:Randomized} we consider the corresponding Dol\'eans-Dade exponential $\hat\kappa^{t,\hat\nu}=(\hat\kappa_s^{t,\hat\nu})_{s\in[t,T]}$ defined as in \eqref{Doleans} and we introduce the probability measure $\hat\P^{t,\hat \nu}$ on $(\hat\Omega,\hat\Fc_T^{t,\hat W,\hat\pi, \hat \theta})$ as $d\hat\P^{t,\hat \nu}=\hat\kappa_T^{t,\hat \nu}\,d\hat\P$. Finally, we denote by $\hat\E^{t,\hat\nu}$ the expectation with respect to $\hat\P^{t,\hat\nu}$.

For every $t\in[0,T]$ and $a\in\Lambda$, we denote by $\hat I^{t,a}=(\hat I_s^{t,a})_{s\in[t,T]}$ the stochastic process taking values in $\Lambda$ defined as (notice that, when $\Lambda$ is a subset of a vector space, we can write \eqref{FSDEI} also as $\hat I_s^{t,a} = a + \int_{t}^{s}\int_\Lambda(b- \hat I_{r-}^{t,a}) \, \hat \theta(dr \,db)$, $s\in[t,T]$)
\begin{equation}\label{FSDEI}
\hat I_s^{t,a} \ = \ \sum_{n \geq 1} a\,1_{[t,\hat T_n)}(s) + \sum_{\substack{n \geq 1\\t<\hat T_n}} \hat \eta_n\,1_{[\hat T_n , \hat T_{n+1})}(s), \qquad\qquad \text{for all }t \leq s\leq T,
\end{equation}
where we recall that $(\hat T_n,\hat\eta_n)_{n\geq1}$ is the marked point process associated with the random measure $\hat \theta$, in particular we have $\hat \theta(dt\,da)=\sum_{n \geq 1} \delta_{(\hat T_n,\hat\eta_n)}(dt\,da)$.

Now, for every $(t,x,a)\in[0,T]\times H\times\Lambda$, we consider the following equation:
\begin{equation}\label{FSDEX}
\begin{cases}
\vspace{1mm}\displaystyle d\hat X_s \ = \ A\hat X_s\,ds + b(s,\hat X_s,\hat I_s)\,ds + \sigma(s,\hat X_s,\hat I_s)\,d\hat W_s  \\
\vspace{1mm}\displaystyle\qquad\quad\;\; + \, \int_{U \setminus \{0\}} \gamma(s,\hat X_s,\hat I_{s-}, z)\,\big(\hat\pi(ds\,dz)-\lambda_\pi(dz)\,ds\big), \qquad\qquad t\leq s\leq T,\\
\displaystyle \hat X_t \ = \ x.
\end{cases}
\end{equation}
We have the following result.
\begin{Proposition}
Under assumptions {\bf (A$_{\text{\textnormal{\textsc M}}}$)} and {\bf (A${}_{\mathcal R}$)}\textup{-(i)}, for every $(t,x,a)\in[0,T]\times H\times\Lambda$, there exists a unique mild solution $\hat X^{t,x,a}=(\hat X_s^{t,x,a})_{s\in[t,T]}$ to equation \eqref{FSDEX}, such that, for every $p\geq1$,
\begin{equation}\label{Est_X_rand_Markov}
\hat\E\Big[\sup_{s\in[t,T]}|\hat X_s^{t,x,a}|^p\Big] \ \leq \ C_p\,\big(1 + |x|^p\big),
\end{equation}
for some positive constant $C_p$, independent of $t$, $x$, $a$.
\end{Proposition}
\textbf{Proof.}
The proof can be done proceeding along the same lines as in the proof of Proposition \ref{P:SDE}.
\ep

\paragraph{BSDE with non-positive jumps.} We introduce the following additional notations.

\begin{itemize}
\item ${\bf S^2(t,T)}$ denotes the set of c\`adl\`ag $\hat\F^{t,\hat W,\hat\pi,\hat\theta}$-adapted processes $Y\colon[t,T]\times\hat\Omega\rightarrow\R$ satisfying
\[
\|Y\|_{_{{\bf S^2(t,T)}}}^2 \ := \ \hat\E\Big[ \sup_{t\leq s\leq T} |Y_s|^2 \Big] \ < \ \infty.
\]
\item ${\bf L^p(\hat W;t,T)}$, $p$ $\geq$ $1$, denotes the set of $\Pc(\hat\F^{t,\hat W,\hat\pi,\hat\theta})$-measurable processes
$Z\colon[t,T]\times\hat\Omega\rightarrow\Xi$ satisfying
\[
\|Z\|_{_{\bf L^p(\hat W)}}^p \ := \ \hat\E\bigg[\bigg(\int_t^T |Z_s|_\Xi^2\,ds\bigg)^{\frac{p}{2}}\bigg] \ < \ \infty.
\]
\item ${\bf L^p(\hat\pi;t,T)}$, $p$ $\geq$ $1$, denotes the set of $\Pc(\hat\F^{t,\hat W,\hat\pi,\hat\theta})\otimes\Bc(U)$-measurable maps $L\colon[t,T]\times\hat\Omega\times U\rightarrow\R$ satisfying
\[
\| L \|_{_{{\bf L^p(\hat\pi)}}}^p \ := \ \hat\E\bigg[ \bigg(\int_t^T\int_U |L_s(z)|^2\, \lambda_\pi(dz)\,ds\bigg)^{\frac{p}{2}} \bigg] \ < \ \infty.
\]
\item ${\bf L^p(\hat\theta;t,T)}$, $p$ $\geq$ $1$, denotes the set of $\Pc(\hat\F^{t,\hat W,\hat\pi,\hat\theta})\otimes\Bc(\Lambda)$-measurable maps $R\colon[t,T]\times\hat\Omega\times \Lambda\rightarrow\R$ satisfying
\[
\| R \|_{_{{\bf L^p(\hat\theta)}}}^p \ := \ \hat\E\bigg[ \bigg(\int_t^T\int_\Lambda  |R_s(b)|^2\, \lambda_0(db)\,ds\bigg)^{\frac{p}{2}} \bigg] \ < \ \infty.
\]
\item ${\bf K^2(t,T)}$ denotes the set of non-decreasing $\Pc(\hat\F^{t,\hat W,\hat\pi,\hat\theta})$-measurable processes $K\in{\bf S^2(t,T)}$ satisfying $K_t=0$.
\end{itemize}

For every $(t,x,a)\in[0,T]\times H\times\Lambda$, we introduce the following backward stochastic differential equation with non-positive jumps:
\begin{align}
Y_s \ &= \ g(\hat X_T^{t,x,a}) + \int_s^T f(r,\hat X_r^{t,x,a},\hat I_r^{t,a}) dr    +  K_T - K_s   -  \int_s^T\int_\Lambda R_r(b) \hat\theta(dr,db) \label{BSDEM}  \\
&\quad \   - \int_s^T  Z_r d\hat W_r - \int_s^T\int_{U\setminus \{0\}} L_r(z)\,(\hat\pi(dr\,dz)-\lambda_\pi(dz)\,dr), \;\;\;   t \leq s \leq T, \, \hat\P\text{-a.s.} \nonumber \\
R_s(b) \ &\leq \ 0, \;\;\;\;\; ds\otimes d\hat \P\otimes\lambda_0(db)\text{-a.e. on }[t,T]\times\hat\Omega\times \Lambda. \label{ConstraintM}
\end{align}

\begin{Definition}
Given $(t,x,a)\in[0,T]\times H\times\Lambda$, a \textbf{minimal solution} to equation \eqref{BSDEM}-\eqref{ConstraintM} is a quintuple $(Y,Z,L,R, K)\in{\bf S^2(t,T)}\times{\bf L^2(\hat W;t,T)}\times{\bf L^2(\hat\pi;t,T)}\times{\bf L^2(\hat\theta;t,T)}\times{\bf K^2(t,T)}$
satisfying \eqref{BSDEM}-\eqref{ConstraintM} such that for any other quintuple $(\tilde Y,\tilde Z,\tilde L,\tilde R,\tilde K)\in{\bf S^2(t,T)}$ $\times{\bf L^2(\hat W;t,T)}\times{\bf L^2(\hat\pi;t,T)}\times{\bf L^2(\hat\theta;t,T)}\times{\bf K^2(t,T)}$ satisfying \eqref{BSDEM}-\eqref{ConstraintM}, we have
\[
Y_s \ \leq \ \tilde Y_s,\quad  t \leq s \leq T, \,\,\hat\P\textup{-a.s.}
\]
\end{Definition}

We can now state the two main results of this section: the first result is the probabilistic representation formula (or non-linear Feynman-Kac formula) for the value function $v$ defined in \eqref{value}; the second result is the so-called \emph{randomized dynamic programming principle} for $v$.

\begin{Theorem}
\label{ThmFeynman-Kac}
Under assumptions {\bf (A$_{\text{\textnormal{\textsc M}}}$)} and {\bf (A${}_{\mathcal R}$)}\textup{-(i)}, for every $(t,x,a)\in[0,T]\times H\times\Lambda$ there exists a unique minimal solution $(Y^{t,x,a},Z^{t,x,a},L^{t,x,a},R^{t,x,a}, K^{t,x,a})\in{\bf S^2(t,T)}\times{\bf L^2(\hat W;t,T)}\times{\bf L^2(\hat\pi;t,T)}\times{\bf L^2(\hat\theta;t,T)}\times{\bf K^2(t,T)}$ to \eqref{BSDEM}-\eqref{ConstraintM}, satisfying
\begin{equation}\label{Y_formula_Markov}
v(s,\hat X_s^{t,x,a}) \ = \ Y_s^{t,x,a}, \qquad \hat\P\text{-a.s.},\,t\leq s\leq T
\end{equation}
and, in particular,
\begin{equation}\label{Y_formula_0_Markov}
v(t,x) \ = \ \hat\E[Y_t^{t,x,a}],
\end{equation}
with $\hat\E[Y_t^{t,x,a}]=Y_t^{t,x,a}$, $\hat\P$-a.s..
\end{Theorem}
\textbf{Proof.}
We firstly define the value function of the so-called randomized stochastic optimal control problem:
\[
\hat v^{\mathcal R}(t,x,a) \ = \ \sup_{\hat\nu\in\hat\Vc_t} \hat\E^{\hat\nu}\bigg[\int_t^T f(s,\hat X_s^{t,x,a},\hat I_s^{t,a})\,ds + g(X_T^{t,x,a})\bigg], \qquad (t,x,a)\in[0,T]\times H\times\Lambda.
\]
Now, we apply Theorems \ref{Main_thm} and \ref{ThmExistence} to our original and randomized control problems. To this end, notice that the control problems in Theorems \ref{Main_thm} and \ref{ThmExistence} are formulated on the time interval $[0,T]$, while our control problems are formulated on the time interval $[t,T]$. Then, taking into account of this time change, we can apply Theorems \ref{Main_thm} and \ref{ThmExistence} interpreting, for what concerns our \emph{original} stochastic control problem, $t$, $x$, $(W_s-W_t)_{s\geq t}$, the restriction of $\pi$ to $[t,\infty)\times U$, $\Ac_t$, $(X_s^{t,x,\alpha})_{s\in[t,T]}$, $v(t,x)$ as follows: $0$, $x_0$, $(W_t)_{t\geq0}$, $\pi$ on $[0,\infty)\times U$, $\Ac$, $(X_s^{x_0,\alpha})_{s\in[t,T]}$, $V_0$ in subsection \ref{SubS:Standard}; similarly, concerning our \emph{randomized} stochastic control problem, we have that $t$, $x$, $a$, $(\hat W_s-\hat W_t)_{s\geq t}$, the restriction of $\hat\pi$ to $[t,\infty)\times U$, the restriction of $\hat\theta$ to $[t,\infty)\times\Lambda$, $\Vc_t$, $(\hat X_s^{t,x,a})_{s\in[t,T]}$, $(\hat I_s^{t,a})_{s\in[t,T]}$ $\hat v^{\mathcal R}(t,x,a)$ correspond to $0$, $x_0$, $a_0$, $(\hat W_t)_{t\geq0}$, $\hat\pi$ on $[0,\infty)\times U$, $\hat\theta$ on $[0,\infty)\times\Lambda$, $\Vc$, $(\hat X_t)_{t\in[0,T]}$, $(\hat I_t)_{t\in[0,T]}$, $\hat V_0^{\mathcal R}$ in Section \ref{S:Randomized}. Then, by Theorem \ref{Main_thm} we deduce that
\[
v(t,x) \ = \ \hat v^{\mathcal R}(t,x,a), \qquad \forall\,(t,x,a)\in[0,T]\times H\times\Lambda.
\]
In addition, by Theorem \ref{ThmExistence} we deduce that there exists a unique minimal solution $(Y^{t,x,a},Z^{t,x,a},$ $L^{t,x,a},R^{t,x,a}, K^{t,x,a})\in{\bf S^2(t,T)}\times{\bf L^2(\hat W;t,T)}\times{\bf L^2(\hat\pi;t,T)}\times{\bf L^2(\hat\theta;t,T)}\times{\bf K^2(t,T)}$ to \eqref{BSDEM}-\eqref{ConstraintM}, satisfying \eqref{Y_formula_0_Markov}, so, in particular,
\[
v(t,x) \ = \ Y_t^{t,x,a}, \qquad \hat\P\text{-a.s.}
\]
for all $(t,x,a)\in[0,T]\times H\times\Lambda$. It remains to prove \eqref{Y_formula_Markov}. To this end, we begin noting that, for every $(t,x,a)\in[0,T]\times H\times\Lambda$, the flow property holds: for every $s\in[t,T]$ we have $(\hat X_r^{s,\hat X_s^{t,x,a},\hat I_s^{t,a}},\hat I_r^{s,\hat I_s^{t,a}}) = (\hat X_r^{t,x,a},\hat I_r^{t,a})$, $\hat\P$-a.s., for any $r\in[s,T]$. Indeed, the flow property for $\hat I^{t,a}$ follows directly from its definition in \eqref{FSDEI}, while the flow property for $\hat X^{t,x,a}$ is a consequence of the uniqueness of the solution to equation \eqref{FSDEX}.
Let us now consider the  penalized backward stochastic differential equation associated with \eqref{BSDEM}-\eqref{ConstraintM}:
\begin{align}
Y_s^n \ &= \ g(\hat X_T^{t,x,a}) + \int_s^T f(r,\hat X_r^{t,x,a},\hat I_r^{t,a}) dr + n\int_s^T \int_\Lambda \big(R_r^n(b)\big)_+ \lambda_\theta(db) dr \label{penBSDEM} \\
&\quad \   - \int_s^T  Z_r^n d\hat W_r - \int_s^T\int_\Lambda R_r^n(b) \hat\theta(dr,db) - \int_s^T\int_{U\setminus \{0\}} L_r^n(z)\,(\hat\pi(dr\,dz)-\lambda_\pi(dz)\,dr), \notag
\end{align}
for all $t\leq s\leq T$, $\hat\P$-a.s..
For every $(t,x,a)\in[0,T]\times H\times\Lambda$, we deduce from Proposition \ref{P:BSDEpen} the existence of a unique solution $(Y^{n,t,x,a},Z^{n,t,x,a},L^{n,t,x,a},R^{n,t,x,a})\in{\bf S^2(t,T)}\times{\bf L^2(\hat W;t,T)}\times{\bf L^2(\hat\pi;t,T)}\times{\bf L^2(\hat\theta;t,T)}$ to \eqref{penBSDEM}. Then, we define the deterministic function $v^n\colon[0,T]\times H\times\Lambda\rightarrow\R$ as (notice that $\hat\E[Y_t^{n,t,x,a}]=Y_t^{n,t,x,a}$, $\hat\P$-a.s., since the random variable $Y_t^{n,t,x,a}$ is $\hat\Fc_t^{t,\hat W,\hat\pi,\hat\theta}$-measurable)
\begin{equation}\label{vn}
\hat v^n(t,x,a) \ := \ \hat\E[Y_t^{n,t,x,a}], \qquad (t,x,a)\in[0,T]\times H\times\Lambda.
\end{equation}
Now, using the flow property and the uniqueness of the solution for the backward stochastic differential equation \eqref{penBSDEM}, we find: for every $s\in[t,T]$, we have $Y_r^{n,s,\hat X_s^{t,x,a},\hat I_s^{t,a}}=Y_r^{n,t,x,a}$, $\hat\P$-a.s., for any $r\in[s,T]$. This implies, from \eqref{vn}, that
\begin{equation}\label{Y_formula_Markov_pen}
\hat v^n(s,\hat X_s^{t,x,a}) \ = \ Y_s^{t,x,a}, \qquad \hat\P\text{-a.s.},\,t\leq s\leq T.
\end{equation}
Finally, by item (i) in Theorem \ref{ThmExistence} we have that $Y_t^{n,t,x,a}$ converges $\hat\P$-a.s. to $Y_t^{t,x,a}$, which implies that $\hat v^n$ converges pointwise to $\hat v^{\mathcal R}$. So, in particular, letting $n\rightarrow\infty$ in equality \eqref{Y_formula_Markov_pen}, we see that \eqref{Y_formula_Markov} holds.
\ep

\begin{Theorem}
\label{ThmDPP}
Let assumptions {\bf (A$_{\text{\textnormal{\textsc M}}}$)} and {\bf (A${}_{\mathcal R}$)}\textup{-(i)} hold.

\vspace{1mm}

\noindent \textup{1)} For every $R>0$, there exists a modulus of continuity $\omega_R$ such that
\[
|v(t,x) - v(t',x)| \ \leq \ \omega_R(|t - t'|),
\]
for all $t,t'\in[0,T]$, $|x|\leq R$.

\vspace{1mm}

\noindent \textup{2)} The randomized dynamic programming principle holds: for every $t\in[0,T]$ and any $\hat\F^{t,\hat W,\hat\pi,\hat\theta}$-stopping time $\hat\tau$ taking values in $[t,T]$, we have
\begin{equation}\label{RandDynProgPr_Markov}
v(t,x) \ = \ \sup_{\hat\nu\in\hat\Vc_t} \hat\E^{t,\hat\nu}
\bigg[\int_t^{\hat\tau} f(s,\hat X_s^{t,x,a},\hat I_s^{t,a})\,ds + v(\hat\tau,\hat X_{\hat\tau}^{t,x,a})\bigg].
\end{equation}
\end{Theorem}
\textbf{Proof.}
We firstly prove a preliminary result, namely the randomized dynamic programming principle for deterministic times: for every $t\in[0,T]$ and any $t'\in[t,T]$,
\begin{equation}\label{RandDynProgPr_Markov_determ}
v(t,x) \ = \ \sup_{\hat\nu\in\hat\Vc_t} \hat\E^{t,\hat\nu}
\bigg[\int_t^{t'} f(s,\hat X_s^{t,x,a},\hat I_s^{t,a})\,ds + v(t',\hat X_{t'}^{t,x,a})\bigg].
\end{equation}
Following the same arguments as in the proof of Theorem \ref{ThmDPP}, we see that we can apply Theorem \ref{ThmExistence} to our backward stochastic differential equation \eqref{BSDEM}-\eqref{ConstraintM}. So, in particular, by \eqref{RandDynProgPr} we have: for every $t\in[0,T]$ and any $\hat\F^{t,\hat W,\hat\pi,\hat\theta}$-stopping time $\hat\tau$ taking values in $[t,T]$,
\begin{equation}\label{Y=Y_tau}
Y_t^{t,x,a} \ = \ \sup_{\hat\nu\in\hat\Vc_t} \hat\E^{t,\hat\nu}
\bigg[\int_t^{\hat\tau} f(s,\hat X_s^{t,x,a},\hat I_s^{t,a})\,ds + Y_{\hat\tau}^{t,x,a}\bigg].
\end{equation}
Now, by \eqref {Y=Y_tau} with $\hat\tau=t'$, together with \eqref{Y_formula_Markov}, we see that \eqref{RandDynProgPr_Markov} follows.

\vspace{1mm}

\noindent\emph{Proof of }1). We proceed as in the proof of Lemma 4.3 in \cite{swiechzabczyk16}. More precisely, fix $R>0$, $0\leq t<t'\leq T$, and $|x|\leq R$. Then, by \eqref{RandDynProgPr_Markov_determ} we have
\begin{equation}\label{vt-t'}
|v(t,x) - v(t',x)| \ \leq \ \sup_{\hat\nu\in\hat\Vc_t} \hat\E^{t,\hat\nu}
\bigg[\int_t^{t'} \big|f(s,\hat X_s^{t,x,a},\hat I_s^{t,a})\big|\,ds + \big|v(t',\hat X_{t'}^{t,x,a}) - v(t',x)\big|\bigg].
\end{equation}
Now, notice that proceeding along the same lines as in the proof of estimate (3.13) of Theorem 3.4 in \cite{swiechzabczyk16}, we can prove that the following estimate holds:
\begin{equation}\label{Estimatet-t'}
\hat\E\Big[\sup_{t\leq s\leq t'}\big|\hat X_s^{t,x,a} - x\big|^2\big] \ \leq \ \omega_x(t' - t),
\end{equation}
for some modulus $\omega_x$. Then, using the assumptions on $f$ in {\bf (A$_{\text{\textnormal{\textsc M}}}$)}-(v), estimates \eqref{Est_X_rand_Markov} and \eqref{Estimatet-t'}, inequality \eqref{v_unif}, and estimate (D.1) in \cite{fabbrigozziswiech17}, we obtain from \eqref{vt-t'}:
\[
|v(t,x) - v(t',x)| \ \leq \ \tilde C\,(t' - t)\,(1 + |x|) + \sup_{\hat\nu\in\hat\Vc_t} \hat\E^{t,\hat\nu}\Big[\omega_v\big(\big|\hat X_{t'}^{t,x,a} - x\big|_{-1}\big)\Big] \ \leq \ \omega_R(|t - t'|),
\]
for some constant $\tilde C\geq0$ and some modulus $\omega_R$.

\vspace{1mm}

\noindent\emph{Proof of }2). From item 1) and inequality \eqref{v_unif}, it follows that $v$ is continuous on $[0,T]\times H$ (taking on $H$ the usual norm $|\cdot|$). As a consequence, the stochastic process $(v(s,\hat X_s^{t,x,a}))_{s\in[t,T]}$ has c\`{a}dl\`{a}g paths. Since $(Y_s^{t,x,a})_{s\in[t,T]}$ also has c\`{a}dl\`{a}g paths, we see that the two stochastic processes $(v(s,\hat X_s^{t,x,a}))_{s\in[t,T]}$ and $(Y_s^{t,x,a})_{s\in[t,T]}$ are $\hat\P$-indistinguishable, since by \eqref{Y_formula_Markov} are one the modification of the other. In other words, it holds that
\begin{equation}\label{Y_formula_Markov_bis}
v(s,\hat X_s^{t,x,a}) \ = \ Y_s^{t,x,a}, \qquad t\leq s\leq T,\,\hat\P\text{-a.s.}
\end{equation}
In particular, given any $\hat\F^{t,\hat W,\hat\pi,\hat\theta}$-stopping time $\hat\tau$ taking values in $[t,T]$, we deduce from \eqref{Y_formula_Markov_bis} that
\[
v(\hat\tau,\hat X_{\hat\tau}^{t,x,a}) \ = \ Y_{\hat\tau}^{t,x,a}, \qquad \hat\P\text{-a.s.}
\]
Then, by \eqref{Y=Y_tau} we see that \eqref{RandDynProgPr_Markov} holds.
\ep

\subsection{Viscosity property of the value function $v$}

We now exploit the randomized dynamic programming principle \eqref{RandDynProgPr_Markov} in order to prove that the value function $v$ in \eqref{value} is a viscosity solution to the following Hamilton-Jacobi-Bellman equation:
\begin{equation}
\!\!\!\!\!\!\!\!\!\!\begin{cases}
v_t + \langle Ax,D_xv\rangle + \sup_{a\in\Lambda} \Big\{\frac{1}{2}\text{Tr}\big(\sigma(t,x,a)\sigma^*(t,x,a) D_x^2v\big) +\langle b(t,x,a),D_xv\rangle +f(t,x,a) \\
+\!\int_{U \setminus \{0\}} (v(t,x+ \gamma(t,x,a,z))-v(t,x)-D_x v(t,x) \gamma(t,x,a,z))\lambda_\pi(dz)\!\Big\} = 0,  \;\text{on }(0,T)\!\times\!H,   \\
v(T,x) = g(x), \hspace{10cm} x\in H.\label{E:HJB}
\end{cases}
\end{equation}

We adopt the definition of  viscosity solution given in \cite{swiechzabczyk16}, Definition 5.2, which requires the following notions.

\begin{Definition}
Let $u\colon(0,T)\times H\rightarrow\R$.

\vspace{1mm}

\noindent We say that $u$ is $B$\textbf{-upper semicontinuous} if, for all $(t,x)\in (0,T)\times H$,
\[
\limsup_{\substack{m\rightarrow+\infty\\(t_m,x_m)\in (0,T)\times H}} \hspace{-.5cm} u(t_m,x_m) \ \leq \ u(t,x)
\]
whenever $t_m\rightarrow t$, $x_m \rightharpoonup x $, $B x_m \rightarrow B x$.

\vspace{1mm}

\noindent We say that $u$ is $B$\textbf{-lower semicontinuous} if, for all $(t,x)\in (0,T)\times H$,
\[
\liminf_{\substack{m\rightarrow+\infty\\(t_m,x_m)\in (0,T)\times H}} \hspace{-.5cm} u(t_m,x_m) \ \geq \ u(t,x)
\]
whenever $t_m\rightarrow t$, $x_m \rightharpoonup x $, $B x_m \rightarrow B x$.

\vspace{1mm}

\noindent We say that $u$ is $B$\textbf{-continuous} if it is both $B$-upper semicontinuous and $B$-lower semicontinuous.
\end{Definition}

\begin{Definition}
\label{D:Test}
A function $\psi\colon(0,T)\times H\rightarrow\R$ is a \textbf{test function} if $\psi(t,x) = \varphi(t,x) +\delta(t,x)h(|x|)$, where:
\begin{enumerate}
\item[\textup{(i)}] $\varphi_t$, $D_x\varphi$, $D_x^2\varphi$, $A^*D_x\varphi$, $\delta_t$, $D_x\delta$, $D_x^2\delta$, $A^*D_x\delta$ are uniformly continuous on $(\varepsilon,T-\varepsilon)\times H$, for every $\eps>0$; in addition, $\varphi$ is $B$-lower semicontinuous; finally, $\delta \geq 0$, bounded, and $B$-continuous.
\item[\textup{(ii)}] $h$ is even, $h'$ and $h''$ are uniformly continuous on $\R$, $h'(r)\geq0$ for every $r>0$.
\end{enumerate}
\end{Definition}

\begin{Remark}
{\rm
Notice that a test function $\psi$ satisfies the following property: for every $\eps>0$, there exists a constant $C_\eps\geq0$ such that $|\psi(t,x)|\leq C_\eps(1+|x|^2)$ on $(\eps,T-\eps)\times H$.
\epR}
\end{Remark}

\begin{Definition}
\label{D:Visc}
\textup{(i)} A $B$-upper semicontinuous function $u\colon(0,T)\times H\rightarrow\R$ is a \textbf{viscosity supersolution} of \eqref{E:HJB} if whenever
\[
(u-\psi)(t,x) = \min_{(0,T)\times H}(u-\psi)
\]
for $(t,x)\in(0,T)\times H$ and $\psi(s,y)=\varphi(s,y)+\delta(s,y)h(|y|)$ a test function, then
\begin{align*}
&\psi_t(t,x)
- \langle x,A^*D_x\varphi(t,x) + h(|x|) A^*D_x\delta(t,x)\rangle \\
&+ \sup_{a\in\Lambda} \bigg(\frac{1}{2}\textup{Tr}\big(\sigma(t,x,a)\sigma^*(t,x,a) D_x^2\psi(t,x)\big)+ \langle b(t,x,a),D_x\psi(t,x)\rangle + f(t,x,a) \\
& +\int_{U \setminus \{0\}} (\psi(t,x)(t,x+ \gamma(t,x,a,z))-\psi(t,x)(t,x)-D_x \psi(t,x)(t,x) \gamma(t,x,a,z))\lambda_\pi(dz)\bigg)\leq 0.
\end{align*}
\textup{(ii)} A $B$-lower semicontinuous function $u\colon(0,T)\times H\rightarrow\R$ is a \textbf{viscosity subsolution} of \eqref{E:HJB} if whenever
\[
(u+\psi)(t,x) = \max_{(0,T)\times H}(u+\psi)
\]
for $(t,x)\in(0,T)\times H$ and $\psi(s,y)=\varphi(s,y)+\delta(s,y)h(|y|)$ a test function,
then
\begin{align*}
&-\psi_t(t,x)
+ \langle x,A^*D_x\varphi(t,x) + h(|x|)A^*D_x\delta(t,x) \rangle \\
&+ \sup_{a\in\Lambda} \bigg(-\frac{1}{2}\textup{Tr}\big(\sigma(t,x,a)\sigma^*(x,a) D_x^2\psi(t,x)\big)-\langle b(t,x,a),D_x\psi(t,x)\rangle + f(t,x,a)\\
& -\int_{U \setminus \{0\}} (\psi(t,x)(t,x+ \gamma(t,x,a,z))-\psi(t,x)(t,x)-D_x \psi(t,x)(t,x) \gamma(t,x,a,z))\lambda_\pi(dz)\bigg) \geq 0.
\end{align*}
\textup{(iii)} A function $u\colon(0,T)\times H\rightarrow\R$ is a \textbf{viscosity solution} of \eqref{E:HJB} if it is both a viscosity subsolution and a  viscosity supersolution of \eqref{E:HJB}.
\end{Definition}

In order to prove that  $v$ is a viscosity solution to equation \eqref{E:HJB} we will need  the following technical result.

\begin{Lemma}
\label{L:Dynkin}
Let assumption {\bf (A$_{\text{\textnormal{\textsc M}}}$)} hold. Let $\psi=\varphi+\delta h(|\cdot|)$ be a test function. Fix $t,t'\in(0,T)$, with $t<t'$, and let $\hat\tau$ be a $\hat\F^{t,\hat W,\hat\pi,\hat\theta}$-stopping time taking values in $[t,t']$. Then, for any $(x,a)\in H \times \Lambda$, $\hat \nu \in \hat {\mathcal V}_t$,
\begin{align}
\label{E:Ito}
& \hat\E^{t,\hat\nu}\left[\psi(\hat\tau,\hat X_{\hat\tau}^{t,x,a})\right] \ \geq \ \psi(t,x) + \hat\E^{t,\hat\nu}\Big[ \int_t^{\hat\tau} \psi_t(r,\hat X_r^{t,x,a})dr \notag\\
&-\int_t^{\hat\tau} \langle \hat X_r^{t,x,a},A^*D_x\psi(r,\hat X_r^{t,x,a})+ h(|\hat X_r^{t,x,a}|) A^* D_x \delta(r,\hat X_r^{t,x,a})\rangle dr\notag\\
&  + \frac{1}{2}\int_t^{\hat\tau} \textup{Tr}\big[\sigma(r,\hat X_r^{t,x,a},\hat I_r^{t,a})\sigma^*(r,\hat X_r^{t,x,a},\hat I_r^{t,a})D_x^2 \psi(r,\hat X_r^{t,x,a})\big] dr \notag \\
& + \int_t^{\hat\tau} \langle b(r,\hat X_r^{t,x,a},\hat I_r^{t,a}),D_x \psi(r,\hat X_r^{t,x,a})\rangle dr +\int_t^{\hat\tau} \int_{U \setminus \{0\}} \Big(\psi\big(r,\hat X_r^{t,x,a}+ \gamma(r,\hat X_r^{t,x,a},\hat I_r^{t,a}, z)\big)\notag\\
& \qquad - \psi(r,\hat X_r^{t,x,a})- D_x\psi(r,\hat X_r^{t,x,a}) \gamma(r,\hat X_r^{t,x,a},\hat I_r^{t,a}, z)\Big) \lambda_{\pi}(dz)dr\Big].
\end{align}
\end{Lemma}
\textbf{Proof.}
The proof can be done proceeding along the same lines as in the proof of Lemma 5.3  in \cite{swiechzabczyk16}, the only difference being the presence of the pure jump process $\hat I^{t,a}$. For this reason, here we just give an outline. The proof consists in approximating the process $\hat X^{t,x,a}$ by means of a sequence of more regular processes $\hat X^{n,t,x,a}$, which are obtained replacing the operator $A$ in equation \eqref{FSDEX} by its Yosida approximations $(A_n)_n$. It is well-known, see e.g. Theorem 27.2 in \cite{metivier82}, that $\psi(\cdot,\hat X_{\cdot}^{n,t,x,a})$ satisfies an It\^o formula. Then, using convergence results of $\hat X^{n,t,x,a}$ towards $\hat X^{t,x,a}$, which can be found for instance in Proposition 1.115 of \cite{fabbrigozziswiech17}, and taking the expectation under $\hat\P^{t,\hat\nu}$,  we deduce  \eqref{E:Ito} using that $\langle -A \hat X_r^{t,x,a},\delta(r,\hat X_r^{t,x,a})\frac{h'(|\hat X_r^{t,x,a}|)}{|\hat X_r^{t,x,a}|} \hat X_r^{t,x,a}\rangle \geq 0$.
\ep

\begin{Proposition}
\label{P:ViscPropvm}
Let assumptions {\bf (A$_{\text{\textnormal{\textsc M}}}$)} and {\bf (A${}_{\mathcal R}$)}\textup{-(i)} hold. The value function $v$ defined in \eqref{value} is a viscosity solution to equation \eqref{E:HJB}.
\end{Proposition}
\textbf{Proof.} We split the proof into two steps.

\noindent\emph{Proof of the viscosity subsolution property of $v$.} Let $(t,x,a)\in(0,T)\times H\times\Lambda$ and let $\psi(s,y) = \varphi(s,y) + \delta(s,y)h(|y|)$ be a test function such that
$(v+\psi)(t,x) = \max_{(0,T)\times H}(v+\psi)$.
We shall prove that
\begin{align*}
&-\psi_t(t,x) + \langle x,A^*D_x\psi(t,x)\rangle + h(|x|) A^\ast D \delta(t,x) \\
&+ \sup_{a \in \Lambda} \Big\{-\frac{1}{2}{\rm Tr}\big(\sigma(t,x,a)\sigma^*(t,x,a) D_x^2 \psi(t,x)\big)  -\langle b(t,x,a),D_x \psi(t,x)\rangle + f(t,x,a) \\
&-\int_{U \setminus \{0\}} \big(\psi(t,x + \gamma(t,x,a,z))-\psi(t,x)-D_x\psi(t,x)\gamma(t,x,a,z)\big) \lambda_{\pi}(dz)\Big\} \ \geq \ 0.
\end{align*}
We assume, without loss of generality, that
\begin{equation}\label{1}
v(t,x)+\psi(t,x) \ = \ 0,
\end{equation}
so, in particular,
\begin{equation}\label{2}
v(s,y) + \psi(s,y) \ \leq \ 0, \qquad \forall\,(s,y) \in (0,\,T) \times H.
\end{equation}
For any $\eta >0$, we define $\beta(\eta):=\sup_{(s,y)\in \partial B(t,x; \eta)}(v + \psi)(s,y)$, where
\begin{align*}
B(t,x; \eta) \ &= \ \big\{(s,y) \in (0,T) \times H: \max \{|x-y|, |t-s|\}< \eta\big\}, \\
\partial B(t,x; \eta) \ &= \ \big\{(s,y) \in (0,T) \times H: \max \{|x-y|, |t-s|\}= \eta\big\}.
\end{align*}
Notice that $\beta(\eta)<0$, for any $\eta>0$. Let us proceed by contradiction, assuming that
\begin{align*}
&-\psi_t(t,x) + \langle x,A^*D_x\psi(t,x)\rangle + h(|x|) A^\ast D \delta(t,x) \\
&+ \sup_{a \in \Lambda} \Big\{-\frac{1}{2}{\rm Tr}\big(\sigma(t,x,a)\sigma^*(t,x,a) D_x^2 \psi(t,x)\big)  -\langle b(t,x,a),D_x \psi(t,x)\rangle + f(t,x,a) \\
&-\int_{U \setminus \{0\}} \big(\psi(t,x + \gamma(t,x,a,z))-\psi(t,x)-D_x\psi(t,x)\gamma(t,x,a,z)\big) \lambda_{\pi}(dz)\Big\} \ < \ 0.
\end{align*}
Using the Lipschitz property of $b$, $\sigma$, $\gamma$, and the uniform continuity of $f$, when on $H$ we consider the standard topology induced by the norm $|\cdot|$ (notice that $b$, $\sigma$, $f$ satisfy the mentioned properties when on $(H,|\cdot|_{-1})$, and hence they satisfy the same properties on $(H,|\cdot|)$), and using also the uniform continuity of $\psi_t$, $A^*D_x\psi$, $D_x\psi$, and $D_x^2\psi$, we have that, given $\eta \in (0, 2 (T-t))$, there exists $\varepsilon \in (0,\,-\beta(\eta)/(T-t)]$, with $\eps<T$, such that
\begin{align}
&-\psi_t(s,y) + \langle y,A^*D_x\psi(s,y)\rangle + h(|y|) A^\ast D \delta(s,y) \notag \\
&+ \sup_{a \in \Lambda} \Big\{-\frac{1}{2}{\rm Tr}\big(\sigma(s,y,a)\sigma^*(s,y,a) D_x^2 \psi(s,y)\big)  -\langle b(s,y,a),D_x \psi(s,y)\rangle + f(s,y,a) \label{E:eps} \\
&-\int_{U \setminus \{0\}} \big(\psi(s,y + \gamma(s,y,a,z))-\psi(s,y)-D_x\psi(s,y)\gamma(s,y,a,z)\big) \lambda_{\pi}(dz)\Big\} \ \leq \ -\eps, \notag
\end{align}
for any $(s,y)\in(0,T)\times H$ with $|s-t|,|y-x|\leq\eta$. Define
\[
\hat\tau \ := \ \inf \big\{s\in[t,T]:
(s,\hat X_s^{t,x,a}) \notin B(t,x; \eta/2)
\big\}, \qquad \hat\theta \ := \ \hat\tau
\wedge T,
\]
where $\inf\emptyset=\infty$. Since the stochastic process $(\hat X_s^{t,x,a})_{s\in[t,T]}$ is c\`{a}dl\`{a}g, it is in particular right-continuous at time $t$. As a consequence, $\hat\theta>t$, $\hat\P$-a.s..

For every $\eps>0$, by the randomized dynamic programming principle \eqref{RandDynProgPr_Markov}, it follows that there exists $\hat\nu^\eps\in\hat\Vc_t$ such that
\[
v(t,x) \ \leq \ \hat\E^{t,\hat\nu^\eps}\bigg[\int_t^{\hat\theta} f(r,\hat X_r^{t,x,a},\hat I^{t,a}_{r})\, dr + v(\hat\theta,\hat X^{t,x,a}_{\hat\theta}) \bigg]+\frac{\varepsilon}{2}(T-t),
\]
which in turn yields, by \eqref{1}-\eqref{2},
\[
-\psi(t,x) \ \leq \ \hat\E^{t,\hat\nu^\eps}\bigg[\int_t^{\hat\theta} f(r,\hat X_r^{t,x,a},\hat I^{t,a}_{r})\, dr  -\psi(\hat\theta,\hat X^{t,x,a}_{\hat\theta}) + \beta(\delta)\,1_{\{\hat\tau \leq T\}}\bigg] + \frac{\varepsilon}{2}(T-t).
\]	
By applying Lemma \ref{L:Dynkin}, the previous inequality yields
 \begin{align*}
-\frac{\varepsilon}{2}(T-t) \ &\leq \ \hat\E^{t,\hat\nu^\eps}\bigg[\int_t^{\hat\theta} \langle\hat X_r^{t,x,a},A^*D_x\psi(r,\hat X_r^{t,x,a})+ h(|\hat X_r^{t,x,a}|) A^* D_x \delta(r,\hat X_r^{t,x,a})\rangle dr\bigg] \\
&\quad \ + \hat\E^{t,\hat\nu^\eps}\bigg[ \int_t^{\hat\theta} \bigg(-\psi_t(r,\hat X_r^{t,x,a}) - \langle b(r,\hat X_r^{t,x,a},\hat I_r^{t,a}),D_x \psi(r,\hat X_r^{t,x,a})\rangle  \\
&\quad \ - \frac{1}{2}\textup{Tr}\big[\sigma(r,\hat X_r^{t,x,a},\hat I_r^{t,a})\sigma^*(r,\hat X_r^{t,x,a},\hat I_r^{t,a})D_x^2 \psi(r,\hat X_r^{t,x,a})\big] +f(r,\hat X_r^{t,x,a},\hat I^{t,a}_{r}) \notag \\
&\quad \ + \beta(\delta)\,\hat\P^{t,\hat\nu^\eps}(\hat\tau \leq T) - \int_{U \setminus \{0\}} \Big(\psi(r,\hat X_r^{t,x,a}+ \gamma(r,\hat X_r^{t,x,a},\hat I_r^{t,a}, z)\notag\\
& \quad \ -\psi(r,\hat X_r^{t,x,a})- D_x\psi(r,\hat X_r^{t,x,a}) \gamma(r,\hat X_r^{t,x,a},\hat I_r^{t,a}, z)\Big) \lambda_{\pi}(dz)\bigg)dr\bigg]\\
&\leq \ - \varepsilon\,(T-t)\,\hat\P^{t,\hat\nu^\eps}(\hat\tau \leq T) -\varepsilon\,\hat\E^{t,\hat\nu^\eps}[\hat\theta -t] \ \leq \ -\varepsilon\,(T-t),
\end{align*}
where we have  used \eqref{E:eps} and the fact that $\hat\tau \leq \frac{\eta}{2} \leq T$. This yields a contradiction and concludes the proof.

\vspace{2mm}

\noindent \emph{Proof of the viscosity supersolution property of $v$.} Let $(t,x,a)\in(0,T)\times H\times\Lambda$ and let $\psi(s,y) = \varphi(s,y) + \delta(s,y)h(|y|)$ be a test function such that $(v-\psi)(t,x) = \min_{(0,T)\times H}(v-\psi)$.
We shall prove that
\begin{align*}
&\psi_t(t,x)
- \langle x,A^*D_x\varphi(t,x) + h(|x|) A^*D_x\delta(t,x)\rangle \\
&+ \sup_{a\in\Lambda} \bigg(\frac{1}{2}\textup{Tr}\big(\sigma(t,x,a)\sigma^*(t,x,a) D_x^2\psi(t,x)\big)+ \langle b(t,x,a),D_x\psi(t,x)\rangle + f(t,x,a) \\
& +\int_{U \setminus \{0\}} (\psi(t,x)(t,x+ \gamma(t,x,a,z))-\psi(t,x)(t,x)-D_x \psi(t,x)(t,x) \gamma(t,x,a,z))\lambda_\pi(dz)\bigg) \ \leq \ 0.
\end{align*}
We assume that
\begin{equation}\label{1bis}
v(t,x) - \psi(t,x) \ = \ 0,
\end{equation}
so, in particular,
\begin{equation}\label{2bis}
v(s,y) - \psi(s,y) \ \geq \ 0, \qquad \forall\,(s,y) \in (0,\,T) \times H.
\end{equation}
Let $h >0$, $\eta>0$, and set
$$
\hat\tau \ := \ \inf\big\{s\in[t,T]\colon |\hat X_{s}^{t,x,a} - x| > \eta\big\}, \qquad \hat\theta \ := \ \hat\tau \wedge (t+h) \wedge \hat T_1,
$$
where we recall that $(\hat T_n,\hat\eta_n)_{n\geq1}$ is the marked point process associated with the random measure $\hat \theta$ (in particular we have $\hat \theta(dt\,da)=\sum_{n \geq 1} \delta_{(\hat T_n,\hat\eta_n)}(dt\,da)$). So, in particular, $\hat T_1$ is the first jump time of the stochastic process $\hat I^{t,a}$ defined in \eqref{FSDEI}.

By the randomized dynamic programming principle \eqref{RandDynProgPr_Markov}, we have
\begin{align*}
v(t,x) \ &\geq \ \hat\E^{t,\hat\nu}\bigg[\int_t^{\hat\theta} f(r,\hat X_r^{t,x,a},\hat I^{t,a}_{r})\, dr + v(\hat\theta, \hat X^{t,x,a}_{\hat\theta}) \bigg],\qquad \forall\,\hat\nu \in \hat{\mathcal V}_t,
\end{align*}
which in turn yields, by \eqref{1bis}-\eqref{2bis},
\begin{align*}
\psi(t,x) \ &\geq \ \hat\E^{t,\hat\nu}\bigg[\int_t^{\hat\theta} f(r,\hat X_r^{t,x,a},\hat I^{t,a}_{r})\, dr + \psi(\hat\theta, \hat X^{t,x,a}_{\hat\theta})\bigg],\qquad \forall\,\hat\nu \in \hat{\mathcal V}_t.
\end{align*}
We take $\hat \nu =1$, so that in the above inequality $\hat\E^{t,\hat\nu}$ coincides with the expectation $\hat\E$ under $\hat\P$. Applying Lemma \ref{L:Dynkin}, we obtain
\begin{align}\label{final_sub}
&0 \geq \hat\E\bigg[ \frac{1}{h}\int_t^{\hat\theta} \psi_t(r,\hat X_r^{t,x,a})dr - \frac{1}{h}\int_t^{\hat\theta}\langle \hat X_r^{t,x,a},A^*D_x\psi(r,\hat X_r^{t,x,a})+ h(|\hat X_r^{t,x,a}|) A^* D_x \delta(r,\hat X_r^{t,x,a})\rangle dr \notag\\
&+ \frac{1}{h}\int_t^{\hat\theta} f(r,\hat X_r^{t,x,a},\hat I^{t,a}_{r})\, dr + \frac{1}{2}\int_t^{\hat\theta} \textup{Tr}\big[\sigma(r,\hat X_r^{t,x,a},\hat I_r^{t,a})\sigma^*(r,\hat X_r^{t,x,a},\hat I_r^{t,a})D_x^2 \psi(r,\hat X_r^{t,x,a})\big] dr \notag \\
&+ \frac{1}{h} \int_t^{\hat\theta} \langle b(r,\hat X_r^{t,x,a},\hat I_r^{t,a}),D_x \psi(r,\hat X_r^{t,x,a})\rangle dr + \frac{1}{h}\int_t^{\hat\theta} \int_{U \setminus \{0\}} \Big(\psi(r,\hat X_r^{t,x,a}+ \gamma(r,\hat X_r^{t,x,a},\hat I_r^{t,a}, z))\notag\\
& - \psi(r,\hat X_r^{t,x,a})- D_x\psi(r,\hat X_r^{t,x,a}) \gamma(r,\hat X_r^{t,x,a},\hat I_r^{t,a}, z)\Big) \lambda_{\pi}(dz)dr\bigg].
\end{align}
Now we notice that, $\hat \P$-a.s., $\hat I_r^{t,a} = a$ and $\hat X^{t,x,a}$ is right-continuous at $t$ (indeed, it is a c\`adl\`ag process). Thus, by the mean value theorem, the random variable inside the expectation $\hat \E$ in \eqref{final_sub} converges $\hat \P$-a.s. to
\begin{align*}
&\psi_t(t,x)-\langle x,\,A^\ast D_x\psi(t,x) + h(|x|)\,A^\ast D_x\delta(t,x)\rangle \\
&+ \langle b(t,x,a),\,D_x\psi(t,x)\rangle + \frac{1}{2}\textup{Tr}\big[\sigma(t,x,a)\sigma^*(t,x,a)D_x^2 \psi(t,x)\big] +f(t,x,a)  \\
&+ \int_{U \setminus \{0\}}\big(\psi(t,x + \gamma(t,x,a,z))-\psi(t,x)+ \gamma(t,x,a,z)\big)D_x\psi(t,x)\,\lambda_\pi(dz)
\end{align*}
when $h$ goes to zero. Then, by the Lebesgue dominated convergence theorem, we obtain from \eqref{final_sub}
\begin{align*}
&\psi_t(t,x)-\langle x,\,A^\ast D_x\psi(t,x) + h(|x|)\,A^\ast D_x\delta(t,x)\rangle \\
&+ \langle b(t,x,a),\,D_x\psi(t,x)\rangle + \frac{1}{2}\textup{Tr}\big[\sigma(t,x,a)\sigma^*(t,x,a)D_x^2 \psi(t,x)\big] + f(t,x,a) \\
&+ \int_{U \setminus \{0\}}\big(\psi(t,x + \gamma(t,x,a,z))-\psi(t,x)+ \gamma(t,x,a,z)\big)D_x\psi(t,x)\,\lambda_\pi(dz) \ \leq \ 0.
\end{align*}
The claim follows from the arbitrariness of $a\in\Lambda$.
\ep

\begin{Remark}\label{R:Uniq}
{\rm
Concerning the uniqueness of viscosity solutions to the Hamilton-Jacobi-Bellman equation \eqref{E:HJB}, a positive result follows from the comparison principle in \cite{swiechzabczyk13}, Theorem 6.2, under the additional assumptions that $f$ and $g$ are bounded and $\Lambda$ is compact, from which we deduce that the value function $v$ in \eqref{value} is the unique viscosity solution in the class of bounded and uniformly continuous solutions on $[0,T]\times H_{-1}$.
\epR}
\end{Remark}

\vspace{9mm}

\small
\bibliographystyle{plain}
\bibliography{biblio}

\end{document}